# Calculating Covering Constants for Mappings in Euclidean Spaces Using Mordukhovich Coderivatives with Applications

Jinlu Li

Department of Mathematics
Shawnee State University
Portsmouth, OH 45662 USA
Email: jli@shawnee.edu

**Abstract** In this paper, we calculate the covering constants for single-valued mappings in Euclidean space by using Mordukhovich derivatives (or coderivatives). At first, we prove the guideline for calculating the Fréchet derivatives of single-valued mappings by their partial derivatives. Then, by using the connections between Fréchet derivatives and Mordukhovich derivatives (or coderivatives) of single-valued mappings in Banach spaces, we derive the useful rules for calculating the Mordukhovich derivatives of single-valued mappings in Euclidean spaces. For practicing these rules, we find the precise solutions of the Fréchet derivatives and Mordukhovich derivatives for some single-valued mappings between Euclidean spaces. By using these solutions, we find or estimate the covering constants for the considered mappings. As applications of the results about the covering constants involved in the Arutyunov Mordukhovich and Zhukovskiy Parameterized Coincidence Point Theorem, we solve some parameterized equations.

**Mathematics subject classification** (**2010**)   49J52, 49J53, 47H10, 90C31

## 1. Introduction

The concept of Fréchet derivative of a single-valued mapping is a generalization of the derivative of a real-valued function of a single real variable in Euclidean spaces in calculus to the case of a vector-valued function of multiple real variables in normed spaces, in particular, in Banach spaces. The Fréchet derivative has been widely used in the calculus of variations and nonlinear functional analysis in Banach spaces (See [8, 9, 11, 22]).

One more step further, in the theory of generalized differentiational analysis in Banach spaces, the concept of Mordukhovich derivative (or coderivative) for set-valued mappings lays the foundation and plays fundamental and crucial role in the theory of set-valued and variational analysis (See [18−21, 23]). However, even for single-valued mappings, it is very difficult to precisely find the Mordukhovich derivatives of the given mappings in Banach spaces, in general. Fortunately, the results of Theorem 1.38 in [18] provides close connections between the Fréchet derivative and the Mordukhovich derivative (or coderivative) of single-valued mappings in Banach spaces. By this theorem, the Mordukhovich derivative can be considered to be the generalization of the Fréchet derivative of single-valued mappings.

The theory of generalized differentiation in set-valued analysis has been rapidly developed during the last years and has been widely applied to many fields such as optimization, control theory, game theory, variational analysis, and so forth (see [1−7, 12−21, 24]). It is well-known that, in contrast to the differentiation theory in calculus when the generalized differentiation is applied to a specific problem with a precisely given mapping, one needs to find the exact solution of the Mordukhovich derivative of the considered mapping. For example, the concept of covering constant plays an important and crucial role in Arutyunov Mordukhovich and Zhukovskiy Parameterized Coincidence Point Theorem. However, the covering constants are defined by the Mordukhovich derivatives. So, if one wants to find the covering constant for a mapping, you have firstly to find the Mordukhovich derivatives of the considered mapping. This idea leads an important motivation for authors to find some formulas for both the Fréchet derivatives and Mordukhovich derivatives of mappings in Banach spaces. For example, in [12–13], the precise solutions for the standard metric project operator have been proved in both Hilbert spaces and uniformly convex and uniformly smooth Banach spaces.

In this paper, we concentrate to study the covering constants for single-valued mappings in Euclidean spaces by using the Fréchet derivatives and Mordukhovich derivatives, which is implemented as the following steps.

(i) At first, we prove the principle of the Fréchet differentiability of single-valued mappings in Euclidean spaces.
(ii) By the second order approximations of mappings, we derive the guidelines for calculating the Fréchet derivatives.
(iii) We use the relationship between Fréchet derivatives and Mordukhovich derivatives (It is given by Theorem 1.38 in [18]) to deduce the algorithm for calculating the Mordukhovich derivatives of single-valued mappings in Euclidean spaces.
(iv) We provide examples of mappings between Euclidean spaces, which includes polynomial mappings, rational mappings, exponential mappings and logarithm mappings. For these mappings, we find their Fréchet derivatives, Mordukhovich derivatives and covering constants.
(v) By the results of covering constants obtained in Step iv, we apply the Arutyunov Mordukhovich and Zhukovskiy Parameterized Coincidence Point Theorem to solve systems of parameterized equations.

**2. Review for Fréchet Derivative and Mordukhovich Derivative (Coderivative)**

In this section, we briefly review the concepts and properties of Fréchet derivatives and Mordukhovich derivatives (or coderivatives) of single-valued mappings in Banach spaces, in particular, in Euclidean spaces, which will be used in this paper. See [18–21] for more details,

Let $(X, \|\cdot\|_X)$ and $(Y, \|\cdot\|_Y)$ be real Banach spaces, which have topological dual spaces $(X^*, \|\cdot\|_{X^*})$ and $(Y^*, \|\cdot\|_{Y^*})$, respectively. Let $\langle \cdot, \cdot \rangle_X$ denote the real canonical pairing between $X^*$ and $X$ and $\langle \cdot, \cdot \rangle_Y$ the real canonical pairing between $Y^*$ and $Y$. Let $\theta_X$ and $\theta_Y$ denote the origins in $X$ and $Y$, respectively. For any $x \in X$ and $r > 0$, let $B_X(x, r)$ and $S_X(x, r)$ denote the closed ball and sphere in $X$ centered at point $x$ with radius $r$, respectively.

Let $g: X \to Y$ be a single-valued mapping and let $\bar{x} \in X$. If there is a linear and continuous mapping $\nabla g(\bar{x}): X \to Y$ such that

$$\lim_{x \to \bar{x}} \frac{g(x) - g(\bar{x}) - \nabla g(\bar{x})(x - \bar{x})}{\|x - \bar{x}\|_X} = \theta_Y,$$

then $g$ is said to be Fréchet differentiable at $\bar{x}$ and $\nabla g(\bar{x})$ is called the Fréchet derivative of $g$ at $\bar{x}$.

The Mordukhovich derivative (coderivative) is initially defined for set-valued mappings in Banach spaces. Since in this paper, we only deal with the single-valued mappings in Euclidean spaces, so we only review the concept of Mordukhovich derivative (coderivative) for single-valued mappings in Banach spaces (See [18–21] for more details). Let $\Delta$ be a nonempty subset of $X$ and let $g: \Delta \to Y$ be a single-valued mapping. For $x \in \Delta$ and $y = g(x)$, a set-valued mapping $\widehat{D}^*g(x,y): Y^* \to X^*$ is defined, for any $y^* \in Y^*$, by (See Definitions 1.13 and 1.32 in Chapter 1 in [15])

$$\widehat{D}^*g(x,y)(y^*) = \left\{ z^* \in X^*: \limsup_{\substack{(u,g(u)) \to (x,g(x)) \\ u \in \Delta}} \frac{\langle z^*, u-x \rangle_X - \langle y^*, g(u)-g(x) \rangle_Y}{\|u-x\|_X + \|g(u)-g(x)\|_Y} \leq 0 \right\}.$$

If this set-valued mapping $\widehat{D}^*g(x,y): Y^* \to X^*$ satisfies that

$$\widehat{D}^*g(x,y)(y^*) \neq \emptyset, \text{ for every } y^* \in Y^*, \tag{2.1}$$

then $g$ is said to be Mordukhovich differentiable (codifferentiable) at the point $x$ and $\widehat{D}^*g(x,y)$ is called the Mordukhovich derivative (which is also called Mordukhovich coderivative, or coderivative) of $g$ at point $x$. For this single-valued mapping $g: \Delta \to Y$, at point $x$, we write

$$\widehat{D}^*g(x,y)(y^*) \equiv \widehat{D}^*g(x)(y^*), \text{ for every } y^* \in Y^*.$$

Furthermore, if $g: \Delta \to Y$ is a single-valued continuous mapping. Then, by the above definition, the Mordukhovich derivative of $g$ at point $x$ is calculated by

$$\widehat{D}^*g(x)(y^*) = \left\{ z^* \in X^*: \limsup_{\substack{u \to x \\ u \in \Delta}} \frac{\langle z^*, u-x \rangle_X - \langle y^*, g(u)-g(x) \rangle_Y}{\|u-x\|_X + \|g(u)-g(x)\|_Y} \leq 0 \right\}, \text{ for any } y^* \in Y^*. \tag{2.2}$$

The following theorem shows the connection between Fréchet derivatives and Mordukhovich derivatives for sing-valued mappings. The results of the following theorem provide a powerful tool to calculate the Mordukhovich derivatives by the Fréchet derivatives of single-valued mappings.

**Theorem 1.38 in [15]**. *Let $X$ be a Banach space with dual space $X^*$ and let $g: X \to Y$ be a single-valued mapping. Suppose that $g$ is Fréchet differentiable at $x \in X$ with $y = g(x)$. Then, the Mordukhovich derivative of $g$ at $x$ satisfies the following equation*

$$\widehat{D}^*g(x)(y^*) = \{(\nabla g(x))^*(y^*)\}, \text{ for all } y^* \in Y^*.$$

One of the important applications of Mordukhovich derivatives of set-valued mappings is to define the covering constants for set-valued mappings. The covering constant for $\Phi: X \rightrightarrows Y$ at point $(\bar{x}, \bar{y}) \in \text{gph } \Phi$ is defined by (see (2.4) in [1], also see [18])

$$\hat{\alpha}(\Phi, \bar{x}, \bar{y}) := \sup_{\eta > 0} \inf \{ \|z^*\|_{X^*}: z^* \in \widehat{D}^*\Phi(x,y)(w^*), x \in \mathbb{B}_X(\bar{x}, \eta), y \in \Phi(x) \cap \mathbb{B}_Y(\bar{y}, \eta), \|w^*\|_{Y^*} = 1 \}. \tag{2.4}$$

Here, $\|\cdot\|_{X^*}$ and $\|\cdot\|_{Y^*}$ denote the norms in $X^*$ and $Y^*$, respectively. $\mathbb{B}_X(\bar{x}, \eta)$ is the closed ball in $X$ centered at $\bar{x}$ with radius $\eta$, and $\mathbb{B}_Y(\bar{y}, \eta)$ is the closed ball in $Y$ centered at $\bar{y}$ with radius $\eta$.

In particular, let $g: X \to Y$ be a single-valued mapping. For any $\bar{x}, \bar{y} \in X$ with $\bar{y} = g(\bar{x})$, (2.4) becomes

$$\hat{\alpha}(g, \bar{x}, \bar{y}) = \sup_{\eta > 0} \inf \{ \|z^*\|_{X^*}: z^* \in \widehat{D}^*g(x)(w^*), x \in \mathbb{B}_X(\bar{x}, \eta), g(x) \in \mathbb{B}_Y(\bar{y}, \eta), \|w^*\|_{Y^*} = 1 \}. \tag{2.5}$$

In particular, in this paper, we consider finite dimensional Hilbert spaces (Euclidean spaces) as special cases of Banach spaces. Let $(H, \|\cdot\|)$ be a real Hilbert space with inner product $\langle \cdot, \cdot \rangle$ and origin $\theta$. For $r > 0$ and $x \in H$, let $\mathbb{B}(x, r)$ be the closed ball in $H$ with radius $r$ and centered at $x$. Let $g: H \to H$ be a single-valued mapping and let $\bar{x} \in X$. Then, the Fréchet derivative of $g$ at $\bar{x}$ is defined by $\nabla g(\bar{x}): H \to H$ such that

$$\lim_{x \to \bar{x}} \frac{g(x) - g(\bar{x}) - \nabla g(\bar{x})(x - \bar{x})}{\|x - \bar{x}\|} = \theta.$$

Let $\Delta$ be a nonempty subset of $H$ and let $g: \Delta \to H$ be a single-valued continuous mapping. For $x \in \Delta$ and $y = g(x)$, by (2.2), the Mordukhovich derivative of $g$ at point $x$ becomes

$$\widehat{D}^* g(x)(y) = \left\{ z \in H : \limsup_{\substack{u \to x \\ u \in \Delta}} \frac{\langle z, u - x \rangle - \langle y, g(u) - g(x) \rangle}{\|u - x\| + \|g(u) - g(x)\|} \leq 0 \right\}, \text{ for any } y \in H. \tag{2.6}$$

In this case, for any $\bar{x}, \bar{y} \in H$ with $\bar{y} = g(\bar{x})$, by (2.4), the covering constant for $g$ at point $(\bar{x}, \bar{y})$ is that

$$\hat{a}(g, \bar{x}, \bar{y}) = \sup_{\eta > 0} \inf \{ \|z\| : z \in \widehat{D}^* g(x)(w), x \in \mathbb{B}(\bar{x}, \eta), g(x) \in \mathbb{B}(\bar{y}, \eta), \|w\| = 1 \}. \tag{2.7}$$

Next, we review the concepts of covering properties.

**Definition 1.51 in [21]** (covering properties) Let $X$, $Y$ be Banach spaces. Let $\mathbb{B}_X$ and $\mathbb{B}_Y$ denote the unit balls in $X$, $Y$, respectively. Let $F: X \rightrightarrows Y$ with dom $F \neq \emptyset$.

(i) Let $U$ and $V$ be nonempty subsets in $X$ and $Y$, respectively. We say that $F$ enjoys the covering property on $U$ relative to $V$ if there is $\gamma > 0$ such that

$$F(x) \cap V + \gamma r \mathbb{B}_Y \subset F(x + r \mathbb{B}_X), \text{ whenever } x + r \mathbb{B}_X \subset U, \text{ as } r > 0. \tag{2.8}$$

(ii) Given $(\bar{x}, \bar{y}) \in \text{gph } F$, we say that $F$ has the local covering property around $(\bar{x}, \bar{y})$ with modulus $\gamma > 0$ if there is a neighborhood $U$ of $\bar{x}$ and a neighborhood $V$ of $\bar{y}$ such that the above inclusion in (i) holds. The supremum of all such moduli $\{\gamma\}$ is called the exact covering bound of $F$ around $(\bar{x}, \bar{y})$, which is denoted by

$$\text{cov} F(\bar{x}, \bar{y}) = \sup\{\gamma : \gamma \text{ satisfies (2.8) for some } U \subset X, V \subset Y\}.$$

**Theorem 4.1 in [21]** (neighborhood characterization of local covering) *Let $X$, $Y$ be Asplund spaces. Let $F: X \rightrightarrows Y$ be a set-valued mapping. Assume that $F$ is closed-graph around $(\bar{x}, \bar{y}) \in \text{gph } F$. Then the following are equivalent*:

(a) *$F$ enjoys the local covering property around $(\bar{x}, \bar{y})$ (that is, $\text{cov } F(\bar{x}, \bar{y}) > 0$).*
(b) *One has $\hat{a}(F, \bar{x}, \bar{y}) > 0$.*

*Moreover, the exact covering bound of $F$ around $(\bar{x}, \bar{y})$ is computed by*

$$\text{cov} F(\bar{x}, \bar{y}) = \hat{a}(F, \bar{x}, \bar{y}).$$

## 3. Fréchet and Mordukhovich Derivatives (Coderivatives) of Mappings from $\mathbb{R}^n$ to $\mathbb{R}^m$

Let $n \geq 1$ and let $(\mathbb{R}^n, \|\cdot\|_n)$ be the standard $n$-d Euclidean space with the ordinal Hilbert $L_2$-norm $\|\cdot\|_n$ and row vectors. Let $\theta$ denote the origin of $\mathbb{R}^n$. Let $m, n \geq 1$ and let $f: \mathbb{R}^n \to \mathbb{R}^m$ be a single-valued mapping with the following representation.

$$f((x_1, x_2, \ldots, x_n)) = \big(f_1(x_1, x_2, \ldots, x_n), f_2(x_1, x_2, \ldots, x_n), \ldots, f_m(x_1, x_2, \ldots, x_n)\big), \text{ for } (x_1, x_2, \ldots, x_n) \in \mathbb{R}^n.$$

Where, for $i = 1, 2, \ldots, m$, $f_i(x_1, x_2, \ldots, x_n): \mathbb{R}^n \to \mathbb{R}$ is a real valued multivariable function defined on $\mathbb{R}^n$ with variables $x_1, x_2, \ldots, x_n$. We write the above equation as $f = (f_1, f_2, \ldots, f_m)$. Let $(z_1, z_2, \ldots, z_n) \in \mathbb{R}^n$. For $i = 1, 2, \ldots, n$, the partial derivative of $f_i$ with respect to the variable $x_j$ at $(z_1, z_2, \ldots, z_n)$ is

$$\frac{\partial f_i}{\partial x_j}(z_1, z_2, \ldots, z_n), \text{ for every } j = 1, 2, \ldots, n. \tag{3.1}$$

By equation (3.1), for any $i = 1, 2, \ldots, m$, the existence of $f_i(z_1, z_2, \ldots, z_n)$ means that,

$$\lim_{x_j \to z_j} \frac{f_i(x_1, x_2, \ldots, x_n) - f_i(z_1, z_2, \ldots, z_n) - \frac{\partial f_i}{\partial x_j}(z_1, z_2, \ldots, z_n)(x_j - z_j)}{x_j - z_j} = 0, \text{ for } j = 1, 2, \ldots, n. \tag{3.2}$$

The limit (3.2) is equivalent to, for any $i = 1, 2, \ldots, m$,

$$\lim_{x_j \to z_j} \frac{f_i(x_1, x_2, \ldots, x_n) - f_i(z_1, z_2, \ldots, z_n)}{x_j - z_j} = -\frac{\partial f_i}{\partial x_j}(z_1, z_2, \ldots, z_n), \text{ for } j = 1, 2, \ldots, n. \tag{3.3}$$

Furthermore, if $f$ satisfies the following conditions

$$\lim_{x \to z} \frac{f_i(x_1, x_2, \ldots, x_n) - f_i(z_1, z_2, \ldots, z_n) - \left(\sum_{j=1}^n \frac{\partial f_i}{\partial x_j}(z_1, z_2, \ldots, z_n)(x_j - z_j)\right)}{\|x - z\|_n} = 0, \text{ for } i = 1, 2, \ldots, m, \tag{3.4}$$

then, $f$ is differentiable and has the linear approximation at point $z = (z_1, z_2, \ldots, z_n)$, which is the first order of the Taylor polynomial of the mapping $f$ from $\mathbb{R}^n$ to $\mathbb{R}^m$. We have some sufficient conditions for $f$ to have the linear approximation at point $z$. For every $i = 1, 2, \ldots, m$, the second order partial derivative of $f_i$ at point $z = (z_1, z_2, \ldots, z_n)$ with respect to $x_j$ and $x_k$ is

$$\frac{\partial^2 f_i}{\partial x_j \partial x_k}(z_1, z_2, \ldots, z_n), \text{ for } j, k = 1, 2, \ldots, n.$$

**Fact 3.1.** *Let $f = (f_1, f_2, \ldots, f_m). : \mathbb{R}^n \to \mathbb{R}^m$ be a single-valued mapping. Let $z \in \mathbb{R}^n$. If there is a ball $B$ with radius $r > 0$ and centered at $z$ such that, for every $i = 1, 2, \ldots, m$, the real valued function $f_i: \mathbb{R}^n \to \mathbb{R}$ is twice differentiable (all second partial derivatives of $f_i$ exist) at every point $y \in B$, that is*

$$\frac{\partial^2 f_i}{\partial x_j \partial x_k}(z_1, z_2, \ldots, z_n) \text{ exists for any } z \in B, \text{ for } j, k = 1, 2, \ldots, n.$$

*then $f$ has the linear approximation at this point $z$.*

*Proof.* This is a well-known result. It is proved by using the second order Taylor polynomial of each real valued function $f_i: \mathbb{R}^n \to \mathbb{R}$, for $i = 1, 2, \ldots, n$. The details are omitted here.  □

**Theorem 3.2.** *Let $f = (f_1, f_2, \ldots, f_m): \mathbb{R}^n \to \mathbb{R}^m$ be a single-valued mapping. Let $(z_1, z_2, \ldots, z_n) \in \mathbb{R}^n$. Suppose that, for every $i = 1, 2, \ldots, m$, $\frac{\partial f_i}{\partial x_j}(z_1, z_2, \ldots, z_n)$ exists, for every $j = 1, 2, \ldots, n$. Suppose that $f$ has the linear approximation at point $z$. Then,*

(a) *$f$ is Fréchet differentiable at $z$ and the Fréchet derivative of $f$ at $z$ is the following $n \times m$ matrix,*

$$\nabla f(z) = \begin{pmatrix} \frac{\partial f_1}{\partial x_1}(z_1, z_2, \ldots, z_n) & \cdots & \frac{\partial f_m}{\partial x_1}(z_1, z_2, \ldots, z_n) \\ \vdots & \ddots & \vdots \\ \frac{\partial f_1}{\partial x_n}(z_1, z_2, \ldots, z_n) & \cdots & \frac{\partial f_m}{\partial x_n}(z_1, z_2, \ldots, z_n) \end{pmatrix}, \quad (3.5)$$

(b) $f$ is Mordukhovich differentiable at point $z$ that is the Jacobian matrix of $f$ at $z$

$$\widehat{D}^* f(x) = \nabla f(z)^T = \begin{pmatrix} \frac{\partial f_1}{\partial x_1}(z_1, z_2, \ldots, z_n) & \cdots & \frac{\partial f_1}{\partial x_n}(z_1, z_2, \ldots, z_n) \\ \vdots & \ddots & \vdots \\ \frac{\partial f_m}{\partial x_1}(z_1, z_2, \ldots, z_n) & \cdots & \frac{\partial f_m}{\partial x_n}(z_1, z_2, \ldots, z_n) \end{pmatrix}. \quad (3.6)$$

*Proof.* Proof of (a). Notice that, in this paper, $\mathbb{R}^n$ has row vectors. By the meanings of (3.2) and (3.3), $f$ has the linear approximation at this point $z$, by (3.4), we calculate

$$\lim_{x \to z} \frac{f(x_1, x_2, \ldots, x_n) - f(z_1, z_2, \ldots, z_n) - \nabla f(z)(x - z)}{\|x - z\|_n}$$

$$= \lim_{x \to z} \left( \frac{\big(f_1(x_1,x_2,\ldots,x_n), f_2(x_1,x_2,\ldots,x_n), \ldots, f_m(x_1,x_2,\ldots,x_n)\big) - \big(f_1(z_1,z_2,\ldots,z_n), f_2(z_1,z_2,\ldots,z_n), \ldots, f_m(z_1,z_2,\ldots,z_n)\big)}{\|x - z\|_n} \right.$$

$$\left. - \frac{(x_1 - z_1, \ldots, x_n - z_n) \begin{pmatrix} \frac{\partial f_1}{\partial x_1}(z_1,z_2,\ldots,z_n) & \cdots & \frac{\partial f_m}{\partial x_1}(z_1,z_2,\ldots,z_n) \\ \vdots & \ddots & \vdots \\ \frac{\partial f_1}{\partial x_n}(z_1,z_2,\ldots,z_n) & \cdots & \frac{\partial f_m}{\partial x_n}(z_1,z_2,\ldots,z_n) \end{pmatrix}}{\|x - z\|_n} \right)$$

$$= \lim_{x \to z} \left( \frac{f_1(x_1,x_2,\ldots,x_n) - f_1(z_1,z_2,\ldots,z_n) - \left(\sum_{j=1}^n (x_j - z_j)\frac{\partial f_1}{\partial x_j}(z_1,z_2,\ldots,z_n)\right)}{\|x - z\|_n}, \ldots, \right.$$

$$\left. \frac{f_m(x_1, x_2, \ldots, x_n) - f_m(z_1, z_2, \ldots, z_n) - \left(\sum_{j=1}^n (x_j - z_j)\frac{\partial f_m}{\partial x_1}(z_1, z_2, \ldots, z_n)\right)}{\|x - z\|_n} \right)$$

$= (0, \ldots, 0)$.

Part (b) can be proved by Theorem 1.38 in [18] and Part (a) of this theorem. □

The results of the following proposition can be reduced by Theorem 3.2 immediately.

**Proposition 3.3.** Let $f = (f_1, f_2, \ldots, f_m): \mathbb{R}^n \to \mathbb{R}^m$ be a single-valued mapping.

(a) *Suppose that, for every $i = 1, 2, \ldots, m$, $f_i(x_1, x_2, \ldots, x_n)$ is a polynomial function with respect to $x_1, x_2, \ldots, x_n$, then $f$ is Fréchet differentiable at every point in $\mathbb{R}^n$, and $f$ is Mordukhovich differentiable on $\mathbb{R}^n$;*

(b) *Suppose that, for every $i = 1, 2, \ldots, m$, $f_i(x_1, x_2, \ldots, x_n)$ is a rational function with respect to $x_1, x_2, \ldots, x_n$. Let $z \in \mathbb{R}^n$. If $z$ is not a zero point for the denominator of every $f_i$, then $f$ is Fréchet differentiable at $z$; and therefore, $f$ is Mordukhovich differentiable at $z$.*

## 4. Covering Constants for Single-Valued Mappings from $\mathbb{R}^n$ to $\mathbb{R}^m$ with $n < m$

In this section, as applications of Lemma 3.1 and Theorems 3.2 and 3.3 in the previous section, we use the Fréchet and Mordukhovich derivatives to prove that, for any single-valued mapping $f$ from $\mathbb{R}^n$ to $\mathbb{R}^m$, if $n < m$, then the covering constant for $f$ satisfies that $\hat{\alpha}(f, \bar{z}, \bar{w}) = 0$, for $\bar{z} \in \mathbb{R}^n, \bar{w} \in \mathbb{R}^m$ with $f(\bar{z}) = \bar{w}$.

**Theorem 4.1**. *Let $n < m$. Let $f = (f_1, f_2, \ldots, f_m): \mathbb{R}^n \to \mathbb{R}^m$ be a single-valued mapping. Let $\bar{z} \in \mathbb{R}^n$ and $\bar{w} \in \mathbb{R}^m$ with $\bar{w} = f(\bar{z})$. Suppose that $f$ is Fréchet differentiable at and around $\bar{z}$. Then,*

$$\hat{\alpha}(f, \bar{z}, \bar{w}) = 0.$$

*Proof.* Suppose that the Fréchet derivative of $f$ at $z \in \mathbb{R}^n$ nearing to $\bar{z}$ is given by the following $n \times m$ matrix,

$$\nabla f(z) = \begin{pmatrix} \frac{\partial f_1}{\partial x_1}(z_1, z_2, \ldots, z_n) & \cdots & \frac{\partial f_m}{\partial x_1}(z_1, z_2, \ldots, z_n) \\ \vdots & \ddots & \vdots \\ \frac{\partial f_1}{\partial x_n}(z_1, z_2, \ldots, z_n) & \cdots & \frac{\partial f_m}{\partial x_n}(z_1, z_2, \ldots, z_n) \end{pmatrix}. \tag{4.1}$$

By Theorem 1.38 in [18] and the condition (4.1) in this proposition, we have that $f$ is Mordukhovich differentiable at point $z$ and the Mordukhovich derivative of $f$ at $z$ is the following $m \times n$ matrix (Jacobian matrix of $f$ at $z$)

$$\widehat{D}^* f(x) = \nabla f(z)^T = \begin{pmatrix} \frac{\partial f_1}{\partial x_1}(z_1, z_2, \ldots, z_n) & \cdots & \frac{\partial f_1}{\partial x_n}(z_1, z_2, \ldots, z_n) \\ \vdots & \ddots & \vdots \\ \frac{\partial f_m}{\partial x_1}(z_1, z_2, \ldots, z_n) & \cdots & \frac{\partial f_m}{\partial x_n}(z_1, z_2, \ldots, z_n) \end{pmatrix}. \tag{4.2}$$

Let $z = (z_1, z_2, \ldots, z_n) \in \mathbb{R}^n$. Let $x = (x_1, x_2, \ldots, x_n) \in \mathbb{R}^n$ and $u = (u_1, u_2, \ldots, u_m) \in \mathbb{R}^m$. We consider the following system of $n$ linear equations with respect to $m$ variables $u_1, u_2, \ldots, u_m$.

$$(0, 0, \ldots, 0) = (u_1, u_2, \ldots, u_m) \begin{pmatrix} \frac{\partial f_1}{\partial x_1}(z_1, z_2, \ldots, z_n) & \cdots & \frac{\partial f_1}{\partial x_n}(z_1, z_2, \ldots, z_n) \\ \vdots & \ddots & \vdots \\ \frac{\partial f_m}{\partial x_1}(z_1, z_2, \ldots, z_n) & \cdots & \frac{\partial f_m}{\partial x_n}(z_1, z_2, \ldots, z_n) \end{pmatrix}. \tag{4.3}$$

Since $n < m$, the system (4.3) of linear equations has nonzero solutions. For any given nonzero solution $u = (u_1, u_2, \ldots, u_m)$ of (4.3), let $y_i = \frac{u_i}{\sqrt{\sum_{j=1}^m u_j^2}}$, for $i = 1, 2, \ldots, m$. Then $y = (y_1, y_2, \ldots, y_m) \in \mathbb{R}^m$ and $y$ is also a nonzero solution of (4.3) satisfying $\|y\|_m = 1$. This implies that, for $z = (z_1, z_2, \ldots, z_n) \in \mathbb{R}^n$, $\widehat{D}^* f(z)$ exists with representation (4.2). Then, there is $y = (y_1, y_2, \ldots, y_m) \in \mathbb{R}^m$ with $\|y\|_m = 1$ such that

$$\widehat{D}^* f(z)(y) = (0, 0, \ldots, 0). \tag{4.4}$$

By (4.2), we calculate

$$\hat{\alpha}(f, \bar{z}, \bar{w}) = \sup_{\eta > 0} \inf \{\|x\|_n : x = \widehat{D}^* f(z)(y), z \in \mathbb{B}(\bar{z}, \eta), f(z) \in \mathbb{B}(\bar{w}, \eta), \|y\|_m = 1\}$$

$$\leq \sup_{\eta>0} \inf\{\|x\|_n \colon x = \widehat{D}^*f(z)(y), z \in \mathbb{B}(\bar{z},\eta), f(z) \in \mathbb{B}(\bar{w},\eta), \|y\|_m = 1, y \text{ satisfies } (4.4)\}$$

$$= \sup_{\eta>0} \inf\{\|\theta\|_n \colon \theta = \widehat{D}^*f(z)(y), z \in \mathbb{B}(\bar{z},\eta), f(z) \in \mathbb{B}(\bar{w},\eta), \|y\|_m = 1, y \text{ satisfies } (4.4)\}$$

$$= \sup_{\eta>0} \inf\{0 \colon \theta = \widehat{D}^*f(z)(y), z \in \mathbb{B}(\bar{z},\eta), f(z) \in \mathbb{B}(\bar{w},\eta), \|y\|_m = 1, y \text{ satisfies } (4.4)\}$$

$$= 0. \qquad \square$$

When we consider single-valued mappings are special cases of set-valued mappings with values of singletons, by the connection between covering constant and covering property (Theorem 4.1 in [21]), we have the following consequence of Theorem 4.1.

**Corollary 4.2.** *Let $n < m$. Let $f \colon \mathbb{R}^n \to \mathbb{R}^m$ be a single-valued mapping. Let $\bar{z} \in \mathbb{R}^n$ and $\bar{w} \in \mathbb{R}^m$ with $\bar{w} = f(\bar{z})$. Suppose that $f$ is Fréchet differentiable at and around $\bar{z}$. Then, for any nonempty subsets $U$ and $V$ in $X$ and $Y$, respectively, $f$ does not enjoy the covering property on $U$ relative to $V$. That is,*

$$\mathrm{cov} f(\bar{z}, \bar{w}) = 0.$$

*Proof.* By Theorem 4.1, we have $\hat{\alpha}(f, \bar{z}, \bar{w}) = 0$. The Fréchet differentiability of $f$ at and around $\bar{z}$ implies the continuity of $f$ at and around $\bar{z}$. Hence, $f$ has the closed graph property of at point $(\bar{z}, \bar{w})$. Then, by Theorem 4.1 in [21] and Theorem 4.1, we have $\mathrm{cov} f(\bar{z}, \bar{w}) = \hat{\alpha}(f, \bar{z}, \bar{w}) = 0$. By Definition 1.51 in [21] (covering properties), we obtain that for any nonempty subsets $U$ and $V$ in $X$ and $Y$, respectively, $f$ does not enjoy the covering property on $U$ relative to $V$. $\square$

Next, we provide some examples to demonstrate the results of Theorem 4.1.

**4.3. A linear and continuous norm preserving mapping from $\mathbb{R}^2$ to $\mathbb{R}^3$.** We define $f \colon \mathbb{R}^2 \to \mathbb{R}^3$, for every $x = (x_1, x_2) \in \mathbb{R}^2$, by

$$f(x) = \left(x_1, \frac{x_2}{\sqrt{2}}, \frac{x_2}{\sqrt{2}}\right).$$

This is a linear and continuous mapping, for $z = (z_1, z_2) \in \mathbb{R}^2$, we have

$$\nabla f(z) \colon \mathbb{R}^2 \to \mathbb{R}^3 \quad \text{and} \quad \widehat{D}^*f(z) \colon \mathbb{R}^3 (\text{It is } (\mathbb{R}^3)^*) \to \mathbb{R}^2 (\text{It is } (\mathbb{R}^2)^*).$$

Then, $f$ has the following properties.

(a) $f$ is norm preserving, this is,

$$\|f(x)\|_3 = \|x\|_2, \quad \text{for any } x \in \mathbb{R}^2.$$

(b) $f$ is Fréchet differentiable and Mordukhovich differentiable on $\mathbb{R}^2 \setminus \{\theta\}$. For any $z = (z_1, z_2) \in \mathbb{R}^3$, we have

$$\nabla f(z) = \begin{pmatrix} 1 & 0 & 0 \\ 0 & \frac{1}{\sqrt{2}} & \frac{1}{\sqrt{2}} \end{pmatrix} \quad \text{and} \quad \widehat{D}^*f(z) = \begin{pmatrix} 1 & 0 \\ 0 & \frac{1}{\sqrt{2}} \\ 0 & \frac{1}{\sqrt{2}} \end{pmatrix}.$$

(c) The covering constant for $f$ is constant on $\mathbb{R}^3$ satisfying

$$\hat{\alpha}(f,\bar{z},\bar{w}) = 0, \text{ for any } \bar{z} = (\bar{z}_1,\bar{z}_2) \in \mathbb{R}^3 \text{ with } \bar{w} = f(\bar{z}) \in \mathbb{R}^2.$$

*Proof.* The proofs of (a) and (b) are straight forward and it is omitted here. We only prove (c). Let $z = (z_1, z_2) \in \mathbb{R}^2$. Let $x = (x_1, x_2) \in \mathbb{R}^2$ and $y = (y_1, y_2, y_3) \in \mathbb{R}^3$. If $x = \widehat{D}^* f(z)(y)$, by part (b), we have that

$$(x_1, x_2) = (y_1, y_2, y_3) \begin{pmatrix} 1 & 0 \\ 0 & \frac{1}{\sqrt{2}} \\ 0 & \frac{1}{\sqrt{2}} \end{pmatrix} = \left(y_1, \frac{y_2}{\sqrt{2}}, \frac{y_3}{\sqrt{2}}\right).$$

This reduces that

$$x = \widehat{D}^* f(z)(y) \quad \Longrightarrow \quad \|x\|_2^2 = y_1^2 + \left(\frac{y_2+y_3}{\sqrt{2}}\right)^2. \tag{4.5}$$

Let $\bar{z} = (\bar{z}_1, \bar{z}_2) \in \mathbb{R}^2$ with $\bar{w} = f(\bar{z}) \in \mathbb{R}^3$. By (4.5), we calculate

$$\hat{\alpha}(f,\bar{z},\bar{w}) = \sup_{\eta>0} \inf\{\|x\|_2 : x = \widehat{D}^* f(z)(y), z \in \mathbb{B}(\bar{z},\eta), f(z) \in \mathbb{B}(\bar{w},\eta), \|y\|_3 = 1\}$$

$$= \sup_{\eta>0} \inf\left\{\sqrt{y_1^2 + \left(\frac{y_2+y_3}{\sqrt{2}}\right)^2} : x = \widehat{D}^* f(z)(y), z \in \mathbb{B}(\bar{z},\eta), f(z) \in \mathbb{B}(\bar{w},\eta), \|y\|_3 = 1\right\}$$

$$\leq \sup_{\eta>0} \inf\left\{\sqrt{y_1^2 + \left(\frac{y_2+y_3}{\sqrt{2}}\right)^2} : x = \widehat{D}^* f(z)(y), z \in \mathbb{B}(\bar{z},\eta), f(z) \in \mathbb{B}(\bar{w},\eta), \|y\|_3 = 1, y_1 = 0, y_2 = -y_3 = \frac{1}{\sqrt{2}}\right\}$$

$$= 0. \qquad \square$$

**4.4. A continuous mapping from $\mathbb{R}^2$ to $\mathbb{R}^3$.** Define $f: \mathbb{R}^2 \to \mathbb{R}^3$, for $x = (x_1, x_2) \in \mathbb{R}^2$, by

$$f(x) = (x_1, x_1 x_2, x_2).$$

This is a continuous mapping. $f$ has the following properties.

(a) $f$ is a norm expanding mapping, this is,

$$\|f(x)\|_3 \geq \|x\|_2, \text{ for any } x \in \mathbb{R}^2.$$

(b) $f$ is Fréchet differentiable and Mordukhovich differentiable on $\mathbb{R}^2 \setminus \{\theta\}$. For any $z = (z_1, z_2) \in \mathbb{R}^3$, we have

$$\nabla f(z) = \begin{pmatrix} 1 & z_2 & 0 \\ 0 & z_1 & 1 \end{pmatrix} \quad \text{and} \quad \widehat{D}^* f(z) = \begin{pmatrix} 1 & 0 \\ z_2 & z_1 \\ 0 & 1 \end{pmatrix}.$$

(c) The covering constant for $f$ is constant on $\mathbb{R}^3$ satisfying

$$\hat{\alpha}(f,\bar{z},\bar{w}) = 0, \text{ for any } \bar{z} = (\bar{z}_1,\bar{z}_2) \in \mathbb{R}^3 \text{ with } \bar{w} = f(\bar{z}) \in \mathbb{R}^2.$$

*Proof.* The proofs of (a) and (b) are straight forward and it is omitted here. We only prove (c). Let $z =$

$(z_1, z_2) \in \mathbb{R}^2$. Let $x = (x_1, x_2) \in \mathbb{R}^2$ and $y = (y_1, y_2, y_3) \in \mathbb{R}^3$. If $x = \widehat{D}^* f(z)(y)$, by part (b), we have

$$(x_1, x_2) = (y_1, y_2, y_3) \begin{pmatrix} 1 & 0 \\ z_2 & z_1 \\ 0 & 1 \end{pmatrix} = (y_1 + y_2 z_2, y_2 z_1 + y_3).$$

This induces that

$$x = \widehat{D}^* f(z)(y) \implies \|x\|_2^2 = (y_1 + y_2 z_2)^2 + (y_2 z_1 + y_3)^2. \tag{4.6}$$

Let $\bar{z} = (\bar{z}_1, \bar{z}_2) \in \mathbb{R}^2$ with $\bar{w} = f(\bar{z}) \in \mathbb{R}^3$. By (4.6), we calculate

$$\hat{\alpha}(f, \bar{z}, \bar{w}) = \sup_{\eta > 0} \inf\{\|x\|_2 : x = \widehat{D}^* f(z)(y), z \in \mathbb{B}(\bar{z}, \eta), f(z) \in \mathbb{B}(\bar{w}, \eta), \|y\|_3 = 1\}$$

$$\leq \sup_{\eta > 0} \inf\left\{\sqrt{(y_1 + y_2 z_2)^2 + (y_2 z_1 + y_3)^2} : x = \widehat{D}^* f(z)(y), z \in \mathbb{B}(\bar{z}, \eta), f(z) \in \mathbb{B}(\bar{w}, \eta), \|y\|_3 = 1, y_1 + y_2 z_2 = y_2 z_1 + y_3 = 0\right\}$$

$$\leq \sup_{\eta > 0} \inf\{0 : x = \widehat{D}^* f(z)(y), z \in \mathbb{B}(\bar{z}, \eta), f(z) \in \mathbb{B}(\bar{w}, \eta), \|y\|_3 = 1, y_1 + y_2 z_2 = y_2 z_1 + y_3 = 0\}$$

$= 0.$ □

## 5. Covering Constants for Some Norm Preserving Mappings

In this section, we first consider the single-valued mapping $f: \mathbb{R}^2 \to \mathbb{R}^2$ defined by (5.1), which is studied in Examples 2 in [4] and 4.2 in [2]. We will calculate the Fréchet and Mordukhovich derivatives of $f$, by which we will directly find the exact covering constant for $f$. In this section, $f: \mathbb{R}^2 \to \mathbb{R}^2$ is defined by

$$f(x_1, x_2) = \left(\frac{x_1^2 - x_2^2}{\sqrt{x_1^2 + x_2^2}}, \frac{2 x_1 x_2}{\sqrt{x_1^2 + x_2^2}}\right), \text{ for } (x_1, x_2) \in \mathbb{R}^2 \setminus \{\theta\} \text{ with } f(\theta) = \theta. \tag{5.1}$$

It is clear that $\|f(\theta)\| = \|\theta\|$. For any $x = (x_1, x_2) \in \mathbb{R}^2$ with $(x_1, x_2) \neq \theta$, we calculate

$$\|f((x_1, x_2))\|^2 = \left(\frac{x_1^2 - x_2^2}{\sqrt{x_1^2 + x_2^2}}\right)^2 + \left(\frac{2 x_1 x_2}{\sqrt{x_1^2 + x_2^2}}\right)^2 = x_1^2 + x_2^2 = \|(x_1, x_2)\|^2.$$

This implies that the mapping $f: \mathbb{R}^2 \to \mathbb{R}^2$ is norm preserving. That is,

$$\|f(x)\| = \|x\|, \text{ for any } x \in \mathbb{R}^2.$$

By the norm preserving property of $f$, it follows immediately that $f$ is a continuous mapping on $\mathbb{R}^2$, which means that $f$ is continuous at origin $\theta$ in $\mathbb{R}^2$. This mapping has been used several times (at least two times as what I known). In this section, we will directly prove the following properties:

$$\hat{\alpha}(f, \bar{z}, f(\bar{z})) = 1, \text{ for any } \bar{z} \in \mathbb{R}^2. \tag{5.2}$$

In [2], in order to prove (5.2), $f$ is converted to a mapping $\gamma$ in the complex plane $\mathbb{C}$ by

$$\gamma(z) = \frac{z^2}{\|z\|}, \text{ for any } z \in \mathbb{C} \setminus \{\theta\} \text{ with } \gamma(\theta) = \theta.$$

Then, it was geometrically proved that the exact covering bound of $\gamma$ satisfies that $\text{cov}\gamma(x, y) = 1$. for $x, y \in \mathbb{C}$ with $y = \gamma(x)$. Then, by the connection between $\text{cov}\gamma(x, y)$ and $\alpha(f, x, y)$ proved in Theorem 4.1, Corollaries 4.2, 4.3 and Proposition 3.6 in [22], we have that

$$\hat{\alpha}(f, \bar{x}, \bar{y}) = \text{cov}\, f(\bar{x}, \bar{y}) = \alpha(f, \bar{x}, \bar{y}) = 1, \text{ for } \bar{x}, \bar{y} \in \mathbb{R}^2 \text{ with } \bar{y} = f(\bar{x}).$$

which implies (5.2). In this section, we will give a direct proof for (5.2) with more details and properties of this mapping $f$. Next, we first prove that at the origin $\theta$ in $\mathbb{R}^2$, $f$ is neither Fréchet differentiable, nor Mordukhovich differentiable.

Notice that after we prove that, $f$ is not Mordukhovich differentiable at point $\theta$, then, by the property that Fréchet differentiable at point $\theta$ implies that $f$ is Mordukhovich differentiable at point $\theta$ (Theorem 1.38 in [18]), we obtain that $f$ is not Fréchet differentiable at point $\theta$. However, in order to study the details for $f$ is not Fréchet Fréchet differentiable at point $\theta$, we will give a direct proof.

**Lemma 5.1.** *$f$ is not Fréchet differentiable at $\theta$.*

*Proof.* Assume, by the way of contradiction, that $f$ is Fréchet differentiable at $\theta$ and $\nabla f(\theta): \mathbb{R}^2 \to \mathbb{R}^2$ exists, which must be a linear and continuous mapping on $\mathbb{R}^2$. Then, $\nabla f(\theta)$ must be precisely represented by a real $2 \times 2$ matrix $\begin{pmatrix} a & b \\ c & d \end{pmatrix}$, for some real numbers $a$, $b$, $c$, and $d$.

Assume $a \neq 0$. In this case, we take a special direction in the following limit for $u \to \theta$ by $u_1 = -at$ with $t \downarrow 0$, $u_2 = 0$. By the assumption that $\nabla f(\theta) = \begin{pmatrix} a & b \\ c & d \end{pmatrix}$ and by $f(\theta) = \theta$, it follows that

$$\lim_{u=(u_1,u_2) \to \theta} \frac{f(u)-f(\theta)- \nabla f(z)(u-\theta)}{\|u-z\|}$$

$$= \lim_{t \downarrow 0} \frac{\left(\frac{a^2 t^2}{\sqrt{a^2 t^2}},\ 0\right)-(-at,\ 0)\begin{pmatrix} a & b \\ c & d \end{pmatrix}}{\|(-at,0)\|}$$

$$= \lim_{t \downarrow 0} \frac{(|a|t,\ 0)+(a^2 t,\ abt)}{|a|t}$$

$$= (1+|a|,\ (\text{sign}\,a)b)$$

$$\neq \theta.$$

This contradicts to the assumption that $\nabla f(\theta) = \begin{pmatrix} a & b \\ c & d \end{pmatrix}$.

Next, suppose $a = 0$. In this case, we take a special direction in the limit (5.3) for $u \to \theta$ by $u_1 \downarrow 0$, $u_2 = 0$, $u_3 = z_3$ and $u_4 = z_4$. It follows that

$$\lim_{u=(u_1,u_2) \to \theta} \frac{f(u)-f(\theta)- \nabla f(z)(u-\theta)}{\|u-z\|}$$

$$= \lim_{u_1 \downarrow 0} \frac{\left(\frac{u_1^2}{\sqrt{u_1^2}},\ 0\right)-(u_1,\ 0)\begin{pmatrix} 0 & b \\ c & d \end{pmatrix}}{\|(u_1,0)\|}$$

$$= \lim_{u_1 \downarrow 0} \frac{(1, 0)-(0, bu_1)}{u_1}$$

$$= (1, b)$$

$$\neq \theta.$$

This contradicts the assumption that $\nabla f(\theta) = \begin{pmatrix} a & b \\ c & d \end{pmatrix}$. It completes the proof. □

**Proposition 5.2.** *f is not Mordukhovich differentiable at $\theta$. More precisely speaking, we have*

$$\widehat{D}^* f(\theta)(y) = \emptyset, \quad \text{for any } y \in \mathbb{R}^2 \setminus \{\theta\}.$$

*Proof.* Let $y = (y_1, y_2) \in \mathbb{R}^2$ and $x = (x_1, x_2) \in \mathbb{R}^2$, suppose that $x \in \widehat{D}^* f(\theta)(y)$. We calculate the following limits.

$$0 \geq \limsup_{u \to \theta} \frac{\langle (x_1,x_2), (u_1,u_2)-\theta \rangle - \langle (y_1,y_2), f(u)-f(\theta) \rangle}{\|u-\theta\| + \|f(u)-f(\theta)\|}$$

$$= \limsup_{u \to \theta} \frac{\langle (x_1,x_2), (u_1,u_2) \rangle - \langle (y_1,y_2), \left( \frac{u_1^2-u_2^2}{\sqrt{u_1^2+u_2^2}}, \frac{2u_1 u_2}{\sqrt{u_1^2+u_2^2}} \right) \rangle}{\|u\| + \|f(u)\|}$$

$$= \limsup_{u \to \theta} \frac{x_1 u_1 + x_2 u_2 - y_1 \frac{u_1^2-u_2^2}{\sqrt{u_1^2+u_2^2}} - y_2 \frac{2u_1 u_2}{\sqrt{u_1^2+u_2^2}}}{2\|u\|}. \tag{5.3}$$

Case 1. $x_1 > y_1$. In this case, we take a special direction in the limit (5.3) for $u \to \theta$ by $u_2 = 0$ and $u_1 \downarrow 0$. It follows that

$$\limsup_{u \to \theta} \frac{x_1 u_1 + x_2 u_2 - y_1 \frac{u_1^2-u_2^2}{\sqrt{u_1^2+u_2^2}} - y_2 \frac{2u_1 u_2}{\sqrt{u_1^2+u_2^2}}}{2\|u\|}$$

$$\geq \lim_{u_2=0 \text{ and } u_1 \downarrow 0} \frac{x_1 u_1 - y_1 u_1}{2|u_1|}$$

$$= \lim_{u_2=0 \text{ and } u_1 \downarrow 0} \left( \frac{x_1}{2} - \frac{y_1}{2} \right)$$

$$= \frac{x_1}{2} - \frac{y_1}{2} > 0.$$

This implies that

$$x_1 > y_1 \quad \Rightarrow \quad x \notin \widehat{D}^* f(\theta)(y).$$

Case 2. $-x_1 > y_1$. In this case, we take a special direction in the limit (5.3) for $u \to \theta$ by $u_2 = 0$ and $u_1 \uparrow 0$. We have

$$\operatorname*{limsup}_{u\to\theta}\frac{x_1u_1+x_2u_2-y_1\frac{u_1^2-u_2^2}{\sqrt{u_1^2+u_2^2}}-y_2\frac{2u_1u_2}{\sqrt{u_1^2+u_2^2}}}{2\|u\|}$$

$$\geq \lim_{u_2=0 \text{ and } u_1\uparrow 0}\frac{x_1u_1-y_1|u_1|}{2|u_1|}$$

$$= \lim_{u_2=0 \text{ and } u_1\uparrow 0}\left(\frac{-x_1}{2}-\frac{y_1}{2}\right)$$

$$= -\frac{x_1}{2}-\frac{y_1}{2} > 0.$$

This implies that

$$-x_1 > y_1 \quad \Rightarrow \quad x \notin \widehat{D}^*f(\theta)(y).$$

Case 3. $x_1 + x_2 > \sqrt{2}y_2$. In this case, we take a special direction in the limit (5.3) for $u \to \theta$ by $u_1 = u_2$ and $u_1 \downarrow 0$. We have

$$\operatorname*{limsup}_{u\to\theta}\frac{x_1u_1+x_2u_2-y_1\frac{u_1^2-u_2^2}{\sqrt{u_1^2+u_2^2}}-y_2\frac{2u_1u_2}{\sqrt{u_1^2+u_2^2}}}{2\|u\|}$$

$$\geq \lim_{u_1=u_2 \text{ and } u_1\downarrow 0}\frac{x_1u_1+x_2u_1-\frac{2y_2u_1}{\sqrt{2}}}{2\sqrt{2}u_1}$$

$$= \lim_{u_1=u_2 \text{ and } u_1\downarrow 0}\frac{x_1+x_2-\sqrt{2}y_2}{2\sqrt{2}}$$

$$= \frac{x_1+x_2-\sqrt{2}y_2}{2\sqrt{2}} > 0.$$

This implies that

$$x_1 + x_2 > \sqrt{2}y_2 \quad \Rightarrow \quad x \notin \widehat{D}^*f(\theta)(y).$$

Case 4. $-x_1 - x_2 > \sqrt{2}y_2$. In this case, we take a special direction in the limit (5.3) for $u \to \theta$ by $u_1 = u_2$ and $u_1 \uparrow 0$. We have

$$\operatorname*{limsup}_{u\to\theta}\frac{x_1u_1+x_2u_2-y_1\frac{u_1^2-u_2^2}{\sqrt{u_1^2+u_2^2}}-y_2\frac{2u_1u_2}{\sqrt{u_1^2+u_2^2}}}{2\|u\|}$$

$$\geq \lim_{u_1=u_2 \text{ and } u_1\uparrow 0}\frac{x_1u_1+x_2u_2-\frac{2y_2|u_1|}{\sqrt{2}}}{2\sqrt{2}|u_1|}$$

$$= \lim_{u_1=u_2 \text{ and } u_1\uparrow 0}\frac{-x_1-x_2-\sqrt{2}y_2}{2\sqrt{2}}$$

$$= \frac{-x_1-x_2-\sqrt{2}y_2}{2\sqrt{2}} > 0.$$

This implies that

$$-x_1 - x_2 > \sqrt{2}y_2 \implies x \notin \widehat{D}^*f(\theta)(y).$$

Case 5. $-x_1 + x_2 > -\sqrt{2}y_2$. In this case, we take a special direction in the limit (5.3) for $u \to \theta$ by $u_1 = -u_2$ and $u_2 \downarrow 0$. It follows that

$$\limsup_{u \to \theta} \frac{x_1 u_1 + x_2 u_2 - y_1 \frac{u_1^2 - u_2^2}{\sqrt{u_1^2 + u_2^2}} - y_2 \frac{2 u_1 u_2}{\sqrt{u_1^2 + u_2^2}}}{2 \|u\|}$$

$$\geq \lim_{u_1 = -u_2 \text{ and } u_2 \downarrow 0} \frac{x_1 u_1 + x_2 u_2 + \sqrt{2} y_2 u_2}{2\sqrt{2}|u_2|}$$

$$= \lim_{u_1 = -u_2 \text{ and } u_2 \downarrow 0} \frac{-x_1 u_2 + x_2 u_2 + \sqrt{2} y_2 u_2}{2\sqrt{2} u_2}$$

$$= \lim_{u_1 = -u_2 \text{ and } u_2 \downarrow 0} \left( \frac{-x_1}{2\sqrt{2}} + \frac{x_2}{2\sqrt{2}} + \frac{\sqrt{2} y_2}{2\sqrt{2}} \right)$$

$$= \frac{-x_1}{2\sqrt{2}} + \frac{x_2}{2\sqrt{2}} + \frac{\sqrt{2} y_2}{2\sqrt{2}} > 0.$$

This implies that

$$-x_1 + x_2 > -\sqrt{2}y_2 \implies x \notin \widehat{D}^*f(\theta)(y).$$

Case 6. $x_1 - x_2 > -\sqrt{2}y_2$. In this case, we take a special direction in the limit (5.3) for $u \to \theta$ by $u_1 = -u_2$ and $u_2 \uparrow 0$. It follows that

$$\limsup_{u \to \theta} \frac{x_1 u_1 + x_2 u_2 - y_1 \frac{u_1^2 - u_2^2}{\sqrt{u_1^2 + u_2^2}} - y_2 \frac{2 u_1 u_2}{\sqrt{u_1^2 + u_2^2}}}{2 \|u\|}$$

$$\geq \lim_{u_1 = -u_2 \text{ and } u_2 \uparrow 0} \frac{x_1 u_1 + x_2 u_2 + \sqrt{2} y_2 |u_2|}{2\sqrt{2}|u_2|}$$

$$= \lim_{u_1 = -u_2 \text{ and } u_2 \uparrow 0} \frac{-x_1 u_2 + x_2 u_2 + \sqrt{2} y_2 |u_2|}{2\sqrt{2}|u_2|}$$

$$= \lim_{u_1 = -u_2 \text{ and } u_2 \uparrow 0} \left( \frac{x_1}{2\sqrt{2}} - \frac{x_2}{2\sqrt{2}} + \frac{\sqrt{2} y_2}{2\sqrt{2}} \right)$$

$$= \frac{x_1}{2\sqrt{2}} - \frac{x_2}{2\sqrt{2}} + \frac{\sqrt{2} y_2}{2\sqrt{2}} > 0.$$

This implies that

$$x_1 - x_2 > -\sqrt{2}y_2 \implies x \notin \widehat{D}^*f(\theta)(y).$$

By summarizing the above 6 cases, we have that the following inequalities are necessary conditions for $x \in \widehat{D}^*f(\theta)(y)$.

(i) $x_1 \leq y_1$, (ii) $-x_1 \leq y_1$, (iii) $x_1 + x_2 \leq \sqrt{2}y_2$

(iv) $-x_1 - x_2 \leq \sqrt{2}y_2$, (v) $-x_1 + x_2 \leq -\sqrt{2}y_2$, (vi) $x_1 - x_2 \leq -\sqrt{2}y_2$

From the above inequalities (i–vi), we obtain more precise necessary conditions for $x \in \widehat{D}^*f(\theta)(y)$.

By (iii–vi), we have

(I) $\qquad\qquad\qquad x \in \widehat{D}^*f(\theta)(y) \quad \Longrightarrow \quad y_2 = 0.$

Substituting $y_2 = 0$ into (v–vi), we get

(II) $\qquad\qquad\qquad x \in \widehat{D}^*f(\theta)(y) \quad \Longrightarrow \quad x_1 - x_2 = 0$, that is $x_1 = x_2$.

Substituting $y_2 = 0$ and $x_1 = x_2$ into (iii–iv), we get

(III) $\qquad\qquad\qquad x \in \widehat{D}^*f(\theta)(y) \quad \Longrightarrow \quad x_1 = x_2 = 0.$

Substituting $x_1 = 0$ into (i) or (ii), we get

(IV) $\qquad\qquad\qquad x \in \widehat{D}^*f(\theta)(y) \quad \Longrightarrow \quad y_1 \geq 0.$

Under the necessary conditions (I–III), assume $y_1 > 0$. In this case, we take a special direction in the limit (5.3) for $u \to \theta$ by $u_1 = 0$ and $u_2 \uparrow 0$. We have

$$\limsup_{u \to \theta} \frac{x_1 u_1 + x_2 u_2 - y_1 \frac{u_1^2 - u_2^2}{\sqrt{u_1^2 + u_2^2}} - y_2 \frac{2u_1 u_2}{\sqrt{u_1^2 + u_2^2}}}{2\|u\|}$$

$$= \limsup_{u \to \theta} \frac{-y_1 \frac{u_1^2 - u_2^2}{\sqrt{u_1^2 + u_2^2}}}{2\|u\|}$$

$$\geq \lim_{u_1 = 0 \text{ and } u_2 \uparrow 0} \frac{y_1 |u_2|}{2|u_2|}$$

$$= \frac{y_1}{2} > 0.$$

This implies that

$$y_1 > 0 \quad \Longrightarrow \quad x \notin \widehat{D}^*f(\theta)(y).$$

By the necessary condition (IV), this implies that

(V) $\qquad\qquad\qquad x \in \widehat{D}^*f(\theta)(y) \quad \Longrightarrow \quad y_1 = 0.$

Under the necessary conditions (I–V), we obtain that, for $y = (y_1, y_2) \in \mathbb{R}^2$ and $x = (x_1, x_2) \in \mathbb{R}^2$,

$$x \in \widehat{D}^*f(\theta)(y) \quad \Longrightarrow \quad x = y = \theta.$$

This proposition is proved. $\qquad\square$

Next, we show that $f$ is Fréchet differentiable and Mordukhovich differentiable on $\mathbb{R}^2 \setminus \{\theta\}$.

**Theorem 5.3.** *Let $z = (z_1, z_2) \in \mathbb{R}^2 \setminus \{\theta\}$. Then $f$ is Fréchet differentiable and Mordukhovich differentiable at $z$ such that*

$$\nabla f(z) = \begin{pmatrix} \dfrac{(z_1^2 + 3z_2^2)z_1}{(z_1^2+z_2^2)\sqrt{z_1^2+z_2^2}} & \dfrac{2z_2^2 z_2}{(z_1^2+z_2^2)\sqrt{z_1^2+z_2^2}} \\ \dfrac{-(3z_1^2+z_2^2)z_2}{(z_1^2+z_2^2)\sqrt{z_1^2+z_2^2}} & \dfrac{2z_1^2 z_1}{(z_1^2+z_2^2)\sqrt{z_1^2+z_2^2}} \end{pmatrix} \quad and \quad \widehat{D}^* f(z) = \begin{pmatrix} \dfrac{(z_1^2 + 3z_2^2)z_1}{(z_1^2+z_2^2)\sqrt{z_1^2+z_2^2}} & \dfrac{-(3z_1^2+z_2^2)z_2}{(z_1^2+z_2^2)\sqrt{z_1^2+z_2^2}} \\ \dfrac{2z_2^2 z_2}{(z_1^2+z_2^2)\sqrt{z_1^2+z_2^2}} & \dfrac{2z_1^2 z_1}{(z_1^2+z_2^2)\sqrt{z_1^2+z_2^2}} \end{pmatrix}.$$

*However, more precisely speaking, for $x = (x_1, x_2), y = (y_1, y_2) \in \mathbb{R}^2$, if $x = \widehat{D}^* f(z)(y)$, then*

$$x_1 = y_1 \frac{(z_1^2 + 3z_2^2)z_1}{(z_1^2+z_2^2)\sqrt{z_1^2+z_2^2}} + y_2 \frac{2z_2^2 z_2}{(z_1^2+z_2^2)\sqrt{z_1^2+z_2^2}},$$

$$x_2 = y_1 \frac{(-3z_1^2 - z_2^2)z_2}{(z_1^2+z_2^2)\sqrt{z_1^2+z_2^2}} + y_2 \frac{2z_1^2 z_1}{(z_1^2+z_2^2)\sqrt{z_1^2+z_2^2}}.$$

*Proof.* Proof of (a). Let $f = (f_1, f_2)$ with

$$f_1(x_1, x_2) = \frac{x_1^2 - x_2^2}{\sqrt{x_1^2 + x_2^2}}, \text{ for } (x_1, x_2) \in \mathbb{R}^2 \setminus \{\theta\} \text{ with } f_1(\theta) = 0,$$

and

$$f_2(x_1, x_2) = \frac{2x_1 x_2}{\sqrt{x_1^2 + x_2^2}}, \text{ for } (x_1, x_2) \in \mathbb{R}^2 \setminus \{\theta\} \text{ with } f_2(\theta) = 0.$$

Let $z = (z_1, z_2) \in \mathbb{R}^2 \setminus \{\theta\}$. There is $r > 0$ such that for any $x = (x_1, x_2) \in \mathbb{R}^2$,

$$\|x - z\| < r \implies x \neq \theta.$$

This implies that both $f_1$ and $f_2$ are twice differentiable on the ball with radius $r > 0$ and centered at $z$. Then by Lemma 3.1 and Theorem 3.2, $f$ is Fréchet differentiable at $z$ and the Fréchet derivative of $f$ at $z$ is the following 2×2 matrix,

$$\nabla f(z) = \begin{pmatrix} \dfrac{(z_1^2 + 3z_2^2)z_1}{(z_1^2+z_2^2)\sqrt{z_1^2+z_2^2}} & \dfrac{2z_2^2 z_2}{(z_1^2+z_2^2)\sqrt{z_1^2+z_2^2}} \\ \dfrac{-(3z_1^2+z_2^2)z_2}{(z_1^2+z_2^2)\sqrt{z_1^2+z_2^2}} & \dfrac{2z_1^2 z_1}{(z_1^2+z_2^2)\sqrt{z_1^2+z_2^2}} \end{pmatrix}.$$

Proof of part (b). By Theorem 1.38 in [18] and Part (a), we have that $f$ is Mordukhovich differentiable at $z$ and it satisfies that

$$\widehat{D}^* f(z)(y) = (\nabla f(x))^*(y), \text{ for any } y \in \mathbb{R}^2.$$

By the representation of $\nabla f(x)$ in Part (a), this gives the solution of $\widehat{D}^* f(z)$ listed in Part (b). □

**Proposition 5.4.** *Let $z = (z_1, z_2) \in \mathbb{R}^2 \setminus \{\theta\}$ and $x = (x_1, x_2), y = (y_1, y_2) \in \mathbb{R}^2$, if $x = \widehat{D}^* f(z)(y)$, then $x$ and $y$ satisfy the following conditions.*

(i) $$x_1^2 + x_2^2 = y_1^2 + y_2^2 + 12\left(y_1 \frac{z_1 z_2}{z_1^2+z_2^2} - y_2 \frac{z_1^2-z_2^2}{2(z_1^2+z_2^2)}\right)^2.$$

(ii) $$\|x\| \geq \|y\|.$$

(iii) $$y_1 z_1 z_2 = y_2 \frac{z_1^2-z_2^2}{2} \implies \|x\| = \|y\|.$$

*Proof.* We only prove (i). It is because that Parts (ii) and (iii) follow from Part (i) immediately. By **Theorem 5.3** and $x \in \widehat{D}^* f(z)(y)$, we have

$$x_1 = y_1 \left(\frac{(z_1^2+3z_2^2)z_1}{(z_1^2+z_2^2)\sqrt{z_1^2+z_2^2}}\right) + y_2 \left(\frac{2z_2^2 z_2}{(z_1^2+z_2^2)\sqrt{z_1^2+z_2^2}}\right),$$

and

$$x_2 = y_1 \left(\frac{-(3z_1^2+z_2^2)z_2}{(z_1^2+z_2^2)\sqrt{z_1^2+z_2^2}}\right) + y_2 \left(\frac{2z_1^2 z_1}{(z_1^2+z_2^2)\sqrt{z_1^2+z_2^2}}\right).$$

Then, we calculate

$$x_1^2 + x_2^2 = \left(y_1 \left(\frac{(z_1^2+3z_2^2)z_1}{(z_1^2+z_2^2)\sqrt{z_1^2+z_2^2}}\right) + y_2 \left(\frac{2z_2^2 z_2}{(z_1^2+z_2^2)\sqrt{z_1^2+z_2^2}}\right)\right)^2 + \left(y_1 \left(\frac{-(3z_1^2+z_2^2)z_2}{(z_1^2+z_2^2)\sqrt{z_1^2+z_2^2}}\right) + y_2 \left(\frac{2z_1^2 z_1}{(z_1^2+z_2^2)\sqrt{z_1^2+z_2^2}}\right)\right)^2$$

$$= y_1^2 \left(\left(\frac{(z_1^2+3z_2^2)z_1}{(z_1^2+z_2^2)\sqrt{z_1^2+z_2^2}}\right)^2 + \left(\frac{-(3z_1^2+z_2^2)z_2}{(z_1^2+z_2^2)\sqrt{z_1^2+z_2^2}}\right)^2\right) + y_2^2 \left(\left(\frac{2z_2^2 z_2}{(z_1^2+z_2^2)\sqrt{z_1^2+z_2^2}}\right)^2 + \left(\frac{2z_1^2 z_1}{(z_1^2+z_2^2)\sqrt{z_1^2+z_2^2}}\right)^2\right)$$

$$+ 2y_1 y_2 \left(\left(\frac{-(3z_1^2+z_2^2)z_2}{(z_1^2+z_2^2)\sqrt{z_1^2+z_2^2}}\right)\left(\frac{2z_2^2 z_2}{(z_1^2+z_2^2)\sqrt{z_1^2+z_2^2}}\right) + \left(\frac{-(3z_1^2+z_2^2)z_2}{(z_1^2+z_2^2)\sqrt{z_1^2+z_2^2}}\right)\left(\frac{2z_1^2 z_1}{(z_1^2+z_2^2)\sqrt{z_1^2+z_2^2}}\right)\right)$$

$$= y_1^2 \left(\frac{(z_1^2-z_2^2)^2}{(z_1^2+z_2^2)^2} - \frac{4(z_1^2-z_2^2)^2}{(z_1^2+z_2^2)^2} + \frac{4(z_1^2+z_2^2)}{z_1^2+z_2^2}\right) + 2y_1 y_2 \left(\frac{2(z_1^2-z_2^2)z_1 z_2}{(z_1^2+z_2^2)^2} - \frac{4(z_1^2-z_2^2)z_1 z_2}{(z_1^2+z_2^2)^2} - \frac{4(z_1^2-z_2^2)z_1 z_2}{(z_1^2+z_2^2)^2}\right)$$

$$+ y_2^2 \left(4\left(\frac{z_1^4+2z_1^2 z_2^2+z_2^4}{(z_1^2+z_2^2)^2} - \frac{3z_1^2 z_2^2}{(z_1^2+z_2^2)^2}\right)\right)$$

$$= y_1^2 \left(\frac{z_1^4+14z_1^2 z_2^2+z_2^4}{(z_1^2+z_2^2)^2}\right) + 2y_1 y_2 \left(-\frac{6(z_1^2-z_2^2)z_1 z_2}{(z_1^2+z_2^2)^2}\right) + y_2^2 \left(4\left(\frac{z_1^4+2z_1^2 z_2^2+z_2^4}{(z_1^2+z_2^2)^2} - \frac{3z_1^2 z_2^2}{(z_1^2+z_2^2)^2}\right)\right)$$

$$= y_1^2 \left(1 + \frac{12z_1^2 z_2^2}{(z_1^2+z_2^2)^2}\right) + 2y_1 y_2 \left(-\frac{6(z_1^2-z_2^2)z_1 z_2}{(z_1^2+z_2^2)^2}\right) + y_2^2 \left(4\left(\frac{(z_1^2-z_2^2)^2}{(z_1^2+z_2^2)^2} + \frac{z_1^2 z_2^2}{(z_1^2+z_2^2)^2}\right)\right)$$

$$= y_1^2 + y_2^2 + y_1^2 \frac{12z_1^2 z_2^2}{(z_1^2+z_2^2)^2} - 4y_1 y_2 \left(\frac{3(z_1^2-z_2^2)z_1 z_2}{(z_1^2+z_2^2)^2}\right) + 4y_2^2 \left(\frac{(z_1^2-z_2^2)^2}{(z_1^2+z_2^2)^2} + \frac{z_1^2 z_2^2}{(z_1^2+z_2^2)^2} - \frac{1}{4}\right)$$

$$= y_1^2 + y_2^2 + y_1^2 \frac{12z_1^2 z_2^2}{(z_1^2+z_2^2)^2} - 4y_1 y_2 \left(\frac{3(z_1^2-z_2^2)z_1 z_2}{(z_1^2+z_2^2)^2}\right) + 4y_2^2 \left(\frac{4(z_1^2-z_2^2)^2 + 4z_1^2 z_2^2 - (z_1^2+z_2^2)^2}{4(z_1^2+z_2^2)^2}\right)$$

$$= y_1^2 + y_2^2 + y_1^2 \frac{12z_1^2 z_2^2}{(z_1^2+z_2^2)^2} - 4y_1 y_2 \left(\frac{3(z_1^2-z_2^2)z_1 z_2}{(z_1^2+z_2^2)^2}\right) + 4y_2^2 \left(\frac{3(z_1^4+z_2^4) - 6z_1^2 z_2^2}{4(z_1^2+z_2^2)^2}\right)$$

$$= y_1^2 + y_2^2 + 12\left(y_1^2 \frac{z_1^2 z_2^2}{(z_1^2+z_2^2)^2} - y_1 y_2 \frac{(z_1^2-z_2^2)z_1 z_2}{(z_1^2+z_2^2)^2} + y_2^2 \frac{(z_1^2-z_2^2)^2}{4(z_1^2+z_2^2)^2}\right)$$

$$= y_1^2 + y_2^2 + 12\left(y_1 \frac{z_1 z_2}{z_1^2+z_2^2} - y_2 \frac{z_1^2-z_2^2}{2(z_1^2+z_2^2)}\right)^2. \qquad \square$$

**Theorem 5.5.** *Let $f: \mathbb{R}^2 \to \mathbb{R}^2$ be defined by* (5.1). *Then, we have*

$$\hat{\alpha}(f, \bar{z}, f(\bar{z})) = 1, \text{ for any } \bar{z} \in \mathbb{R}^2.$$

*Proof.* For any $\bar{z}, \bar{w} \in \mathbb{R}^2$ with $\bar{w} = f(\bar{z})$ and $\eta > 0$, let $\mathbb{B}(\bar{z}, \eta)$ and $\mathbb{B}(\bar{w}, \eta)$ denote the closed ball in $\mathbb{R}^2$ with radius $\eta$ and centered at point $\bar{z}$ and $\bar{w}$, respectively. By (2.4) or (2.5), we have

$$\hat{\alpha}(f, \bar{z}, \bar{w}) = \sup_{\eta>0} \inf\{\|x\|: x \in \widehat{D}^* f(z)(y), z \in \mathbb{B}(\bar{z}, \eta), f(z) \in \mathbb{B}(\bar{w}, \eta), \|y\| = 1\}$$

$$= \sup_{\eta>0} \inf\{\|x\|: x \in \widehat{D}^* f(z)(y), z \in \mathbb{B}(\bar{z}, \eta), f(z) \in \mathbb{B}(\bar{w}, \eta), \|y\| = 1\}.$$

In case $\theta \in \mathbb{B}(\bar{z}, \eta)$ and $f(\theta) = \theta \in \mathbb{B}(\bar{w}, \eta)$, by Proposition 5.4, we have that

$$\widehat{D}^* f(\theta)(y) = \emptyset, \text{ for any } y \in \mathbb{R}^2 \text{ with } \|y\| = 1.$$

This implies that in the above limit set in, there are some point $z \in \mathbb{B}(\bar{z}, \eta)$ and $y \in \mathbb{R}^2$ with $\|y\| = 1$ such that $\widehat{D}^* f(z)(y) = \emptyset$. Since this Hilbert space $\mathbb{R}^2$ is Asplund and $f: \mathbb{R}^2 \to \mathbb{R}^2$ defined by (1.1) is continuous on $\mathbb{R}^2$, it yields that, for every $\eta > 0$, $f$ is Fréchet differentiable on a dense subset of $\mathbb{B}(\bar{z}, \eta)$. Then, by Theorem 1.38 in [18], $f$ is Mordukhovich differentiable on a dense subset of $\mathbb{B}(\bar{z}, \eta)$, for every $\eta > 0$, which implies that

$$\{\|x\|: x \in \widehat{D}^* f(z)(y), z \in \mathbb{B}(\bar{z}, \eta), f(z) \in \mathbb{B}(\bar{w}, \eta), \|y\| = 1\} \neq \emptyset, \text{ for every } \eta > 0.$$

However, if $\widehat{D}^* f(z)(y) \neq \emptyset$, for $z \in \mathbb{B}(\bar{z}, \eta)$ and $y \in \mathbb{R}^2$ with $\|y\| = 1$, then, for any $x \in \widehat{D}^* f(z)(y)$, by Theorem 4.6, we have $\|x\| \geq \|y\|$. One can show that, for any $z = (z_1, z_2) \in \mathbb{B}(\bar{z}, r)$, the following system of equations has solutions for $y_1, y_2$.

$$\begin{cases} y_1^2 + y_2^2 = 1, \\ y_1 z_1 z_2 = y_2 \frac{z_1^2 - z_2^2}{2}. \end{cases}$$

This implies that, for every $\eta > 0$, we have

$$\left\{\|x\|: x \in \widehat{D}^* f(z)(y), z \in \mathbb{B}(\bar{z}, \eta), f(z) \in \mathbb{B}(\bar{w}, \eta), \|y\| = 1, y_1 z_1 z_2 = y_2 \frac{z_1^2 - z_2^2}{2}\right\} \neq \emptyset.$$

Then, we calculate

$$\hat{\alpha}(f,\bar{z},\bar{w}) = \sup_{\eta>0}\inf\{\|x\|: x \in \widehat{D}^*f(z)(y), z \in \mathbb{B}(\bar{z},\eta), f(z) \in \mathbb{B}(\bar{w},\eta), \|y\| = 1\}$$

$$= \sup_{\eta>0}\inf\left\{\|x\|: x \in \widehat{D}^*f(z)(y), z \in \mathbb{B}(\bar{z},\eta), f(z) \in \mathbb{B}(\bar{w},\eta), \|y\| = 1, y_1 z_1 z_2 = y_2 \frac{z_1^2 - z_2^2}{2}\right\}$$

$$= \sup_{\eta>0}\inf\left\{\|y\|: x \in \widehat{D}^*f(z)(y), z \in \mathbb{B}(\bar{z},\eta), f(z) \in \mathbb{B}(\bar{w},\eta), \|y\| = 1, y_1 z_1 z_2 = y_2 \frac{z_1^2 - z_2^2}{2}\right\}$$

$$= \sup_{\eta>0}\{1\} = 1. \qquad \square$$

Notice that in Theorem 5.3, for any $z = (z_1, z_2) \in \mathbb{R}^2\setminus\{\theta\}$, the existence and solution of $\widehat{D}^*f(z)$ are proved by the Fréchet derivative $\nabla f(z)$ at $z$ and Theorem 1.38 in [18]. It is clearly to see that it is very difficult to solve $\widehat{D}^*f(z)$ without using the Fréchet derivative $\nabla f(z)$ at $z$. Next, we directly prove the results of $\widehat{D}^*f(z)$ in Part (b) of Theorem 5.3. This proof provides some ideas to directly calculate $\widehat{D}^*f(z)$ without using its Fréchet derivative $\nabla f(z)$ and Theorem 1.38 in [18].

**Proposition 5.6.** *Let* $z = (z_1, z_2) \in \mathbb{R}^2\setminus\{\theta\}$ *and* $x = (x_1, x_2), y = (y_1, y_2) \in \mathbb{R}^2$. *If* $x \in \widehat{D}^*f(z)(y)$, *then*

$$x_1 = y_1\left(\frac{-(z_1^2 - z_2^2)z_1}{(z_1^2+z_2^2)\sqrt{z_1^2+z_2^2}} + \frac{2z_1}{\sqrt{z_1^2+z_2^2}}\right) + y_2\left(\frac{-2z_1 z_2 z_1}{(z_1^2+z_2^2)\sqrt{z_1^2+z_2^2}} + \frac{2z_2}{\sqrt{z_1^2+z_2^2}}\right)$$

$$= y_1\frac{(z_1^2 + 3z_2^2)z_1}{(z_1^2+z_2^2)\sqrt{z_1^2+z_2^2}} + y_2\frac{2z_2^2 z_2}{(z_1^2+z_2^2)\sqrt{z_1^2+z_2^2}},$$

*and*

$$x_2 = y_1\left(\frac{-(z_1^2 - z_2^2)z_2}{(z_1^2+z_2^2)\sqrt{z_1^2+z_2^2}} + \frac{-2z_2}{\sqrt{z_1^2+z_2^2}}\right) + y_2\left(\frac{-2z_1 z_2 z_2}{(z_1^2+z_2^2)\sqrt{z_1^2+z_2^2}} + \frac{2z_1}{\sqrt{z_1^2+z_2^2}}\right)$$

$$= y_1\frac{(-3z_1^2 - z_2^2)z_2}{(z_1^2+z_2^2)\sqrt{z_1^2+z_2^2}} + y_2\frac{2z_1^2 z_1}{(z_1^2+z_2^2)\sqrt{z_1^2+z_2^2}}.$$

*Proof.* Let $y = (y_1, y_2) \in \mathbb{R}^2$ and $x = (x_1, x_2) \in \mathbb{R}^2$, in order to check if $x \in \widehat{D}^*f(z)(y)$ or not, we calculate the following limits.

$$\limsup_{u\to z}\frac{\langle x, u-z\rangle - \langle y, f(u)-f(z)\rangle}{\|u-z\| + \|f(u)-f(z)\|}$$

$$= \limsup_{u\to z}\frac{\langle(x_1,x_2),\ (u_1,u_2)-(z_1,z_2)\rangle - \langle(y_1,y_2),\ \left(\frac{u_1^2-u_2^2}{\sqrt{u_1^2+u_2^2}}, \frac{2u_1 u_2}{\sqrt{u_1^2+u_2^2}}\right) - \left(\frac{z_1^2-z_2^2}{\sqrt{z_1^2+z_2^2}}, \frac{2z_1 z_2}{\sqrt{z_1^2+z_2^2}}\right)\rangle}{\|u-z\| + \|f(u)-f(z)\|}$$

$$= \limsup_{u\to z}\frac{x_1(u_1-z_1)+x_2(u_2-z_2) - y_1\left(\frac{u_1^2-u_2^2}{\sqrt{u_1^2+u_2^2}} - \frac{z_1^2-z_2^2}{\sqrt{z_1^2+z_2^2}}\right) - y_2\left(\frac{2u_1 u_2}{\sqrt{u_1^2+u_2^2}} - \frac{2z_1 z_2}{\sqrt{z_1^2+z_2^2}}\right)}{\|u-z\| + \|f(u)-f(z)\|}$$

$$= \limsup_{u \to z} \frac{x_1(u_1-z_1)+x_2(u_2-z_2) - y_1\left(\frac{u_1^2-u_2^2}{\sqrt{u_1^2+u_2^2}} - \frac{z_1^2-z_2^2}{\sqrt{z_1^2+z_2^2}}\right) - y_2\left(\frac{2u_1 u_2}{\sqrt{u_1^2+u_2^2}} - \frac{2z_1 z_2}{\sqrt{z_1^2+z_2^2}}\right)}{\sqrt{(u_1-z_1)^2+(u_2-z_2)^2} + \sqrt{\left(\frac{u_1^2-u_2^2}{\sqrt{u_1^2+u_2^2}} - \frac{z_1^2-z_2^2}{\sqrt{z_1^2+z_2^2}}\right)^2 + \left(\frac{2u_1 u_2}{\sqrt{u_1^2+u_2^2}} - \frac{2z_1 z_2}{\sqrt{z_1^2+z_2^2}}\right)^2}}. \quad (5.4)$$

In the limit (5.4), we consider some special directions for $u \to z$.

(D1) We take a special direction in the limit (3.4) as $u = (z_1, z_2 + s)$ with $s \downarrow 0$. We have

$$\limsup_{u \to z} \frac{\langle x, u-z \rangle - \langle y, f(u)-f(z) \rangle}{\|u-z\| + \|f(u)-f(z)\|}$$

$$\geq \lim_{u = (z_1, z_2+s) \text{ with } s \downarrow 0} \frac{\langle x, u-z \rangle - \langle y, f(u)-f(z) \rangle}{\|u-z\| + \|f(u)-f(z)\|}$$

$$= \lim_{s \downarrow 0} \frac{sx_2 - y_1\left(\frac{z_1^2-(z_2+s)^2}{\sqrt{z_1^2+(z_2+s)^2}} - \frac{z_1^2-z_2^2}{\sqrt{z_1^2+z_2^2}}\right) - y_2\left(\frac{2z_1(z_2+s)}{\sqrt{z_1^2+(z_2+s)^2}} - \frac{2z_1 z_2}{\sqrt{z_1^2+z_2^2}}\right)}{\sqrt{(z_1-z_1)^2+(z_2+s-z_2)^2} + \sqrt{\left(\frac{z_1^2-(z_2+s)^2}{\sqrt{z_1^2+(z_2+s)^2}} - \frac{z_1^2-z_2^2}{\sqrt{z_1^2+z_2^2}}\right)^2 + \left(\frac{2z_1(z_2+s)}{\sqrt{z_1^2+(z_2+s)^2}} - \frac{2z_1 z_2}{\sqrt{z_1^2+z_2^2}}\right)^2}}$$

$$= \lim_{s \downarrow 0} \frac{sx_2 - y_1\left(\frac{z_1^2-(z_2+s)^2}{\sqrt{z_1^2+(z_2+s)^2}} - \frac{z_1^2-z_2^2}{\sqrt{z_1^2+z_2^2}}\right) - y_2\left(\frac{2z_1(z_2+s)}{\sqrt{z_1^2+(z_2+s)^2}} - \frac{2z_1 z_2}{\sqrt{z_1^2+z_2^2}}\right)}{s + \sqrt{\left(\frac{z_1^2-(z_2+s)^2}{\sqrt{z_1^2+(z_2+s)^2}} - \frac{z_1^2-z_2^2}{\sqrt{z_1^2+z_2^2}}\right)^2 + \left(\frac{2z_1(z_2+s)}{\sqrt{z_1^2+(z_2+s)^2}} - \frac{2z_1 z_2}{\sqrt{z_1^2+z_2^2}}\right)^2}}. \quad (5.5)$$

We calculate the two big terms in (5.5).

$$\frac{z_1^2-(z_2+s)^2}{\sqrt{z_1^2+(z_2+s)^2}} - \frac{z_1^2-z_2^2}{\sqrt{z_1^2+z_2^2}}$$

$$= \frac{(z_1^2-z_2^2)-2sz_2-s^2}{\sqrt{(z_1^2+z_2^2)+2sz_2+s^2}} - \frac{z_1^2-z_2^2}{\sqrt{z_1^2+z_2^2}}$$

$$= \frac{z_1^2-z_2^2}{\sqrt{(z_1^2+z_2^2)+2sz_2+s^2}} - \frac{z_1^2-z_2^2}{\sqrt{z_1^2+z_2^2}} + \frac{-2sz_2-s^2}{\sqrt{(z_1^2+z_2^2)+2sz_2+s^2}}$$

$$= \frac{(z_1^2-z_2^2)\left(\sqrt{z_1^2+z_2^2} - \sqrt{(z_1^2+z_2^2)+2sz_2+s^2}\right)}{\sqrt{(z_1^2+z_2^2)+2sz_2+s^2}\sqrt{z_1^2+z_2^2}} + \frac{-2sz_2-s^2}{\sqrt{(z_1^2+z_2^2)+2sz_2+s^2}}$$

$$= \frac{(z_1^2-z_2^2)\left(\left(\sqrt{z_1^2+z_2^2}\right)^2 - \left(\sqrt{(z_1^2+z_2^2)+2sz_2+s^2}\right)^2\right)}{\sqrt{(z_1^2+z_2^2)+2sz_2+s^2}\sqrt{z_1^2+z_2^2}\left(\sqrt{z_1^2+z_2^2} + \sqrt{(z_1^2+z_2^2)+2sz_2+s^2}\right)} + \frac{-2sz_2-s^2}{\sqrt{(z_1^2+z_2^2)+2sz_2+s^2}}$$

$$= \frac{(z_1^2-z_2^2)\Big(\big((z_1^2+z_2^2)\big)-\big((z_1^2+z_2^2)+2sz_2+s^2\big)\Big)}{\sqrt{(z_1^2+z_2^2)+2sz_2+s^2}\sqrt{z_1^2+z_2^2}\Big(\sqrt{z_1^2+z_2^2}+\sqrt{(z_1^2+z_2^2)+2sz_2+s^2}\Big)} + \frac{-2sz_2-s^2}{\sqrt{(z_1^2+z_2^2)+2sz_2+s^2}}$$

$$= \frac{(z_1^2-z_2^2)(-2sz_2-s^2)}{\sqrt{(z_1^2+z_2^2)+2sz_2+s^2}\sqrt{z_1^2+z_2^2}\Big(\sqrt{z_1^2+z_2^2}+\sqrt{(z_1^2+z_2^2)+2sz_2+s^2}\Big)} + \frac{-2sz_2-s^2}{\sqrt{(z_1^2+z_2^2)+2sz_2+s^2}}$$

$$= s\left(\frac{-(z_1^2-z_2^2)(2z_2+s)}{\sqrt{(z_1^2+z_2^2)+2sz_2+s^2}\sqrt{z_1^2+z_2^2}\Big(\sqrt{z_1^2+z_2^2}+\sqrt{(z_1^2+z_2^2)+2sz_2+s^2}\Big)} + \frac{-2z_2-s}{\sqrt{(z_1^2+z_2^2)+2sz_2+s^2}}\right)$$

$$\equiv sA(s). \tag{5.6}$$

Now we calculate other big term in (5.5).

$$\frac{2z_1(z_2+s)}{\sqrt{(z_1)^2+(z_2+s)^2}} - \frac{2z_1 z_2}{\sqrt{z_1^2+z_2^2}}$$

$$= \frac{2z_1 z_2+2sz_1}{\sqrt{(z_1^2+z_2^2)+2sz_2+s^2}} - \frac{2z_1 z_2}{\sqrt{z_1^2+z_2^2}}$$

$$= \frac{2z_1 z_2}{\sqrt{(z_1^2+z_2^2)+2sz_2+s^2}} - \frac{2z_1 z_2}{\sqrt{z_1^2+z_2^2}} + \frac{2sz_1}{\sqrt{(z_1^2+z_2^2)+2sz_2+s^2}}$$

$$= \frac{2z_1 z_2\Big(\sqrt{z_1^2+z_2^2} - \sqrt{(z_1^2+z_2^2)+2s(z_1+z_2)+2s^2}\Big)}{\sqrt{(z_1^2+z_2^2)+2sz_2+s^2}\sqrt{z_1^2+z_2^2}} + \frac{2sz_1}{\sqrt{(z_1^2+z_2^2)+2sz_2+s^2}}$$

$$= \frac{2z_1 z_2\Big((z_1^2+z_2^2)-\big((z_1^2+z_2^2)+2sz_2+s^2\big)\Big)}{\sqrt{(z_1^2+z_2^2)+2sz_2+s^2}\sqrt{z_1^2+z_2^2}\Big(\sqrt{z_1^2+z_2^2}+\sqrt{(z_1^2+z_2^2)+2sz_2+s^2}\Big)} + \frac{2sz_1}{\sqrt{(z_1^2+z_2^2)+2sz_2+s^2}}$$

$$= \frac{2z_1 z_2(-2sz_2-s^2))}{\sqrt{(z_1^2+z_2^2)+2sz_2+s^2}\sqrt{z_1^2+z_2^2}\Big(\sqrt{z_1^2+z_2^2}+\sqrt{(z_1^2+z_2^2)+2sz_2+s^2}\Big)} + \frac{2sz_1}{\sqrt{(z_1^2+z_2^2)+2sz_2+s^2}}$$

$$= s\left(\frac{-2z_1 z_2(2z_2+s)}{\sqrt{(z_1^2+z_2^2)+2sz_2+s^2}\sqrt{z_1^2+z_2^2}\Big(\sqrt{z_1^2+z_2^2}+\sqrt{(z_1^2+z_2^2)+2sz_2+s^2}\Big)} + \frac{2z_1}{\sqrt{(z_1^2+z_2^2)+2sz_2+s^2}}\right)$$

$$= sB(s). \tag{5.7}$$

Substituting (5.6) and (5.7) into the limit (5.5), we have

$$\limsup_{u\to z} \frac{\langle x,\ u-z\rangle - \langle y,\ f(u)-f(z)\rangle}{\|u-z\|+\|f(u)-f(z)\|}$$

$$\geq \lim_{s\downarrow 0} \frac{sx_2 - y_1 sA(s) - y_2 sB(s)}{s+\sqrt{(sA(s))^2+(sB(s))^2}}$$

$$= \lim_{s \downarrow 0} \frac{sx_2 - y_1 sA(s) - y_2 sB(s)}{s + s\sqrt{(A(s))^2 + (B(s))^2}}$$

$$= \lim_{s \downarrow 0} \frac{x_2 - y_1(A(s)) - y_2(B(s))}{1 + \sqrt{(A(s))^2 + (B(s))^2}}. \tag{5.8}$$

In the limit (5.8), by notations (5.6) and (5.7), we have

$$\lim_{s \downarrow 0} A(s) = \lim_{s \downarrow 0} \left( \frac{-(z_1^2 - z_2^2)(2z_2 + s)}{\sqrt{(z_1^2 + z_2^2) + 2sz_2 + s^2} \sqrt{z_1^2 + z_2^2} \left( \sqrt{z_1^2 + z_2^2} + \sqrt{(z_1^2 + z_2^2) + 2sz_2 + s^2} \right)} + \frac{-2z_2 - s}{\sqrt{(z_1^2 + z_2^2) + 2sz_2 + s^2}} \right)$$

$$= \frac{-(z_1^2 - z_2^2) 2z_2}{\sqrt{z_1^2 + z_2^2} \sqrt{z_1^2 + z_2^2} \left( \sqrt{z_1^2 + z_2^2} + \sqrt{z_1^2 + z_2^2} \right)} + \frac{-2z_2}{\sqrt{z_1^2 + z_2^2}}$$

$$= \frac{-(z_1^2 - z_2^2) z_2}{(z_1^2 + z_2^2) \sqrt{z_1^2 + z_2^2}} + \frac{-2z_2}{\sqrt{z_1^2 + z_2^2}}. \tag{5.9}$$

We also have

$$\lim_{s \downarrow 0} B(s) = \lim_{s \downarrow 0} \left( \frac{-2z_1 z_2 (2z_2 + s)}{\sqrt{(z_1^2 + z_2^2) + 2sz_2 + s^2} \sqrt{z_1^2 + z_2^2} \left( \sqrt{z_1^2 + z_2^2} + \sqrt{(z_1^2 + z_2^2) + 2sz_2 + s^2} \right)} + \frac{2z_1}{\sqrt{(z_1^2 + z_2^2) + 2sz_2 + s^2}} \right)$$

$$= \frac{-4z_1 z_2 z_2}{\sqrt{z_1^2 + z_2^2} \sqrt{z_1^2 + z_2^2} \left( \sqrt{z_1^2 + z_2^2} + \sqrt{z_1^2 + z_2^2} \right)} + \frac{2z_1}{\sqrt{z_1^2 + z_2^2}}$$

$$= \frac{-2z_1 z_2 z_2}{(z_1^2 + z_2^2) \sqrt{z_1^2 + z_2^2}} + \frac{2z_1}{\sqrt{z_1^2 + z_2^2}}. \tag{5.10}$$

Substituting (5.9) and (5.10) into (5.8), we have

$$\limsup_{u \to z} \frac{\langle x, u-z \rangle - \langle y, f(u) - f(z) \rangle}{\|u-z\| + \|f(u) - f(z)\|}$$

$$\geq \lim_{s \downarrow 0} \frac{x_2 - y_1(E(s)) - y_2(F(s))}{1 + \sqrt{(E(s))^2 + (F(s))^2}}$$

$$= \frac{x_2 - y_1 \left( \frac{-(z_1^2 - z_2^2) z_2}{(z_1^2 + z_2^2)\sqrt{z_1^2 + z_2^2}} + \frac{-2z_2}{\sqrt{z_1^2 + z_2^2}} \right) - y_2 \left( \frac{-2z_1 z_2 z_2}{(z_1^2 + z_2^2)\sqrt{z_1^2 + z_2^2}} + \frac{2z_1}{\sqrt{z_1^2 + z_2^2}} \right)}{1 + \sqrt{\left( \frac{-(z_1^2 - z_2^2) z_2}{(z_1^2 + z_2^2)\sqrt{z_1^2 + z_2^2}} + \frac{-2z_2}{\sqrt{z_1^2 + z_2^2}} \right)^2 + \left( \frac{-2z_1 z_2 z_2}{(z_1^2 + z_2^2)\sqrt{z_1^2 + z_2^2}} + \frac{2z_1}{\sqrt{z_1^2 + z_2^2}} \right)^2}}.$$

This implies that

$$x_2 - y_1 \left( \frac{-(z_1^2 - z_2^2) z_2}{(z_1^2 + z_2^2)\sqrt{z_1^2 + z_2^2}} + \frac{-2z_2}{\sqrt{z_1^2 + z_2^2}} \right) - y_2 \left( \frac{-2z_1 z_2 z_2}{(z_1^2 + z_2^2)\sqrt{z_1^2 + z_2^2}} + \frac{2z_1}{\sqrt{z_1^2 + z_2^2}} \right) > 0 \Longrightarrow x \notin \widehat{D}^* f(z)(y).$$

This reduces to that

(i) $x \in \widehat{D}^* f(z)(y) \implies x_2 - y_1 \left( \dfrac{-(z_1^2 - z_2^2) z_2}{(z_1^2 + z_2^2)\sqrt{z_1^2 + z_2^2}} + \dfrac{-2z_2}{\sqrt{z_1^2 + z_2^2}} \right) - y_2 \left( \dfrac{-2 z_1 z_2 z_2}{(z_1^2 + z_2^2)\sqrt{z_1^2 + z_2^2}} + \dfrac{2 z_1}{\sqrt{z_1^2 + z_2^2}} \right) \leq 0.$

(D2) We take a new special direction in the limit (3.4) $u = (z_1, z_2 - s)$ with $s \downarrow 0$ That is $s > 0$. With respect to this direction of the limit in (3.4), similarly, to the proof of (ii), we can show that

(j) $x \in \widehat{D}^* f(z)(y) \implies -x_2 - y_1 \left( \dfrac{(z_1^2 - z_2^2) z_2}{(z_1^2 + z_2^2)\sqrt{z_1^2 + z_2^2}} + \dfrac{2 z_2}{\sqrt{z_1^2 + z_2^2}} \right) - y_2 \left( \dfrac{2 z_1 z_2 z_2}{(z_1^2 + z_2^2)\sqrt{z_1^2 + z_2^2}} + \dfrac{-2 z_1}{\sqrt{z_1^2 + z_2^2}} \right) \leq 0.$

Combining (i) and (j), we obtain that

(I) $\quad x \in \widehat{D}^* f(z)(y) \implies x_2 = y_1 \left( \dfrac{-(z_1^2 - z_2^2) z_2}{(z_1^2 + z_2^2)\sqrt{z_1^2 + z_2^2}} + \dfrac{-2 z_2}{\sqrt{z_1^2 + z_2^2}} \right) + y_2 \left( \dfrac{-2 z_1 z_2 z_2}{(z_1^2 + z_2^2)\sqrt{z_1^2 + z_2^2}} + \dfrac{2 z_1}{\sqrt{z_1^2 + z_2^2}} \right).$

We take new special directions in the limit (5.4):

**(D3)** $u = (z_1 + s, z_2)$ with $s \downarrow 0$ That is $s > 0$ and **(D4)** $u = (z_1 - s, z_2)$ with $s \downarrow 0$ That is $s > 0$. Similarly, to the proof of (I), we can show that

(II) $\quad x \in \widehat{D}^* f(z)(y) \implies x_1 = y_1 \left( \dfrac{-(z_1^2 - z_2^2) z_1}{(z_1^2 + z_2^2)\sqrt{z_1^2 + z_2^2}} + \dfrac{2 z_1}{\sqrt{z_1^2 + z_2^2}} \right) + y_2 \left( \dfrac{-2 z_1 z_2 z_1}{(z_1^2 + z_2^2)\sqrt{z_1^2 + z_2^2}} + \dfrac{2 z_2}{\sqrt{z_1^2 + z_2^2}} \right).$

By (I) and (II), we have that $x \in \widehat{D}^* f(z)(y)$ implies that

$$x_1 = y_1 \left( \dfrac{-(z_1^2 - z_2^2) z_1}{(z_1^2 + z_2^2)\sqrt{z_1^2 + z_2^2}} + \dfrac{2 z_1}{\sqrt{z_1^2 + z_2^2}} \right) + y_2 \left( \dfrac{-2 z_1 z_2 z_1}{(z_1^2 + z_2^2)\sqrt{z_1^2 + z_2^2}} + \dfrac{2 z_2}{\sqrt{z_1^2 + z_2^2}} \right),$$

$$x_2 = y_1 \left( \dfrac{-(z_1^2 - z_2^2) z_2}{(z_1^2 + z_2^2)\sqrt{z_1^2 + z_2^2}} + \dfrac{-2 z_2}{\sqrt{z_1^2 + z_2^2}} \right) + y_2 \left( \dfrac{-2 z_1 z_2 z_2}{(z_1^2 + z_2^2)\sqrt{z_1^2 + z_2^2}} + \dfrac{2 z_1}{\sqrt{z_1^2 + z_2^2}} \right). \qquad \square$$

Next, we extend the mapping $f$ from $\mathbb{R}^2$ to $\mathbb{R}^2$ defined in (5.1) to a mapping $g: \mathbb{R}^4 \to \mathbb{R}^4$ as follows,

$$g((x_1, x_2, x_3, x_4)) = \left( \dfrac{x_1^2 - x_2^2}{\sqrt{x_1^2 + x_2^2}}, \dfrac{2 x_1 x_2}{\sqrt{x_1^2 + x_2^2}}, \dfrac{x_3^2 - x_4^2}{\sqrt{x_3^2 + x_4^2}}, \dfrac{2 x_3 x_4}{\sqrt{x_3^2 + x_4^2}} \right), \text{ for } (x_1, x_2, x_3, x_4) \in \mathbb{R}^4. \qquad (5.11)$$

Here, for $i = 1, 3$, if $x_i = x_{i+1} = 0$, then $\left( \dfrac{x_i^2 - x_{i+1}^2}{\sqrt{x_i^2 + x_{i+1}^2}}, \dfrac{2 x_i x_{i+1}}{\sqrt{x_i^2 + x_{i+1}^2}} \right)$ is defined to be $(0, 0)$.

Similar, to the norm preserving property of $f$ studied in the previous sections, the mapping $g: \mathbb{R}^4 \to \mathbb{R}^4$ is also norm preserving satisfying,

$$\|g(x)\| = \|x\|, \text{ for any } x \in \mathbb{R}^4.$$

**Lemma 5.7.** *Let* $z = (z_1, z_2, z_3, z_4) \in \mathbb{R}^4$. *If $z$ satisfies* $(z_1^2 + z_2^2)(z_3^2 + z_4^2) = 0$, *then $g$ is not Fréchet differentiable at $z$.*

*Proof.* The proof of this lemma is similar to the proof of Lemma 5.1. But it is more complicated. So, we prove it here. Let $z = (z_1, z_2, z_3, z_4) \in \mathbb{R}^4$ with $(z_1^2 + z_2^2)(z_3^2 + z_4^2) = 0$.

Suppose that $z_1^2 + z_2^2 = 0$. Assume, by the way of contradiction, that $g$ is Fréchet differentiable at $z$. Then $\nabla g(z) \colon \mathbb{R}^4 \to \mathbb{R}^4$ must be precisely represented by a real $4 \times 4$ matrix as below.

$$\nabla g(z) = \begin{pmatrix} a_{11} & a_{12} & a_{13} & a_{14} \\ a_{21} & a_{22} & a_{23} & a_{24} \\ a_{31} & a_{32} & a_{33} & a_{34} \\ a_{41} & a_{42} & a_{43} & a_{44} \end{pmatrix} \equiv A.$$

By the assumption that $z_1^2 + z_2^2 = 0$, that is, $z = (0, 0, z_3, z_4)$, we have

$$\lim_{u=(u_1,u_2,u_3,u_4) \to \theta} \frac{g(u) - g(z) - \nabla g(z)(u-z)}{\|u-z\|}$$

$$= \lim_{u \to \theta} \frac{\left( \frac{u_1^2 - u_2^2}{\sqrt{u_1^2 + u_2^2}} - 0, \frac{2u_1 u_2}{\sqrt{u_1^2 + u_2^2}} - 0, \frac{u_3^2 - u_4^2}{\sqrt{u_3^2 + u_4^2}} - \frac{z_3^2 - z_4^2}{\sqrt{z_3^2 + z_4^2}}, \frac{2u_3 u_4}{\sqrt{u_3^2 + u_4^2}} - \frac{2z_3 z_4}{\sqrt{z_3^2 + z_4^2}} \right) - (u-z)A}{\|u-z\|}. \tag{5.12}$$

Assume $a_{11} \neq 0$. In this case, we take a special direction in the limit (5.12) for $u \to \theta$ by $u_1 = -a_{11} t$ with $t \downarrow 0$, $u_2 = 0$, $u_3 = z_3$ and $u_4 = z_4$. By the assumption that $\nabla g(z) = A$, it follows that

$$\lim_{u=(u_1,u_2,u_3,u_4) \to \theta} \frac{g(u) - g(z) - \nabla g(z)(u-z)}{\|u-z\|}$$

$$= \lim_{t \downarrow 0} \frac{\left( \frac{a_{11}^2 t^2}{\sqrt{a_{11}^2 t^2}}, 0, 0, 0 \right) - (-a_{11}t, 0, 0, 0)A}{\|(-a_{11}t, 0, 0, 0)\|}$$

$$= \lim_{t \downarrow 0} \frac{(|a_{11}|t, 0, 0, 0) + (a_{11}^2 t, a_{11}a_{12}t, a_{11}a_{13}t, a_{11}a_{14}t)}{|a_{11}|t}$$

$$= (1 + |a_{11}|, (\text{sign}\, a_{11})a_{12}, (\text{sign}\, a_{11})a_{13}, (\text{sign}\, a_{11})a_{14})$$

$$\neq \theta.$$

This contradicts the assumption that $\nabla g(z) = A$ (which implies that $a_{11} = 0$).

Next, suppose $a_{11} = 0$. In this case, we take a special direction in the limit (5.12) for $u \to \theta$ by $u_1 \downarrow 0$, $u_2 = 0$, $u_3 = z_3$ and $u_4 = z_4$. It follows that

$$\lim_{u=(u_1,u_2,u_3,u_4) \to \theta} \frac{g(u) - g(z) - \nabla g(z)(u-z)}{\|u-z\|}$$

$$= \lim_{u_1 \downarrow 0} \frac{\left(\frac{u_1^2}{\sqrt{u_1^2}},\ 0,\ 0,\ 0\right) - (u_1,\ 0,\ 0,\ 0)A}{\|(u_1,0,0,0)\|}$$

$$= \lim_{u_1 \downarrow 0} \frac{(u_1,\ 0,\ 0,\ 0) + (0,\ u_1 a_{12},\ u_1 a_{13},\ u_1 a_{14})}{u_1}$$

$$= (1,\ a_{12},\ u_1 a_{13},\ u_1 a_{14})$$

$$\neq \theta.$$

This contradicts the assumption that $\nabla g(z) = A$. Next, we suppose that $z_3^2 + z_4^2 = 0$. In this case, $z = (z_1, z_2, 0, 0)$, we have

$$\lim_{u=(u_1,u_2,u_3,u_4) \to \theta} \frac{g(u) - g(z) - \nabla g(z)(u-z)}{\|u-z\|}$$

$$= \lim_{u \to \theta} \frac{\left(\frac{u_1^2 - u_2^2}{\sqrt{u_1^2+u_2^2}} - \frac{z_1^2 - z_2^2}{\sqrt{z_1^2+z_2^2}},\ \frac{2u_1 u_2}{\sqrt{u_1^2+u_2^2}} - \frac{2z_1 z_2}{\sqrt{z_1^2+z_2^2}},\ \frac{u_3^2 - u_4^2}{\sqrt{u_3^2+u_4^2}} - 0,\ \frac{2u_3 u_4}{\sqrt{u_3^2+u_4^2}} - 0 \right) - (u-z)A}{\|u-z\|}. \tag{5.13}$$

Assume $a_{33} \neq 0$. In this case, we take a special direction in the limit (5.12) for $u \to \theta$ by $u_3 = -a_{33} t$ with $t \downarrow 0$, $u_4 = 0$, $u_2 = z_2$ and $u_1 = z_1$. By the assumption that $\nabla g(z) = A$, it follows that

$$\lim_{u=(u_1,u_2,u_3,u_4) \to \theta} \frac{g(u) - g(z) - \nabla g(z)(u-z)}{\|u-z\|}$$

$$= \lim_{t \downarrow 0} \frac{\left(0,\ 0,\ \frac{a_{33}^2 t^2}{\sqrt{a_{33}^2 t^2}},\ 0\right) - (0,\ 0,\ -a_{33}t,\ 0)A}{\|(0,0,-a_{33}t,0)\|}$$

$$= \lim_{t \downarrow 0} \frac{(0,\ 0,\ |a_{33}|t,\ 0) + (a_{31} a_{33} t,\ a_{32} a_{33} t,\ a_{33}^2 t,\ a_{34} a_{33} t)}{|a_{33}|t}$$

$$= ((\operatorname{sign} a_{33})a_{31},\ (\operatorname{sign} a_{32})a_{32},\ 1 + |a_{33}|,\ (\operatorname{sign} a_{11})a_{34})$$

$$\neq \theta.$$

This contradicts the assumption that $\nabla g(z) = A$. Next, suppose $a_{33} = 0$. In this case, we take a special direction in the limit (5.13) for $u \to \theta$ by $u_3 \downarrow 0$, $u_4 = 0$, $u_2 = z_2$ and $u_1 = z_1$. It follows that

$$\lim_{u=(u_1,u_2,u_3,u_4) \to \theta} \frac{g(u) - g(z) - \nabla g(z)(u-z)}{\|u-z\|}$$

$$= \lim_{u_3 \downarrow 0} \frac{\left(0,\ 0,\ \frac{u_3^2}{\sqrt{u_3^2}},\ 0\right) - (0,\ 0,\ u_3,\ 0)A}{\|(0,0,u_3,0)\|}$$

$$= \lim_{u_3 \downarrow 0} \frac{(0,\ 0,\ u_3,\ 0) + (u_3 a_{31},\ u_3 a_{32},\ u_3 0,\ u_3 a_{34})}{u_3}$$

$$= (a_{31},\ a_{32},\ 1, a_{34})$$

$\neq \theta.$

This contradicts the assumption that $\nabla g(z) = A$. It completes the proof. □

**Corollary 5.8.** $g$ is not Fréchet differentiable at $\theta$.

*Proof.* In Lemma 5.7, let $z = (z_1, z_2, z_3, z_4) \in \mathbb{R}^4$ with $(z_1^2 + z_2^2) = 0$ and $(z_3^2 + z_4^2) = 0$. □

**Proposition 5.9.** *The mapping* $g \colon \mathbb{R}^4 \to \mathbb{R}^4$ *defined by* (5.11) *satisfies that*

(a) $\widehat{D}^*g(\theta)(\theta) = \{\theta\}$;
(b) $\widehat{D}^*g(\theta)(y) = \emptyset$, *for any* $y \in \mathbb{R}^4 \setminus \{\theta\}$.

*Proof.* The proof of this proposition is similar to the proof of Propositions 5.2. It is omitted here. □

**Theorem 5.10.** *Let* $z = (z_1, z_2, z_3, z_4) \in \mathbb{R}^4 \setminus \{\theta\}$ *with* $(z_1^2 + z_2^2)(z_3^2 + z_4^2) \neq 0$. *Then we have that*

(a) $g$ *is Fréchet differentiable and Mordukhovich differentiable at* $z$ *such that*

$$\nabla g(z) = \begin{pmatrix} \frac{(z_1^2+3z_2^2)z_1}{(z_1^2+z_2^2)\sqrt{z_1^2+z_2^2}} & \frac{2z_2^2 z_2}{(z_1^2+z_2^2)\sqrt{z_1^2+z_2^2}} & 0 & 0 \\ \frac{-(3z_1^2+z_2^2)z_2}{(z_1^2+z_2^2)\sqrt{z_1^2+z_2^2}} & \frac{2z_1^2 z_1}{(z_1^2+z_2^2)\sqrt{z_1^2+z_2^2}} & 0 & 0 \\ 0 & 0 & \frac{(z_3^2+3z_4^2)z_3}{(z_3^2+z_4^2)\sqrt{z_3^2+z_4^2}} & \frac{2z_4^2 z_4}{(z_3^2+z_4^2)\sqrt{z_3^2+z_4^2}} \\ 0 & 0 & \frac{-(3z_3^2+z_4^2)z_4}{(z_3^2+z_4^2)\sqrt{z_3^2+z_4^2}} & \frac{2z_3^2 z_3}{(z_3^2+z_4^2)\sqrt{z_3^2+z_4^2}} \end{pmatrix}.$$

*and* 

$$\widehat{D}^*g(z) = \begin{pmatrix} \frac{(z_1^2+3z_2^2)z_1}{(z_1^2+z_2^2)\sqrt{z_1^2+z_2^2}} & \frac{-(3z_1^2+z_2^2)z_2}{(z_1^2+z_2^2)\sqrt{z_1^2+z_2^2}} & 0 & 0 \\ \frac{2z_2^2 z_2}{(z_1^2+z_2^2)\sqrt{z_1^2+z_2^2}} & \frac{2z_1^2 z_1}{(z_1^2+z_2^2)\sqrt{z_1^2+z_2^2}} & 0 & 0 \\ 0 & 0 & \frac{(z_3^2+3z_4^2)z_3}{(z_3^2+z_4^2)\sqrt{z_3^2+z_4^2}} & \frac{-(3z_3^2+z_4^2)z_4}{(z_3^2+z_4^2)\sqrt{z_3^2+z_4^2}} \\ 0 & 0 & \frac{2z_4^2 z_4}{(z_3^2+z_4^2)\sqrt{z_3^2+z_4^2}} & \frac{2z_3^2 z_3}{(z_3^2+z_4^2)\sqrt{z_3^2+z_4^2}} \end{pmatrix}.$$

(*More precisely speaking*) *Let* $x = (x_1, x_2, x_3, x_4), y = (y_1, y_2, y_3, y_4) \in \mathbb{R}^4$. *If* $x = \widehat{D}^*g(z)(y)$, *then*

$$x_1 = y_1 \frac{(z_1^2+3z_2^2)z_1}{(z_1^2+z_2^2)\sqrt{z_1^2+z_2^2}} + y_2 \frac{2z_2^2 z_2}{(z_1^2+z_2^2)\sqrt{z_1^2+z_2^2}}, \quad x_2 = y_1 \frac{(-3z_1^2-z_2^2)z_2}{(z_1^2+z_2^2)\sqrt{z_1^2+z_2^2}} + y_2 \frac{2z_1^2 z_1}{(z_1^2+z_2^2)\sqrt{z_1^2+z_2^2}},$$

$$x_3 = y_3 \frac{(z_3^2+3z_4^2)z_3}{(z_3^2+z_4^2)\sqrt{z_3^2+z_4^2}} + y_4 \frac{2z_4^2 z_4}{(z_3^2+z_4^2)\sqrt{z_3^2+z_4^2}}, \quad x_4 = y_3 \frac{-(3z_3^2+z_4^2)z_4}{(z_3^2+z_4^2)\sqrt{z_3^2+z_4^2}} + y_4 \frac{2z_3^2 z_3}{(z_3^2+z_4^2)\sqrt{z_3^2+z_4^2}}.$$

*Proof.* The proof of this theorem is similar to the proof of Theorem 5.3. It is omitted here. □

**Theorem 5.11.** *Let $z = (z_1, z_2, z_3, z_4) \in \mathbb{R}^4 \setminus \{0\}$ with $(z_1^2 + z_2^2)(z_3^2 + z_4^2) = 0$. Let $x = (x_1, x_2, x_3, x_4)$ and $y = (y_1, y_2, y_3, y_4) \in \mathbb{R}^4$.*

(a) *Suppose that $z_1^2 + z_2^2 = 0$ ($z_3^2 + z_4^2 \neq 0$).*

(i) *If $y_1^2 + y_2^2 = 0$ and $x \in \widehat{D}^* g(z)(y)$, then $x$ and $y$ satisfy that*

$$x_1 = y_1 = 0 = x_2 = y_2 = 0,$$

$$x_3 = y_3 \frac{(z_3^2 + 3z_4^2)z_3}{(z_3^2 + z_4^2)\sqrt{z_3^2 + z_4^2}} + y_4 \frac{2z_4^2 z_4}{(z_3^2 + z_4^2)\sqrt{z_3^2 + z_4^2}},$$

$$x_4 = y_3 \frac{-(3z_3^2 + z_4^2)z_4}{(z_3^2 + z_4^2)\sqrt{z_3^2 + z_4^2}} + y_4 \frac{2z_3^2 z_3}{(z_3^2 + z_4^2)\sqrt{z_3^2 + z_4^2}}.$$

(ii) *If $y_1^2 + y_2^2 \neq 0$, then $\widehat{D}^* g(z)(y) = \emptyset$.*

(b) *Suppose that $z_3^2 + z_4^2 = 0$ ($z_1^2 + z_2^2 \neq 0$).*

(i) *If $y_3^2 + y_4^2 = 0$ and $x \in \widehat{D}^* g(z)(y)$, then $x$ and $y$ satisfy that*

$$x_1 = y_1 \frac{(z_1^2 + 3z_2^2)z_1}{(z_1^2 + z_2^2)\sqrt{z_1^2 + z_2^2}} + y_2 \frac{2z_2^2 z_2}{(z_1^2 + z_2^2)\sqrt{z_1^2 + z_2^2}},$$

$$x_2 = y_1 \frac{(-3z_1^2 - z_2^2)z_2}{(z_1^2 + z_2^2)\sqrt{z_1^2 + z_2^2}} + y_2 \frac{2z_1^2 z_1}{(z_1^2 + z_2^2)\sqrt{z_1^2 + z_2^2}},$$

$$x_3 = y_3 = 0 = x_4 = y_4 = 0,$$

(ii) *If $y_3^2 + y_4^2 \neq 0$, then $\widehat{D}^* g(z)(y) = \emptyset$.*

*Proof.* Proof of (a). Let $y = (y_1, y_2, y_3, y_4) \in \mathbb{R}^4$ and $x = (x_1, x_2, x_3, x_4) \in \mathbb{R}^4$, in order to check if $x \in \widehat{D}^* f(z)(y)$ or not, we take a special direction in the limit (5.12) as $u = (u_1, u_2, z_3, z_4)$ with $u \to z$, which is equivalent to $(u_1, u_2) \to (z_1, z_2)$. Let $f: \mathbb{R}^2 \to \mathbb{R}^2$ be defined by (1.1) in [15]. By definition (5.11), we have

$$g(u) - g(z) = \left( \frac{u_1^2 - u_2^2}{\sqrt{u_1^2 + u_2^2}}, \frac{2u_1 u_2}{\sqrt{u_1^2 + u_2^2}}, \frac{z_3^2 - z_4^2}{\sqrt{z_3^2 + z_4^2}}, \frac{2z_3 z_4}{\sqrt{z_3^2 + z_4^2}} \right) - \left( \frac{z_1^2 - z_2^2}{\sqrt{z_1^2 + z_2^2}}, \frac{2z_1 z_2}{\sqrt{z_1^2 + z_2^2}}, \frac{z_3^2 - z_4^2}{\sqrt{z_3^2 + z_4^2}}, \frac{2z_3 z_4}{\sqrt{z_3^2 + z_4^2}} \right).$$

This implies that

$$\langle y, g(u) - g(z) \rangle = \left\langle (y_1, y_2), \left( \frac{u_1^2 - u_2^2}{\sqrt{u_1^2 + u_2^2}}, \frac{2u_1 u_2}{\sqrt{u_1^2 + u_2^2}} \right) - \left( \frac{z_1^2 - z_2^2}{\sqrt{z_1^2 + z_2^2}}, \frac{2z_1 z_2}{\sqrt{z_1^2 + z_2^2}} \right) \right\rangle$$

$$= \langle (y_1, y_2), f(u_1, u_2) - f(z_1, z_2) \rangle.$$

With respect to this special limit, the limit (5.12) for $g$ in $\mathbb{R}^4$ can be converted to a limit in $\mathbb{R}^2$ with respect to the mapping $f: \mathbb{R}^2 \to \mathbb{R}^2$ be defined by (5.1).

$$\limsup_{u \to z} \frac{\langle x, u-z \rangle - \langle y, g(u)-g(z) \rangle}{\|u-z\| + \|g(u)-g(z)\|}$$

$$\geq \lim_{u=(u_1,u_2,z_3,z_4) \to z} \frac{\langle x, u-z \rangle - \langle y, g(u)-g(z) \rangle}{\|u-z\| + \|g(u)-g(z)\|}$$

$$= \lim_{(u_1,u_2) \to (z_1,z_2)} \frac{\langle (x_1,x_2), (u_1,u_2)-(z_1,z_2) \rangle - \langle (y_1,y_2), f(u_1,u_2)-f(z_1,z_2) \rangle}{\|(u_1,u_2)-(z_1,z_2)\| + \|f(u_1,u_2)-f(z_1,z_2)\|}.$$

Then, with respect to the case $z_1^2 + z_2^2 = 0$ (This implies that $(z_1, z_2)$ is the origin in $\mathbb{R}^2$), Part (i) can be proved by Propositions 5.2. Part (ii) can be similarly proved. □

**Theorem 5.12.** *Let* $z = (z_1, z_2, z_3, z_4) \in \mathbb{R}^4$ *with* $(z_1^2 + z_2^2)(z_3^2 + z_4^2) \neq 0$. *Let* $x = (x_1, x_2, x_3, x_4)$ *and* $y = (y_1, y_2, y_3, y_4) \in \mathbb{R}^4$. *If* $x \in \widehat{D}^* g(z)(y)$, *then* $x$ *and* $y$ *satisfy the following conditions.*

(i) $\|x\|^2 = \|y\|^2 + 12 \left( y_1 \frac{z_1 z_2}{z_1^2 + z_2^2} - y_2 \frac{z_1^2 - z_2^2}{2(z_1^2 + z_2^2)} \right)^2 + 12 \left( y_3 \frac{z_3 z_4}{z_3^2 + z_4^2} - y_4 \frac{z_3^2 - z_4^2}{2(z_3^2 + z_4^2)} \right)^2.$

(ii) $\|x\| \geq \|y\|.$

(iii) $y_1 z_1 z_2 = y_2 \frac{z_1^2 - z_2^2}{2}$ *and* $y_3 z_3 z_4 = y_4 \frac{z_3^2 - z_4^2}{2} \implies \|x\| = \|y\|.$

*Proof.* Suppose that $x \in \widehat{D}^* g(z)(y)$. By Theorem 3.3 and by the proof of Proposition 5.4, we have

$$x_1^2 + x_2^2 + x_3^2 + x_4^2$$

$$= \left( y_1 \frac{(z_1^2 + 3z_2^2) z_1}{(z_1^2 + z_2^2)\sqrt{z_1^2 + z_2^2}} + y_2 \frac{2z_2^2 z_2}{(z_1^2 + z_2^2)\sqrt{z_1^2 + z_2^2}} \right)^2 + \left( y_1 \frac{(-3z_1^2 - z_2^2) z_2}{(z_1^2 + z_2^2)\sqrt{z_1^2 + z_2^2}} + y_2 \frac{2z_1^2 z_1}{(z_1^2 + z_2^2)\sqrt{z_1^2 + z_2^2}} \right)^2$$

$$+ \left( y_3 \frac{(z_3^2 + 3z_4^2) z_3}{(z_3^2 + z_4^2)\sqrt{z_3^2 + z_4^2}} + y_4 \frac{2z_4^2 z_4}{(z_3^2 + z_4^2)\sqrt{z_3^2 + z_4^2}} \right)^2 + \left( y_3 \frac{-(3z_3^2 + z_4^2) z_4}{(z_3^2 + z_4^2)\sqrt{z_3^2 + z_4^2}} + y_4 \frac{2z_3^2 z_3}{(z_3^2 + z_4^2)\sqrt{z_3^2 + z_4^2}} \right)^2$$

$$= y_1^2 + y_2^2 + 12 \left( y_1 \frac{z_1 z_2}{z_1^2 + z_2^2} - y_2 \frac{z_1^2 - z_2^2}{2(z_1^2 + z_2^2)} \right)^2 + y_3^2 + y_4^2 + 12 \left( y_3 \frac{z_3 z_4}{z_3^2 + z_4^2} - y_4 \frac{z_3^2 - z_4^2}{2(z_3^2 + z_4^2)} \right)^2. \quad \square$$

**Corollary 5.13.** *Let* $z = (z_1, z_2, z_3, z_4) \in \mathbb{R}^4 \setminus \{\theta\}$ *with* $(z_1^2 + z_2^2)(z_3^2 + z_4^2) = 0$. *Let* $x = (x_1, x_2, x_3, x_4)$ *and* $y = (y_1, y_2, y_3, y_4) \in \mathbb{R}^4$.

(a) *Suppose that* $z_1^2 + z_2^2 = 0$ ($z_3^2 + z_4^2 \neq 0$). *If* $y_1^2 + y_2^2 = 0$ *and* $x \in \widehat{D}^* g(z)(y)$, *then*

(i) $\|x\|^2 = \|y\|^2 + 12 \left( y_3 \frac{z_3 z_4}{z_3^2 + z_4^2} - y_4 \frac{z_3^2 - z_4^2}{2(z_3^2 + z_4^2)} \right)^2.$

(ii) $\|x\| \geq \|y\|.$

(iii) $y_3 z_3 z_4 = y_4 \frac{z_3^2 - z_4^2}{2} \implies \|x\| = \|y\|.$

(b) *Suppose that $z_3^2 + z_4^2 = 0$ ($z_1^2 + z_2^2 \neq 0$). If $y_3^2 + y_4^2 = 0$ and $x \in \widehat{D}^*g(z)(y)$, then*

(i) $\|x\|^2 = \|y\|^2 + 12\left(y_1 \frac{z_1 z_2}{z_1^2 + z_2^2} - y_2 \frac{z_1^2 - z_2^2}{2(z_1^2 + z_2^2)}\right)^2$.

(ii) $\|x\| \geq \|y\|$.

(iii) $y_1 z_1 z_2 = y_2 \frac{z_1^2 - z_2^2}{2} \implies \|x\| = \|y\|$.

*Proof.* Notice that if $y_1^2 + y_2^2 = 0$, then $\|y\|^2 = y_3^2 + y_4^2$. If $y_3^2 + y_4^2 = 0$, then $\|y\|^2 = y_1^2 + y_2^2$. Then, this corollary follows from Theorems 5.11 and 5.12. $\square$

**Theorem 5.14.** *Let $g: \mathbb{R}^4 \to \mathbb{R}^4$ be defined by (5.1). Then*

$$\hat{\alpha}(g(\bar{z})) = 1, \text{ for } \bar{z} = (\bar{z}_1, \bar{z}_2, \bar{z}_3, \bar{z}_4) \in \mathbb{R}^4 \setminus \{\theta\}.$$

*Proof.* For any $\bar{z}, \bar{w} \in \mathbb{R}^4$ with $g(\bar{z}) = \bar{w}$ and $\eta > 0$, let $\mathbb{B}(\bar{z}, \eta)$ and $\mathbb{B}(\bar{w}, \eta)$ denote the closed balls in $\mathbb{R}^4$ with radius $\eta$ and centered at point $\bar{z}$ and $\bar{w}$, respectively. By definition, we have

$$\hat{\alpha}(g(\bar{z})) = \sup_{\eta > 0} \inf\{\|x\|: x \in \widehat{D}^*g(z)(y), z \in \mathbb{B}(\bar{z}, \eta), g(z) \in \mathbb{B}(\bar{w}, \eta), \|y\| = 1\}$$

$$= \sup_{\eta > 0} \inf\{\|x\|: x \in \widehat{D}^*g(z)(y), z \in \mathbb{B}(\bar{z}, \eta), g(z) \in \mathbb{B}(\bar{w}, \eta), \|y\| = 1\}. \tag{5.14}$$

In case $\theta \in \mathbb{B}(\bar{z}, \eta)$ and $g(\theta) = \theta \in \mathbb{B}(\bar{w}, \eta)$, then by Proposition 5.4, we have that

$$\widehat{D}^*g(\theta)(y) = \emptyset, \text{ for any } y \in \mathbb{R}^4 \text{ with } \|y\| = 1. \tag{5.15}$$

Let $z = (z_1, z_2, z_3, z_4) \in \mathbb{R}^4 \setminus \{\theta\}$ with $(z_1^2 + z_2^2)(z_3^2 + z_4^2) = 0$. Suppose that $z \in \mathbb{B}(\bar{z}, \eta)$, for some $\eta > 0$. Let $y = (y_1, y_2, y_3, y_4) \in \mathbb{R}^4$ with $\|y\| = 1$. By Theorem 5.11, we have that

$$z_1^2 + z_2^2 = 0 \text{ and } y_1^2 + y_2^2 \neq 0 \implies \widehat{D}^*g(z)(y) = \emptyset. \tag{5.16}$$

And $\quad\quad\quad z_3^2 + z_4^2 = 0 \text{ and } y_3^2 + y_4^2 \neq 0 \implies \widehat{D}^*g(z)(y) = \emptyset. \tag{5.17}$

From (5.14), (5.15), (5.16) and (5.17), we see that in the limit set in (4.4), there are some point $z \in \mathbb{B}(\bar{z}, \eta)$ and $y \in \mathbb{R}^4$ with $\|y\| = 1$ such that $\widehat{D}^*g(z)(y) = \emptyset$. Since $g: \mathbb{R}^4 \to \mathbb{R}^4$ defined by (5.11) is continuous on $\mathbb{R}^4$, it yields that, for every $\eta > 0$,

$$\{z \in \mathbb{B}(\bar{z}, \eta): g(z) \in \mathbb{B}(\bar{w}, \eta), \|y\| = 1\} \neq \emptyset.$$

By Theorems 5.11 and 5.12, this implies that

$$\{\|x\|: x \in \widehat{D}^*g(z)(y), z \in \mathbb{B}(\bar{z}, \eta), g(z) \in \mathbb{B}(\bar{w}, \eta), \|y\| = 1\} \neq \emptyset, \text{ for every } \eta > 0.$$

However, if $\widehat{D}^*g(z)(y) \neq \emptyset$, for $z \in \mathbb{B}(\bar{z}, \eta)$ and $y \in \mathbb{R}^4$ with $\|y\| = 1$, then, for any $x \in \widehat{D}^*g(z)(y)$, by Theorem 5.12 and Corollary 5.13, we have $\|x\| \geq \|y\|$. One can show that, for any $z = (z_1, z_2, z_3, z_4) \in \mathbb{B}(\bar{z}, r)$, the following system of equations has solutions for $y_1, y_2, y_3, y_4$.

$$\begin{cases} y_1^2 + y_2^2 + y_3^2 + y_4^2 = 1, \\ y_1 z_1 z_2 = y_2 \frac{z_1^2 - z_2^2}{2}, \\ y_3 z_3 z_4 = y_4 \frac{z_3^2 - z_4^2}{2}. \end{cases}$$

Meanwhile, following systems of equations have solutions for $y_1, y_2$ and for $y_3, y_4$, respectively:

$$\begin{cases} y_1^2 + y_2^2 = 1, \\ y_1 z_1 z_2 = y_2 \frac{z_1^2 - z_2^2}{2}, \end{cases} \quad \text{and} \quad \begin{cases} y_3^2 + y_4^2 = 1, \\ y_3 z_3 z_4 = y_4 \frac{z_3^2 - z_4^2}{2}. \end{cases}$$

By Theorem 5.12 and Corollary 5.13, we have

$$\sup_{\eta > 0} \inf\{\|x\|: x \in \widehat{D}^* g(z)(y), z \in \mathbb{B}(\bar{z}, \eta), g(z) \in \mathbb{B}(\bar{w}, \eta), \|y\| = 1\}$$

$$= \sup_{\eta > 0} \inf\left\{\|x\|: x \in \widehat{D}^* g(z)(y), z \in \mathbb{B}(\bar{z}, \eta), g(z) \in \mathbb{B}(\bar{w}, \eta), \|y\| = 1, y_1 z_1 z_2 = y_2 \frac{z_1^2 - z_2^2}{2}, y_3 z_3 z_4 = y_4 \frac{z_3^2 - z_4^2}{2}\right\}$$

$$= \sup_{\eta > 0} \inf\left\{\|y\|: x \in \widehat{D}^* f(z)(y), z \in \mathbb{B}(\bar{z}, \eta), g(z) \in \mathbb{B}(\bar{w}, \eta), \|y\| = 1, y_1 z_1 z_2 = y_2 \frac{z_1^2 - z_2^2}{2}, y_3 z_3 z_4 = y_4 \frac{z_3^2 - z_4^2}{2}\right\}$$

$$= \sup_{\eta > 0} \{1\} = 1. \qquad \square$$

Next, we consider a mapping $h: \mathbb{R}^4 \to \mathbb{R}^4$ defined by (5.18) below, which is slightly modified the mapping $g: \mathbb{R}^4 \to \mathbb{R}^4$ defined by (5.1). Define $h: \mathbb{R}^4 \to \mathbb{R}^4$, for any $x = (x_1, x_2, x_3, x_4) \in \mathbb{R}^4$, by

$$h((x_1, x_2, x_3, x_4)) = \begin{cases} \left(\frac{x_1^2 - x_2^2}{\|x\|}, \frac{2x_1 x_2}{\|x\|}, \frac{x_3^2 - x_4^2}{\|x\|}, \frac{2x_3 x_4}{\|x\|}\right), & \text{if } x \neq \theta, \\ \theta, & \text{if } x = \theta. \end{cases} \tag{5.18}$$

We see that $h$ is a continuous mapping on $\mathbb{R}^4$. Similar, to the mapping $g: \mathbb{R}^4 \to \mathbb{R}^4$ defined by (5.1), one can study the properties of $h$ at the origin $\theta$. However, in this paper, we only study the Fréchet differentiability and the Mordukhovich differentiability of $h$ on $\mathbb{R}^4 \setminus \{\theta\}$. We can only estimate the covering constant for $h$. It is because that will be very complicated to precisely find it.

**Theorem 5.15.** *Let* $z = (z_1, z_2, z_3, z_4) \in \mathbb{R}^4 \setminus \{\theta\}$. *Then we have that*

(a) $h$ *is Fréchet differentiable at* $z$ *and*

$$\nabla h(z) = \begin{pmatrix} \frac{(z_1^2 + 3z_2^2 + 2z_3^2 + 2z_4^2)z_1}{\|z\|^3} & \frac{2z_2^2 z_2}{\|z\|^3} & \frac{-(z_3^2 - z_4^2)z_1}{\|z\|^3} & \frac{-2z_3 z_4 z_1}{\|z\|^3} \\ \frac{-(3z_1^2 + z_2^2 + 2z_3^2 + 2z_4^2)z_2}{\|z\|^3} & \frac{2z_1^2 z_1}{\|z\|^3} & \frac{-(z_3^2 - z_4^2)z_2}{\|z\|^3} & \frac{-2z_3 z_4 z_2}{\|z\|^3} \\ \frac{-(z_1^2 - z_2^2)z_3}{\|z\|^3} & \frac{-2z_1 z_2 z_3}{\|z\|^3} & \frac{(2z_1^2 + 2z_2^2 + z_3^2 + 3z_4^2)z_3}{\|z\|^3} & \frac{2z_4^2 z_4}{\|z\|^3} \\ \frac{-(z_1^2 - z_2^2)z_4}{\|z\|^3} & \frac{-2z_1 z_2 z_4}{\|z\|^3} & \frac{-(2z_1^2 + 2z_2^2 + 3z_3^2 + z_4^2)z_4}{\|z\|^3} & \frac{2z_3^2 z_3}{\|z\|^3} \end{pmatrix}.$$

(b) $h$ *is Mordukhovich differentiable at* $z$ *and*

$$\widehat{D}^*h(z) = \begin{pmatrix} \frac{(z_1^2+3z_2^2+2z_3^2+2z_4^2)z_1}{\|z\|^3} & \frac{-(3z_1^2+z_2^2+2z_3^2+2z_4^2)z_2}{\|z\|^3} & \frac{-(z_1^2-z_2^2)z_3}{\|z\|^3} & \frac{-(z_1^2-z_2^2)z_4}{\|z\|^3} \\ \frac{2z_2^2 z_2}{\|z\|^3} & \frac{2z_1^2 z_1}{\|z\|^3} & \frac{-2z_1 z_2 z_3}{\|z\|^3} & \frac{-2z_1 z_2 z_4}{\|z\|^3} \\ \frac{-(z_3^2-z_4^2)z_1}{\|z\|^3} & \frac{-(z_3^2-z_4^2)z_2}{\|z\|^3} & \frac{(2z_1^2+2z_2^2+z_3^2+3z_4^2)z_3}{\|z\|^3} & \frac{-(2z_1^2+2z_2^2+3z_3^2+z_4^2)z_4}{\|z\|^3} \\ \frac{-2z_3 z_4 z_1}{\|z\|^3} & \frac{-2z_3 z_4 z_2}{\|z\|^3} & \frac{2z_4^2 z_4}{\|z\|^3} & \frac{2z_3^2 z_3}{\|z\|^3} \end{pmatrix}.$$

*Proof.* By Lemma 5.1 and Proposition 5.2, the proof of this theorem is similar to the proof of Theorem 5.3. It is omitted here. □

From the representation of the Mordukhovich derivative of $h$, it is very difficult and complicated to calculate the covering constant for $h$ by using the similar techniques used in Theorem 5.15. We only consider some special cases, which may be interesting for some readers.

**Proposition 5.16.** *Let* $\bar{z} = (\bar{z}_1, \bar{z}_2, \bar{z}_3, \bar{z}_4) \neq \theta$ *with* $\bar{w} = h(\bar{z})$. *We have*

(a)   *If* $\bar{z}_i = 0$, *for some* $i = 1, 2, 3, 4$, *then*

(i) $\qquad\qquad\qquad \hat{\alpha}(h, \bar{z}, \bar{w}) \leq \frac{2\bar{z}_{3-i}^3}{\|\bar{z}\|^3}$, *for* $i = 1, 2$,

(ii) $\qquad\qquad\qquad \hat{\alpha}(h, \bar{z}, \bar{w}) \leq \frac{2\bar{z}_{7-i}^3}{\|\bar{z}\|^3}$, *for* $i = 3, 4$.

*In particular, if* $\bar{z}_1 = \bar{z}_2 = 0$ *or* $\bar{z}_3 = \bar{z}_4 = 0$, *then* $\hat{\alpha}(h, \bar{z}, \bar{w}) = 0$.

(b)   *If* $|\bar{z}_1| = |\bar{z}_2| = |\bar{z}_3| = |\bar{z}_4|$, *then*, $\hat{\alpha}(h, \bar{z}, \bar{w}) \leq \frac{1}{\sqrt{2}}$.

*Proof.* Proof of (a). Suppose $\bar{z}_1 = 0$. By the continuity of $h$ around $\bar{z} = (\bar{z}_1, \bar{z}_2, \bar{z}_3, \bar{z}_4) \neq \theta$, for any $\eta > 0$, we can have $z \in \mathbb{R}^4 \setminus \{\theta\}$ such that

$$z = (z_1, z_2, z_3, z_4) \in \mathbb{B}(\bar{z}, \eta) \setminus \{\theta\} \text{ satisfying } h(z) \in \mathbb{B}(\bar{w}, \eta) \text{ and } z_1 = 0. \tag{5.19}$$

For any $\eta > 0$, let

$$D(\bar{z}, \eta) = \{z \in \mathbb{R}^4 \setminus \{\theta\}: z \text{ satisfies } (4.9)\}$$

For any given $\eta > 0$, let $z \in \mathbb{R}^4 \setminus \{\theta\}$ satisfying (5.19). Let $x = (x_1, x_2, x_3, x_4), y = (y_1, y_2, y_3, y_4) \in \mathbb{R}^4$. If $x = \widehat{D}^*h(z)(y)$, then, by Theorem 5.15, we have

$$x_1 = y_2 \frac{2z_2^2 z_2}{\|z\|^3}$$

$$x_2 = y_1 \frac{-(z_2^2+2z_3^2+2z_4^2)z_2}{\|z\|^3} + y_3 \frac{-(z_3^2-z_4^2)z_2}{\|z\|^3} + y_4 \frac{-2z_3 z_4 z_2}{\|z\|^3}$$

$$x_3 = y_1 \frac{z_2^2 z_3}{\|z\|^3} + y_3 \frac{(2z_2^2+z_3^2+3z_4^2)z_3}{\|z\|^3} + y_4 \frac{2z_4^2 z_4}{\|z\|^3},$$

$$x_4 = y_1 \frac{z_2^2 z_4}{\|z\|^3} + y_3 \frac{-(2z_2^2+3z_3^2+z_4^2)z_4}{\|z\|^3} + y_4 \frac{2z_3^2 z_3}{\|z\|^3}. \tag{5.20}$$

By (5.19) and (5.20), we calculate

$$\hat{\alpha}(h,\bar{z},\bar{w}) = \sup_{\eta>0} \inf\{\|x\|: x \in \widehat{D}^*h(z)(y), z \in \mathbb{B}(\bar{z},\eta)\setminus\{\theta\}, h(z) \in \mathbb{B}(\bar{w},\eta), \|y\|=1\}$$

$$\leq \sup_{\eta>0} \inf\{\|x\|: x \in \widehat{D}^*h(z)(y), z \in (\mathbb{B}(\bar{z},\eta)\setminus\{\theta\}) \cap D(\bar{z},\eta), h(z) \in \mathbb{B}(\bar{w},\eta), \|y\|=1\}$$

$$\leq \sup_{\eta>0} \inf\{\|x\|: x \in \widehat{D}^*h(z)(y), z \in (\mathbb{B}(\bar{z},\eta)\setminus\{\theta\}) \cap D(\bar{z},\eta), h(z) \in \mathbb{B}(\bar{w},\eta), y_2=\|y\|=1\}$$

$$= \sup_{\eta>0} \inf\left\{\frac{2z_2^3}{\|z\|^3}: x \in \widehat{D}^*h(z)(y), z \in (\mathbb{B}(\bar{z},\eta)\setminus\{\theta\}) \cap D(\bar{z},\eta), h(z) \in \mathbb{B}(\bar{w},\eta), y_2=\|y\|=1\right\}$$

$$\leq \sup_{\eta>0} \inf\left\{\frac{2\bar{z}_2^3}{\|\bar{z}\|^3}: x \in \widehat{D}^*h(\bar{z})(y), \bar{z} \in (\mathbb{B}(\bar{z},\eta)\setminus\{\theta\}) \cap D(\bar{z},\eta), h(\bar{z}) \in \mathbb{B}(\bar{w},\eta), y_2=\|y\|=1\right\}$$

$$= \frac{2\bar{z}_2^3}{\|\bar{z}\|^3}.$$

Rest cases can be similarly proved.

Next, we prove (b). Suppose $|\bar{z}_1|=|\bar{z}_2|=|\bar{z}_3|=|\bar{z}_4|$. Let $x=(x_1,x_2,x_3,x_4) \in \mathbb{R}^4$. We take $y=(0,1,0,0) \in \mathbb{R}^4$ (We may take $y=(0,0,0,1)$). If $x=\widehat{D}^*h(z)(y)$, by Theorem 4.9, we have

$$x_1 = \frac{2z_2^2 z_2}{\|z\|^3}, \quad x_2 = \frac{2z_1^2 z_1}{\|z\|^3}, \quad x_3 = \frac{-2z_1 z_2 z_3}{\|z\|^3}, \quad x_4 = \frac{-2z_1 z_2 z_4}{\|z\|^3}. \tag{5.21}$$

By (5.21), we calculate

$$\hat{\alpha}(h,\bar{z},\bar{w}) = \sup_{\eta>0} \inf\{\|x\|: x \in \widehat{D}^*h(z)(y), z \in \mathbb{B}(\bar{z},\eta)\setminus\{\theta\}, h(z) \in \mathbb{B}(\bar{w},\eta), \|y\|=1\}$$

$$\leq \sup_{\eta>0} \inf\{\|x\|: x \in \widehat{D}^*h(z)(y), z \in (\mathbb{B}(\bar{z},\eta)\setminus\{\theta\}) \cap D(\bar{z},\eta), h(z) \in \mathbb{B}(\bar{w},\eta), y_2=\|y\|=1\}$$

$$\leq \sup_{\eta>0} \inf\left\{\sqrt{\frac{16\bar{z}_1^6}{(2|\bar{z}_1|)^6}}: x \in \widehat{D}^*h(\bar{z})(y), \bar{z} \in (\mathbb{B}(\bar{z},\eta)\setminus\{\theta\}) \cap D(\bar{z},\eta), h(\bar{z}) \in \mathbb{B}(\bar{w},\eta), y_2=\|y\|=1\right\}$$

$$= \frac{1}{\sqrt{2}}. \qquad \square$$

## 6. Covering Constants for Some Single-Valued Mappings from $\mathbb{R}^n$ to $\mathbb{R}^m$ with $n \geq m$

In this section, we provide some concrete examples for calculating the covering constants for some single-valued mappings from $\mathbb{R}^n$ to $\mathbb{R}^m$ with $n \geq m$ with some specific techniques. From these examples, we will see the difficulty for calculating the covering constants.

**6.1. A norm preserving mapping from $\mathbb{R}^3$ to $\mathbb{R}^2$.** We define $f: \mathbb{R}^3 \to \mathbb{R}^2$ by

$$f(x) = \left(\sqrt{x_1^2 + x_2^2}, x_3\right), \text{ for every } x=(x_1,x_2,x_3) \in \mathbb{R}^3.$$

This implies that, for $z=(z_1,z_2,z_3) \in \mathbb{R}^3$, we have that $\nabla f(z): \mathbb{R}^3 \to \mathbb{R}^2$ and $\widehat{D}^*f(z): \mathbb{R}^2 \to \mathbb{R}^3$. Then, $f$ has the following properties.

(a) $f$ is norm preserving, this is, $\|f(x)\|_2 = \|x\|_3$, for any $x \in \mathbb{R}^3$.

(b) $f$ is Fréchet differentiable and Mordukhovich differentiable. For any $z = (z_1, z_2, z_3) \in \mathbb{R}^3$ with $z_1^2 + z_2^2 > 0$, we have

$$\nabla f(z) = \begin{pmatrix} \frac{z_1}{\sqrt{z_1^2+z_2^2}} & 0 \\ \frac{x_2}{\sqrt{z_1^2+z_2^2}} & 0 \\ 0 & 1 \end{pmatrix} \quad \text{and} \quad \widehat{D}^* f(z) = \begin{pmatrix} \frac{z_1}{\sqrt{z_1^2+z_2^2}} & \frac{z_2}{\sqrt{z_1^2+z_2^2}} & 0 \\ 0 & 0 & 1 \end{pmatrix}.$$

(c) The covering constant for $f$ is constant on $\mathbb{R}^3$ satisfying

$$\hat{\alpha}(f, \bar{z}, \bar{w}) = 1, \text{ for any } \bar{z} = (\bar{z}_1, \bar{z}_2, \bar{z}_3) \in \mathbb{R}^3 \text{ with } \bar{z}_1^2 + \bar{z}_2^2 > 0 \text{ and } \bar{w} = f(\bar{z}) \in \mathbb{R}^2.$$

*Proof.* The proofs of (a) and (b) are straight forward and it is omitted here. We only prove (c). Let $z = (z_1, z_2, z_3) \in \mathbb{R}^3$ with $z_1^2 + z_2^2 > 0$. Let $x = (x_1, x_2, x_3) \in \mathbb{R}^3$ and $y = (y_1, y_2) \in \mathbb{R}^2$. If $x = \widehat{D}^* f(z)(y)$, by part (b), we have that

$$(x_1, x_2, x_3) = (y_1, y_2) \begin{pmatrix} \frac{z_1}{\sqrt{z_1^2+z_2^2}} & \frac{z_2}{\sqrt{z_1^2+z_2^2}} & 0 \\ 0 & 0 & 1 \end{pmatrix} = \left( \frac{y_1 z_1}{\sqrt{z_1^2+z_2^2}}, \frac{y_1 z_2}{\sqrt{z_1^2+z_2^2}}, y_2 \right). \tag{6.1}$$

This implies that

$$x = \widehat{D}^* f(z)(y) \implies \|x\|_3 = \|y\|_2. \tag{6.2}$$

Let $\bar{z} = (\bar{z}_1, \bar{z}_2, \bar{z}_3) \in \mathbb{R}^3$ with $\bar{z}_1^2 + \bar{z}_2^2 > 0$ and $\bar{w} = (\bar{w}_1, \bar{w}_2) \in \mathbb{R}^2$ with $\bar{w} = f(\bar{z})$. By the continuity of $f$ on $\mathbb{R}^2$, for any given $\eta > 0$, we have

$$\{z = (z_1, z_2, z_3) \in \mathbb{B}(\bar{z}, \eta) : z_1^2 + z_2^2 > 0, f(z) \in \mathbb{B}(\bar{w}, \eta)\} \neq \emptyset.$$

By (6.1) and (6.2), we calculate

$$\hat{\alpha}(f, \bar{z}, \bar{w}) = \sup_{\eta > 0} \inf\{\|x\|_3 : x = \widehat{D}^* f(z)(y), z = (z_1, z_2, z_3) \in \mathbb{B}(\bar{z}, \eta), z_1^2 + z_2^2 > 0, f(z) \in \mathbb{B}(\bar{w}, \eta), \|y\|_2 = 1\}$$

$$= \sup_{\eta > 0} \inf\{\|y\|_2 : x = \widehat{D}^* f(z)(y), z = (z_1, z_2, z_3) \in \mathbb{B}(\bar{z}, \eta), z_1^2 + z_2^2 > 0, f(z) \in \mathbb{B}(\bar{w}, \eta), \|y\|_2 = 1\}$$

$$= 1. \qquad \square$$

### 6.2. A polynomial mapping from $\mathbb{R}^3$ to $\mathbb{R}^2$.

Define $f \equiv (f_1, f_2) \colon \mathbb{R}^3 \to \mathbb{R}^2$, for every $x = (x_1, x_2, x_3) \in \mathbb{R}^3$, by

$$f(x) = (x_1 x_2, x_1 x_3).$$

Then, $f$ has the following properties.

(a) $f$ is Fréchet differentiable and Mordukhovich differentiable on $\mathbb{R}^3$. For any $z = (z_1, z_2, z_3) \in \mathbb{R}^3$, we have

$$\nabla f(z) = \begin{pmatrix} z_2 & z_3 \\ z_1 & 0 \\ 0 & z_1 \end{pmatrix} \quad \text{and} \quad \widehat{D}^* f(z) = \begin{pmatrix} z_2 & z_1 & 0 \\ z_3 & 0 & z_1 \end{pmatrix}.$$

(b) The covering constant for $f$ on $\mathbb{R}^3$ satisfies

$$\hat{\alpha}(f, \bar{z}, \bar{w}) = |\bar{z}_1|, \text{ for any } \bar{z} = (\bar{z}_1, \bar{z}_2, \bar{z}_3) \in \mathbb{R}^3 \text{ with } \bar{w} = f(\bar{z}) \in \mathbb{R}^2.$$

*Proof.* The proof of (a) is straight forward and it is omitted here. We only prove (b). Let $z = (z_1, z_2, z_3) \in \mathbb{R}^3$. Let $x = (x_1, x_2, x_3) \in \mathbb{R}^3$ and $y = (y_1, y_2) \in \mathbb{R}^2$. If $x = \widehat{D}^* f(z)(y)$, by part (a), we have that

$$(x_1, x_2, x_3) = (y_1, y_2) \begin{pmatrix} z_2 & z_1 & 0 \\ z_3 & 0 & z_1 \end{pmatrix} = (y_1 z_2 + y_2 z_3,\ y_1 z_1,\ y_2 z_1).$$

This implies that if $x = \widehat{D}^* f(z)(y)$, then

$$\|x\|_3^2 = x_1^2 + x_2^2 + x_3^2 = (y_1 z_2 + y_2 z_3)^2 + y_1^2 z_1^2 + y_2^2 z_1^2$$

$$= (y_1 z_2 + y_2 z_3)^2 + (y_1^2 + y_2^2) z_1^2. \tag{6.3}$$

Let $\bar{z} = (\bar{z}_1, \bar{z}_2, \bar{z}_3) \in \mathbb{R}^3$ with $\bar{w} = (\bar{w}_1, \bar{w}_2) \in \mathbb{R}^2$ with $\bar{w} = f(\bar{z})$. By (4.3), we calculate

$$\hat{\alpha}(f, \bar{z}, \bar{w}) = \sup_{\eta > 0} \inf\{\|x\|_3 : x \in \widehat{D}^* f(z)(y), z = (z_1, z_2, z_3) \in \mathbb{B}(\bar{z}, \eta), f(z) \in \mathbb{B}(\bar{w}, \eta), \|y\|_2 = 1\}$$

$$= \sup_{\eta > 0} \inf\left\{\sqrt{(y_1 z_2 + y_2 z_3)^2 + z_1^2} : x = \widehat{D}^* f(z)(y), z = (z_1, z_2, z_3) \in \mathbb{B}(\bar{z}, \eta), f(z) \in \mathbb{B}(\bar{w}, \eta), \|y\|_2 = 1\right\}$$

$$= \sup_{\eta > 0} \inf\left\{\sqrt{(y_1 z_2 + y_2 z_3)^2 + z_1^2} : x = \widehat{D}^* f(z)(y), z = (z_1, z_2, z_3) \in \mathbb{B}(\bar{z}, \eta), f(z) \in \mathbb{B}(\bar{w}, \eta), \|y\|_2 = 1, y_1 z_2 + y_2 z_3 = 0\right\}$$

$$= \sup_{\eta > 0} \inf\{|z_1| : x = \widehat{D}^* f(z)(y), z = (z_1, z_2, z_3) \in \mathbb{B}(\bar{z}, \eta), f(z) \in \mathbb{B}(\bar{w}, \eta), \|y\|_2 = 1, y_1 z_2 + y_2 z_3 = 0\}$$

$$= \lim_{\eta \downarrow 0}\{|z_1| : x = \widehat{D}^* f(z)(y), z = (z_1, z_2, z_3) \in \mathbb{B}(\bar{z}, \eta), f(z) \in \mathbb{B}(\bar{w}, \eta), \|y\|_2 = 1, y_1 z_2 + y_2 z_3 = 0\}$$

$$= |\bar{z}_1|. \qquad \square$$

**6.3. Another polynomial mapping from $\mathbb{R}^3$ to $\mathbb{R}^2$.** Define $f: \mathbb{R}^3 \to \mathbb{R}^2$, for $x = (x_1, x_2, x_3) \in \mathbb{R}^3$, by

$$f(x) = (x_1^2 x_3, x_2^2 x_3).$$

Then, $f$ has the following properties.

(a) $f$ is Fréchet differentiable and Mordukhovich differentiable on $\mathbb{R}^3$. For any $z = (z_1, z_2, z_3) \in \mathbb{R}^3$, we have

$$\nabla f(z) = \begin{pmatrix} 2z_1 z_3 & 0 \\ 0 & 2z_2 z_3 \\ z_1^2 & z_2^2 \end{pmatrix} \quad \text{and} \quad \widehat{D}^* f(z) = \begin{pmatrix} 2z_1 z_3 & 0 & z_1^2 \\ 0 & 2z_2 z_3 & z_2^2 \end{pmatrix}.$$

(b) Let $\bar{z} = (\bar{z}_1, \bar{z}_2, \bar{z}_3) \in \mathbb{R}^3$, $\bar{w} = (\bar{w}_1, \bar{w}_2) \in \mathbb{R}^2$ with $\bar{w} = f(\bar{z})$. Suppose $\bar{z}_1 = \bar{z}_2$. Then

$$\hat{\alpha}(f,\bar{z},\bar{w}) \leq 2|\bar{z}_1\bar{z}_3|.$$

*Proof.* The proof of (a) is straight forward and it is omitted here. We only prove (b). Let $z = (z_1, z_2, z_3) \in \mathbb{R}^3$. Let $x = (x_1, x_2, x_3) \in \mathbb{R}^3$ and $y = (y_1, y_2) \in \mathbb{R}^2$. If $x = \widehat{D}^*f(z)(y)$, by part (a), we have that

$$(x_1, x_2, x_3) = (y_1, y_2)\begin{pmatrix} 2z_1z_3 & 0 & z_1^2 \\ 0 & 2z_2z_3 & z_2^2 \end{pmatrix} = (2y_1z_1z_3,\ 2y_2z_2z_3,\ y_1z_1^2 + y_2z_2^2).$$

This implies that if $x = \widehat{D}^*f(z)(y)$, under the condition $y_1^2 + y_2^2 = 1$, we have

$$\|x\|_3^2 = x_1^2 + x_2^2 + x_3^2$$

$$= 4y_1^2z_1^2z_3^2 + 4y_2^2z_2^2z_3^2 + (y_1z_1^2 + y_2z_2^2)^2$$

$$= 4(y_1^2z_1^2 + y_2^2z_2^2)z_3^2 + (y_1z_1^2 + y_2z_2^2)^2.$$

Let $\bar{z} = (\bar{z}_1, \bar{z}_2, \bar{z}_3) \in \mathbb{R}^3$, $\bar{w} = (\bar{w}_1, \bar{w}_2) \in \mathbb{R}^2$ with $\bar{w} = f(\bar{z})$. Suppose $\bar{z}_1 = \bar{z}_2$. By the above equation, we calculate

$$\hat{\alpha}(f,\bar{z},\bar{w}) = \sup_{\eta>0}\inf\{\|x\|_3 : x \in \widehat{D}^*f(z)(y), z = (z_1,z_2,z_3) \in \mathbb{B}(\bar{z},\eta), f(z) \in \mathbb{B}(\bar{w},\eta), \|y\|_2 = 1\}$$

$$= \sup_{\eta>0}\inf\left\{\sqrt{4(y_1^2z_1^2 + y_2^2z_2^2)z_3^2 + (y_1z_1^2 + y_2z_2^2)^2} : x = \widehat{D}^*f(z)(y), z = (z_1,z_2,z_3) \in \mathbb{B}(\bar{z},\eta), f(z) \in \mathbb{B}(\bar{w},\eta), \|y\|_2 = 1\right\}$$

$$\leq \sup_{\eta>0}\inf\left\{\sqrt{4(y_1^2z_1^2 + y_2^2z_2^2)z_3^2 + (y_1z_1^2 + y_2z_2^2)^2} : z \in \mathbb{B}(\bar{z},\eta), z_1 = z_2, f(z) \in \mathbb{B}(\bar{w},\eta), \|y\|_2 = 1, y_1 + y_2 = 0\right\}$$

$$= \sup_{\eta>0}\inf\{2|z_1z_3| : x = \widehat{D}^*f(z)(y), z = (z_1,z_2,z_3) \in \mathbb{B}(\bar{z},\eta), f(z) \in \mathbb{B}(\bar{w},\eta), \|y\|_2 = 1, y_1 + y_2 = 0\}$$

$$= \lim_{\eta\downarrow 0}\{2|z_1z_3| : x = \widehat{D}^*f(z)(y), z = (z_1,z_2,z_3) \in \mathbb{B}(\bar{z},\eta), f(z) \in \mathbb{B}(\bar{w},\eta), \|y\|_2 = 1, y_1 + y_2 = 0\}$$

$$= 2|\bar{z}_1\bar{z}_3|. \qquad \square$$

**6.4. A continuous rational mapping from $\mathbb{R}^3$ to $\mathbb{R}^2$.** Define $f \equiv (f_1, f_2): \mathbb{R}^3 \to \mathbb{R}^2$, for every $x = (x_1, x_2, x_3) \in \mathbb{R}^3$, by

$$f(x) = \left(\frac{x_1}{1+x_3^2}, \frac{x_2}{1+x_3^2}\right).$$

Then, $f$ has the following properties.

(a) $f$ is Fréchet differentiable and Mordukhovich differentiable on $\mathbb{R}^3$. For any $z = (z_1, z_2, z_3) \in \mathbb{R}^3$, we have

$$\nabla f(z) = \begin{pmatrix} \frac{1}{1+z_3^2} & 0 \\ 0 & \frac{1}{1+z_3^2} \\ \frac{-2z_1z_3}{(1+z_3^2)^2} & \frac{-2z_2z_3}{(1+z_3^2)^2} \end{pmatrix} \quad \text{and} \quad \widehat{D}^*f(z) = \begin{pmatrix} \frac{1}{1+z_3^2} & 0 & \frac{-2z_1z_3}{(1+z_3^2)^2} \\ 0 & \frac{1}{1+z_3^2} & \frac{-2z_2z_3}{(1+z_3^2)^2} \end{pmatrix}.$$

(b) The covering constant for $f$ on $\mathbb{R}^3$ satisfies
$$\hat{\alpha}(f,\bar{z},\bar{w}) = \frac{1}{1+\bar{z}_3^2}, \text{ for any } \bar{z} = (\bar{z}_1, \bar{z}_2, \bar{z}_3) \in \mathbb{R}^3 \text{ with } \bar{w} = f(\bar{z}) \in \mathbb{R}^2.$$

*Proof.* The proof of (a) is straight forward and it is omitted here. We only prove (b). Let $z = (z_1, z_2, z_3) \in \mathbb{R}^3$. Let $x = (x_1, x_2, x_3) \in \mathbb{R}^3$ and $y = (y_1, y_2) \in \mathbb{R}^2$. If $x = \widehat{D}^* f(z)(y)$, by part (a), we have that

$$(x_1, x_2, x_3) = (y_1, y_2) \begin{pmatrix} \frac{1}{1+z_3^2} & 0 & \frac{-2z_1 z_3}{(1+z_3^2)^2} \\ 0 & \frac{1}{1+z_3^2} & \frac{-2z_2 z_3}{(1+z_3^2)^2} \end{pmatrix} = \left( \frac{y_1}{1+z_3^2}, \frac{y_2}{1+z_3^2}, \frac{-2y_1 z_1 z_3 - 2y_2 z_2 z_3}{(1+z_3^2)^2} \right).$$

This implies that if $x = \widehat{D}^* f(z)(y)$, under the condition $y_1^2 + y_2^2 = 1$, we have

$$\|x\|_3^2 = x_1^2 + x_2^2 + x_3^2$$

$$= \frac{y_1^2}{(1+z_3^2)^2} + \frac{y_2^2}{(1+z_3^2)^2} + \frac{4z_3^2 (y_1 z_1 + y_2 z_2)^2}{(1+z_3^2)^4}$$

$$= \frac{1}{(1+z_3^2)^2} + \frac{4z_3^2 (y_1 z_1 + y_2 z_2)^2}{(1+z_3^2)^4}.$$

Let $\bar{z} = (\bar{z}_1, \bar{z}_2, \bar{z}_3) \in \mathbb{R}^3$, $\bar{w} = (\bar{w}_1, \bar{w}_2) \in \mathbb{R}^2$ with $\bar{w} = f(\bar{z})$. By the above equation, we calculate

$$\hat{\alpha}(f,\bar{z},\bar{w}) = \sup_{\eta>0} \inf \{\|x\|_3 : x \in \widehat{D}^* f(z)(y), z = (z_1, z_2, z_3) \in \mathbb{B}(\bar{z},\eta), f(z) \in \mathbb{B}(\bar{w},\eta), \|y\|_2 = 1\}$$

$$= \sup_{\eta>0} \inf \left\{ \sqrt{\frac{1}{(1+z_3^2)^2} + \frac{4z_3^2 (y_1 z_1 + y_2 z_2)^2}{(1+z_3^2)^4}} : x = \widehat{D}^* f(z)(y), z = (z_1, z_2, z_3) \in \mathbb{B}(\bar{z},\eta), f(z) \in \mathbb{B}(\bar{w},\eta), \|y\|_2 = 1 \right\}$$

$$= \sup_{\eta>0} \inf \left\{ \sqrt{\frac{1}{(1+z_3^2)^2} + \frac{4z_3^2 (y_1 z_1 + y_2 z_2)^2}{(1+z_3^2)^4}} : x = \widehat{D}^* f(z)(y), z \in \mathbb{B}(\bar{z},\eta), f(z) \in \mathbb{B}(\bar{w},\eta), \|y\|_2 = 1, y_1 z_1 + y_2 z_2 = 0 \right\}$$

$$= \sup_{\eta>0} \inf \left\{ \frac{1}{1+z_3^2} : x = \widehat{D}^* f(z)(y), z = (z_1, z_2, z_3) \in \mathbb{B}(\bar{z},\eta), f(z) \in \mathbb{B}(\bar{w},\eta), \|y\|_2 = 1, y_1 z_2 + y_2 z_3 = 0 \right\}$$

$$= \lim_{\eta \downarrow 0} \left\{ \frac{1}{1+z_3^2} : x = \widehat{D}^* f(z)(y), z = (z_1, z_2, z_3) \in \mathbb{B}(\bar{z},\eta), f(z) \in \mathbb{B}(\bar{w},\eta), \|y\|_2 = 1, y_1 z_2 + y_2 z_3 = 0 \right\}$$

$$= \frac{1}{1+\bar{z}_3^2}. \qquad \square$$

**6.5. A linear and continuous mapping from $\mathbb{R}^3$ to $\mathbb{R}^2$.** Define $f: \mathbb{R}^3 \to \mathbb{R}^2$, for $x = (x_1, x_2, x_3) \in \mathbb{R}^3$ by

$$f(x) = (x_1, x_2).$$

This is a linear and continuous mapping. Then, $f$ has the following properties.

(a) $f$ is Fréchet differentiable and Mordukhovich differentiable on $\mathbb{R}^3$. For any $z = (z_1, z_2, z_3) \in \mathbb{R}^3$, we have

$$\nabla f(z) = \begin{pmatrix} 1 & 0 \\ 0 & 1 \\ 0 & 0 \end{pmatrix} \quad \text{and} \quad \widehat{D}^* f(z) = \begin{pmatrix} 1 & 0 & 0 \\ 0 & 1 & 0 \end{pmatrix}.$$

(b) The covering constant for $f$ on $\mathbb{R}^3$ satisfies

$$\hat{\alpha}(f, \bar{z}, \bar{w}) = 1, \text{ for any } \bar{z} = (\bar{z}_1, \bar{z}_2, \bar{z}_3) \in \mathbb{R}^3 \text{ with } \bar{w} = f(\bar{z}) \in \mathbb{R}^2.$$

*Proof.* The proof of (a) is straight forward and it is omitted here. We only prove (b). Let $z = (z_1, z_2, z_3) \in \mathbb{R}^3$. Let $x = (x_1, x_2, x_3) \in \mathbb{R}^3$ and $y = (y_1, y_2) \in \mathbb{R}^2$. If $x = \widehat{D}^* f(z)(y)$, by part (a), we have that

$$(x_1, x_2, x_3) = (y_1, y_2) \begin{pmatrix} 1 & 0 & 0 \\ 0 & 1 & 0 \end{pmatrix} = (y_1, y_2, 0).$$

This implies that if $x = \widehat{D}^* f(z)(y)$, then

$$\|x\|_3^2 = x_1^2 + x_2^2 + x_3^2 = y_1^2 + y_2^2 = \|y\|_2^2. \tag{6.4}$$

Let $\bar{z} = (\bar{z}_1, \bar{z}_2, \bar{z}_3) \in \mathbb{R}^3$ with $\bar{w} = (\bar{w}_1, \bar{w}_2) \in \mathbb{R}^2$ with $\bar{w} = f(\bar{z})$. By (6.4), we calculate

$$\hat{\alpha}(f, \bar{z}, \bar{w}) = \sup_{\eta > 0} \inf \{\|x\|_3 : x = \widehat{D}^* f(z)(y), z \in \mathbb{B}(\bar{z}, \eta), f(z) \in \mathbb{B}(\bar{w}, \eta), \|y\|_2 = 1\}$$

$$= \sup_{\eta > 0} \inf \{\|y\|_2 : x = \widehat{D}^* f(z)(y), z \in \mathbb{B}(\bar{z}, \eta), f(z) \in \mathbb{B}(\bar{w}, \eta), \|y\|_2 = 1\}$$

$$= \sup_{\eta > 0} \inf \{1 : x = \widehat{D}^* f(z)(y), z \in \mathbb{B}(\bar{z}, \eta), f(z) \in \mathbb{B}(\bar{w}, \eta), \|y\|_2 = 1\}$$

$$= 1. \qquad \square$$

**6.6. A trigonometric mapping.** We define $f: \mathbb{R}^2 \to \mathbb{R}^2$ by

$$f(x_1, x_2) = (\sin(x_1 + x_2), \cos(x_1 + x_2)), \text{ for any } (x_1, x_2) \in \mathbb{R}^2.$$

Then, $f$ has the following properties.

(a) $f$ is a norm constant mapping with

$$\|f(x)\| = 1, \text{ for any } x \in \mathbb{R}^2.$$

(b) $f$ is Fréchet and Mordukhovich differentiable on $\mathbb{R}^2$. For any $z = (z_1, z_2) \in \mathbb{R}^2$, we have

$$\nabla f(z) = \begin{pmatrix} \cos(z_1 + z_2) & -\sin(z_1 + z_2) \\ \cos(z_1 + z_2) & -\sin(z_1 + z_2) \end{pmatrix} \quad \text{and} \quad \widehat{D}^* f(z) = \begin{pmatrix} \cos(z_1 + z_2) & \cos(z_1 + z_2) \\ -\sin(z_1 + z_2) & -\sin(z_1 + z_2) \end{pmatrix}.$$

(c) The covering constant for $f$ is constant with

$$\hat{\alpha}(f, \bar{z}, \bar{w}) = 0, \text{ for any } \bar{z} = (\bar{z}_1, \bar{z}_2) \in \mathbb{R}^2 \text{ with } \bar{w} = f(\bar{z}).$$

*Proof.* The proofs of (a) and (b) are straight forward and they are omitted here. We only prove (c).

Let $x = (x_1, x_2)$ and $y = (y_1, y_2) \in \mathbb{R}^2$. If $x = \widehat{D}^* f(z)(y)$, then

$$x_1 = y_1 \cos(z_1 + z_2) - y_2 \sin(z_1 + z_2),$$

$$x_2 = y_1 \cos(z_1 + z_2) - y_2 \sin(z_1 + z_2).$$

This implies that

$$x_1^2 + x_2^2 = (y_1 \cos(z_1 + z_2) - y_2 \sin(z_1 + z_2))^2 + (y_1 \cos(z_1 + z_2) - y_2 \sin(z_1 + z_2))^2$$

$$= 2(y_1 \cos(z_1 + z_2) - y_2 \sin(z_1 + z_2))^2.$$

Then we obtain that

$$\|x\| = \sqrt{2}|y_1 \cos(z_1 + z_2) - y_2 \sin(z_1 + z_2)|, \text{ for } x = \widehat{D}^* f(z)(y). \quad (6.5)$$

For any given point $\bar{z} = (\bar{z}_1, \bar{z}_2) \neq \theta$ and $\bar{w} = (\bar{w}_1, \bar{w}_2) \in \mathbb{R}^2$ with $\bar{w} = f(\bar{z})$, by (4.4), we calculate

$$\hat{\alpha}(f, \bar{z}, \bar{w}) = \sup_{\eta > 0} \inf\{\|x\|: x \in \widehat{D}^* f(z)(y), z \in \mathbb{B}(\bar{z}, \eta), f(z) \in \mathbb{B}(\bar{w}, \eta), \|y\| = 1\}.$$

By the continuity of $f$ on $\mathbb{R}^2$, for any given $\eta > 0$, we have

$$\{z \in \mathbb{B}(\bar{z}, \eta): f(z) \in \mathbb{B}(\bar{w}, \eta)\} \neq \emptyset.$$

For any given $\eta > 0$, and $z = (z_1, z_2) \in \mathbb{B}(\bar{z}, \eta)$ and $f(z) \in \mathbb{B}(\bar{w}, \eta)$, the following system of equations has solutions with respect to $y = (y_1, y_2)$:

$$\begin{cases} y_1^2 + y_2^2 = 1, \\ y_1 \cos(z_1 + z_2) - y_2 \sin(z_1 + z_2) = 0. \end{cases} \quad (6.6)$$

By (6.5) and (6.6), we have

$$\hat{\alpha}(f, \bar{z}, \bar{w}) = \sup_{\eta > 0} \inf\{\|x\|: x \in \widehat{D}^* f(z)(y), z \in \mathbb{B}(\bar{z}, \eta), f(z) \in \mathbb{B}(\bar{w}, \eta), \|y\| = 1\}$$

$$\leq \sup_{\eta > 0} \inf\{\|x\|: x \in \widehat{D}^* f(z)(y), z \in \mathbb{B}(\bar{z}, \eta) \setminus \{\theta\}, f(z) \in \mathbb{B}(\bar{w}, \eta), \|y\| = 1, y_1 \cos(z_1 + z_2) - y_2 \sin(z_1 + z_2) = 0\}$$

$$= \sup_{\eta > 0} \inf\{0: x \in \widehat{D}^* f(z)(y), z \in \mathbb{B}(\bar{z}, \eta) \setminus \{\theta\}, f(z) \in \mathbb{B}(\bar{w}, \eta), \|y\| = 1, y_1 \cos(z_1 + z_2) - y_2 \sin(z_1 + z_2) = 0\}$$

$$= 0.$$

This proves that

$$\hat{\alpha}(f, \bar{z}, \bar{w}) = 0, \text{ for any } \bar{z} = (\bar{z}_1, \bar{z}_2) \in \mathbb{R}^2 \text{ with } \bar{w} = f(\bar{z}). \qquad \square$$

**6.7. A polynomial mapping.** We define $f: \mathbb{R}^2 \to \mathbb{R}^2$ by

$$f(x_1, x_2) = (x_1^2 - x_2^2, 2x_1 x_2), \text{ for any } (x_1, x_2) \in \mathbb{R}^2.$$

Then, $f$ has the following properties.

(a) $f$ is a norm-expansion mapping with

$$\|f(x)\| = \|x\|^2, \text{ for any } x \in \mathbb{R}^2.$$

(b) $f$ is Fréchet differentiable and Mordukhovich differentiable $\mathbb{R}^2$. For any $z = (z_1, z_2) \in \mathbb{R}^2$,

$$\nabla f(z) = \begin{pmatrix} 2z_1 & 2z_2 \\ -2z_2 & 2z_1 \end{pmatrix} \quad \text{and} \quad \widehat{D}^* f(z) = \begin{pmatrix} 2z_1 & -2z_2 \\ 2z_2 & 2z_1 \end{pmatrix}.$$

(c) The covering constant for $f$ satisfies

$$\hat{\alpha}(f, \bar{z}, f(\bar{z})) = 2\|\bar{z}\|, \text{ for any } \bar{z} \in \mathbb{R}^2.$$

In particular, we have $\hat{\alpha}(f, \theta, \theta) = 0$.

*Proof.* Proof of (a). This is a polynomial mapping. We calculate the norms of this mapping.

$$\|f(x_1, x_2)\|^2 = x_1^4 - 2x_1^2 x_2^2 + x_2^4 + 4x_1^2 x_2^2 = (x_1^2 + x_2^2)^2 = \|(x_1, x_2)\|^4.$$

This implies that $f$ is a norm-expansion mapping with

$$\|f(x)\| = \|x\|^2, \text{ for any } x \in \mathbb{R}^2.$$

Proof of (b). Since $f$ is a polynomial mapping. By Proposition 3.3, $f$ is Fréchet differentiable and Mordukhovich differentiable at every point $z = (z_1, z_2) \in \mathbb{R}^2$. By Theorem 3.1, we obtain $\nabla f(z)$ and $\widehat{D}^* f(z)$ as given in Part (b).

Proof of (c). Let $z = (z_1, z_2)$, $x = (x_1, x_2)$ and $y = (y_1, y_2) \in \mathbb{R}^2$. If $x = \widehat{D}^* f(z)(y)$, we have that

$$x_1 = 2y_1 z_1 + 2y_2 z_2,$$

$$x_2 = -2y_1 z_2 + 2y_2 z_1.$$

We calculate

$$x_1^2 + x_2^2 = (2y_1 z_1 + 2y_2 z_2)^2 + (-2y_1 z_2 + 2y_2 z_1)^2$$

$$= 4y_1^2 z_1^2 + 8y_1 y_2 z_1 z_2 + 4y_2^2 z_2^2 + 4y_1^2 z_2^2 - 8y_1 y_2 z_1 z_2 + 4y_2^2 z_1^2$$

$$= 4(y_1^2 + y_2^2)(z_1^2 + z_2^2).$$

This implies

$$\|x\| = 2\|y\|\|z\|, \text{ for } x = \widehat{D}^* f(z)(y). \tag{6.7}$$

For any given point $\bar{z} = (\bar{z}_1, \bar{z}_2)$ and $\bar{w} \in \mathbb{R}^2$ with $\bar{w} = f(\bar{z})$, by the continuity of $f$ on $\mathbb{R}^2$, we have that,

$$\{z \in \mathbb{B}(\bar{z}, \eta): f(z) \in \mathbb{B}(\bar{w}, \eta)\} \neq \emptyset, \text{ for any } \eta > 0. \tag{6.8}$$

By the Mordukhovich differentiability of $f$ on $\mathbb{R}^2$ in Part (a) of this example, we have

$$\widehat{D}^* f(z)(y) \neq \emptyset, \text{ for any } y \in \mathbb{R}^2. \tag{6.9}$$

By (6.8) and (6.9), we obtain that, for any $\eta > 0$,

$$\{x = \widehat{D}^*f(z)(y): z \in \mathbb{B}(\bar{z},\eta), f(z) \in \mathbb{B}(\bar{w},\eta), \|y\| = 1\} \neq \emptyset. \tag{6.10}$$

For any positive numbers $p$ and $q$ with $p < q$, we have

$$\{z \in \mathbb{B}(\bar{z},p): f(z) \in \mathbb{B}(\bar{w},p)\} \subseteq \{z \in \mathbb{B}(\bar{z},q): f(z) \in \mathbb{B}(\bar{w},q)\}. \tag{6.11}$$

By (6.10) to (6.11), for any $0 < p < q$, we get that

$$\{x = \widehat{D}^*f(z)(y): z \in \mathbb{B}(\bar{z},p), f(z) \in \mathbb{B}(\bar{w},p), \|y\| = 1\}$$

$$\subseteq \{x = \widehat{D}^*f(z)(y): z \in \mathbb{B}(\bar{z},q), f(z) \in \mathbb{B}(\bar{w},q), \|y\| = 1\}.$$

This implies that $\{x = \widehat{D}^*f(z)(y): z \in \mathbb{B}(\bar{z},\eta), f(z) \in \mathbb{B}(\bar{w},\eta), \|y\| = 1\}$ is a decreasing (inclusion) net with respect to $\eta \downarrow 0$. Hence, the following net

$$\inf\{\|x\|: x = \widehat{D}^*f(z)(y), z \in \mathbb{B}(\bar{z},\eta), f(z) \in \mathbb{B}(\bar{w},\eta), \|y\| = 1\}, \tag{6.12}$$

is an increasing net of nonnegative numbers. By (6.7), we estimate (6.12).

$$\inf\{\|x\|: x = \widehat{D}^*f(z)(y), z \in \mathbb{B}(\bar{z},\eta), f(z) \in \mathbb{B}(\bar{w},\eta), \|y\| = 1\}$$

$$\leq \inf\{\|x\|: x = \widehat{D}^*f(\bar{z})(y), \bar{z} \in \mathbb{B}(\bar{z},\eta), f(\bar{z}) = \bar{w} \in \mathbb{B}(\bar{w},\eta), \|y\| = 1\}$$

$$= \inf\{2\|\bar{z}\|: x = \widehat{D}^*f(\bar{z})(y), \bar{z} \in \mathbb{B}(\bar{z},\eta), f(\bar{z}) = \bar{w} \in \mathbb{B}(\bar{w},\eta), \|y\| = 1\}$$

$$= 2\|\bar{z}\|.$$

This implies that (6.12) is a bounded increasing net of nonnegative numbers. We calculate

$$\hat{\alpha}(f,\bar{z},f(\bar{z})) = \sup_{\eta>0} \inf\{\|x\|: x = \widehat{D}^*f(z)(y), z \in \mathbb{B}(\bar{z},\eta), f(z) \in \mathbb{B}(\bar{w},\eta), \|y\| = 1\}$$

$$= \lim_{\eta \downarrow 0} \inf\{2\|y\|\|z\|, z \in \mathbb{B}(\bar{z},\eta), f(z) \in \mathbb{B}(\bar{w},\eta), \|y\| = 1\}$$

$$= \lim_{\eta \downarrow 0} \inf\{2\|z\|, z \in \mathbb{B}(\bar{z},\eta), f(z) \in \mathbb{B}(\bar{w},\eta), \|y\| = 1\}$$

$$= \lim_{\eta \downarrow 0} \inf\{2\|z\|, z \in \mathbb{B}(\bar{z},\eta), f(z) \in \mathbb{B}(\bar{w},\eta), \|y\| = 1\}$$

$$= 2\|\bar{z}\|. \qquad \square$$

**6.8. An exponential mapping.** We define $f: \mathbb{R}^2 \to \mathbb{R}^2$ by

$$f(x) = (e^{x_1+x_2},\ e^{-x_1-x_2}), \text{ for any } x = (x_1, x_2) \in \mathbb{R}^2.$$

Then, $f$ is a continuous mapping on $\mathbb{R}^2$, which has the following properties.

(a) $f$ is Fréchet differentiable and Mordukhovich differentiable on $\mathbb{R}^2$. For each $z = (z_1, z_2) \in \mathbb{R}^2$, we have

$$\nabla f(z) = \begin{pmatrix} e^{z_1+z_2} & -e^{-z_1-z_2} \\ e^{z_1+z_2} & -e^{-z_1-z_2} \end{pmatrix} \quad \text{and} \quad \widehat{D}^*f(z) = \begin{pmatrix} e^{z_1+z_2} & e^{z_1+z_2} \\ -e^{-z_1-z_2} & -e^{-z_1-z_2} \end{pmatrix}.$$

(b) The covering constant for $f$ is constant satisfying

$$\hat{\alpha}(f, \bar{z}, \bar{w}) = 0, \text{ for any } \bar{z} = (\bar{z}_1, \bar{z}_2) \in \mathbb{R}^2 \text{ with } \bar{w} = f(\bar{z}).$$

*Proof.* The proof of (a) is straight forward and it is omitted here. We prove (b). Let $\bar{z} = (\bar{z}_1, \bar{z}_2) \in \mathbb{R}^2$ with $\bar{w} = f(\bar{z})$. Let $x = (x_1, x_2)$ and $y = (y_1, y_2) \in \mathbb{R}^2$, if $x = \widehat{D}^* f(z)(y)$, by part (a), we have that

$$x_1 = y_1 e^{z_1 + z_2} - y_2 e^{-z_1 - z_2},$$

$$x_2 = y_1 e^{z_1 + z_2} - y_2 e^{-z_1 - z_2}.$$

Hence, if $x = \widehat{D}^* f(z)(y)$, then

$$x_1^2 + x_2^2 = (y_1 e^{z_1 + z_2} - y_2 e^{-z_1 - z_2})^2 + (y_1 e^{z_1 + z_2} - y_2 e^{-z_1 - z_2})^2$$

$$= 2(y_1^2 e^{2(z_1 + z_2)} - 2 y_1 y_2 + y_2^2 e^{-2(z_1 + z_2)}).$$

This implies that

$$\|x\| = \sqrt{2} \sqrt{y_1^2 e^{2(z_1 + z_2)} - 2 y_1 y_2 + y_2^2 e^{-2(z_1 + z_2)}}$$

$$= \sqrt{2} \sqrt{(y_1 e^{z_1 + z_2} - y_2 e^{-z_1 - z_2})^2}$$

$$= \sqrt{2} |y_1 e^{z_1 + z_2} - y_2 e^{-z_1 - z_2}|, \text{ for } x = \widehat{D}^* f(z)(y). \tag{6.13}$$

For any given point $\bar{z} = (\bar{z}_1, \bar{z}_2)$ and $\bar{w} = (\bar{w}_1, \bar{w}_2) \in \mathbb{R}^2$ with $\bar{w} = f(\bar{z})$, we calculate

$$\hat{\alpha}(f, \bar{z}, \bar{w}) = \sup_{\eta > 0} \inf \{\|x\| : x \in \widehat{D}^* f(z)(y), z \in \mathbb{B}(\bar{z}, \eta), f(z) \in \mathbb{B}(\bar{w}, \eta), \|y\| = 1\}.$$

By the continuity of $f$ on $\mathbb{R}^2$, for any given $\eta > 0$, we have

$$\{z \in \mathbb{B}(\bar{z}, \eta) : f(z) \in \mathbb{B}(\bar{w}, \eta)\} \neq \emptyset.$$

For any given $\eta > 0$, and $z = (z_1, z_2) \in \mathbb{B}(\bar{z}, \eta)$ and $f(z) \in \mathbb{B}(\bar{w}, \eta)$, one can prove that the following system of equations has solutions with respect to $y = (y_1, y_2)$:

$$\begin{cases} y_1^2 + y_2^2 = 1, \\ y_1 e^{z_1 + z_2} - y_2 e^{-z_1 - z_2} = 0. \end{cases}$$

By (6.13), this implies that

$$\hat{\alpha}(f, \bar{z}, \bar{w}) = \sup_{\eta > 0} \inf \{\|x\| : x \in \widehat{D}^* f(z)(y), z \in \mathbb{B}(\bar{z}, \eta), f(z) \in \mathbb{B}(\bar{w}, \eta), \|y\| = 1\}$$

$$\leq \sup_{\eta > 0} \inf \{\|x\| : x \in \widehat{D}^* f(z)(y), z = (z_1, z_2) \in \mathbb{B}(\bar{z}, \eta), f(z) \in \mathbb{B}(\bar{w}, \eta), \|y\| = 1, y_1 e^{z_1 + z_2} - y_2 e^{-z_1 - z_2} = 0\}$$

$$\leq \sup_{\eta > 0} \inf \{\sqrt{2} |y_1 e^{z_1 + z_2} - y_2 e^{-z_1 - z_2}| : x \in \widehat{D}^* f(z)(y), z \in \mathbb{B}(\bar{z}, \eta), f(z) \in \mathbb{B}(\bar{w}, \eta),$$
$$\|y\| = 1, y_1 e^{z_1 + z_2} - y_2 e^{-z_1 - z_2} = 0\}$$

$$= \sup_{\eta > 0} \inf \{0 : x \in \widehat{D}^* f(z)(y), z \in \mathbb{B}(\bar{z}, \eta) \setminus \{\theta\}, f(z) \in \mathbb{B}(\bar{w}, \eta), \|y\| = 1, y_1 e^{z_1 + z_2} - y_2 e^{-z_1 - z_2} = 0\}$$

= 0.

This proves that

$$\hat{\alpha}(f, \bar{z}, \bar{w}) = 0, \text{ for any } \bar{z} = (\bar{z}_1, \bar{z}_2) \in \mathbb{R}^2 \text{ with } \bar{w} = f(\bar{z}).$$  □

**6.9. A logarithm mapping.** We define $f: \mathbb{R}^2 \to \mathbb{R}^2$ by

$$f(x) = \left(\ln(1 + x_1^2 + x_2^2), \frac{1}{1+x_1^2+x_2^2}\right), \text{ for } x = (x_1, x_2) \in \mathbb{R}^2.$$

Then, $f$ is a continuous mapping on $\mathbb{R}^2$, which has the following properties.

(a) $f$ is Fréchet differentiable and Mordukhovich differentiable on $\mathbb{R}^2$. For $z = (z_1, z_2) \in \mathbb{R}^2$ we have

$$\nabla f(z) = \begin{pmatrix} \frac{2z_1}{1+z_1^2+z_2^2} & -\frac{2z_1}{(1+z_1^2+z_2^2)^2} \\ \frac{2z_2}{1+z_1^2+z_2^2} & -\frac{2z_2}{(1+z_1^2+z_2^2)^2} \end{pmatrix} \quad \text{and} \quad \widehat{D}^* f(z) = \begin{pmatrix} \frac{2z_1}{1+z_1^2+z_2^2} & \frac{2z_2}{1+z_1^2+z_2^2} \\ -\frac{2z_1}{(1+z_1^2+z_2^2)^2} & -\frac{2z_2}{(1+z_1^2+z_2^2)^2} \end{pmatrix}.$$

(b) The covering constant for $f$ is constant on $\mathbb{R}^2$ satisfying

$$\hat{\alpha}(f, \bar{z}, \bar{w}) = 0, \text{ for any } \bar{z} = (\bar{z}_1, \bar{z}_2) \text{ with } \bar{w} = f(\bar{z}).$$

*Proof.* By Lemma 3.1 and Theorem 3.2, the proof of (a) is straight forward and it is omitted here. We prove (b). Let $\bar{z} = (\bar{z}_1, \bar{z}_2)$ with $\bar{w} = f(\bar{z})$. Let $x = (x_1, x_2)$ and $y = (y_1, y_2) \in \mathbb{R}^2$, if $x = \widehat{D}^* f(z)(y)$, by part (a), we have that

$$x_1 = y_1 \frac{2z_1}{1+z_1^2+z_2^2} - y_2 \frac{2z_1}{(1+z_1^2+z_2^2)^2},$$

$$x_2 = y_1 \frac{2z_2}{1+z_1^2+z_2^2} - y_2 \frac{2z_2}{(1+z_1^2+z_2^2)^2}.$$

This implies that

$$x_1^2 + x_2^2 = \left(y_1 \frac{2z_1}{1+z_1^2+z_2^2} - y_2 \frac{2z_1}{(1+z_1^2+z_2^2)^2}\right)^2 + \left(y_1 \frac{2z_2}{1+z_1^2+z_2^2} - y_2 \frac{2z_2}{(1+z_1^2+z_2^2)^2}\right)^2$$

$$= y_1^2 \frac{4z_1^2}{(1+z_1^2+z_2^2)^2} - y_1 y_2 \frac{8z_1^2}{(1+z_1^2+z_2^2)^3} + y_2^2 \frac{4z_1^2}{(1+z_1^2+z_2^2)^4} + y_1^2 \frac{4z_2^2}{(1+z_1^2+z_2^2)^2} - y_1 y_2 \frac{8z_2^2}{(1+z_1^2+z_2^2)^3} + y_2^2 \frac{4z_2^2}{(1+z_1^2+z_2^2)^4}$$

$$= y_1^2 \frac{4(z_1^2+z_2^2)}{(1+z_1^2+z_2^2)^2} - y_1 y_2 \frac{8(z_1^2+z_2^2)}{(1+z_1^2+z_2^2)^3} + y_2^2 \frac{4(z_1^2+z_2^2)}{(1+z_1^2+z_2^2)^4}$$

$$= \frac{4(z_1^2+z_2^2)}{(1+z_1^2+z_2^2)^2}\left(y_1 - \frac{y_2}{1+z_1^2+z_2^2}\right)^2.$$

This implies that

$$\|x\| = \frac{2\sqrt{z_1^2+z_2^2}}{1+z_1^2+z_2^2}\left|y_1 - \frac{y_2}{1+z_1^2+z_2^2}\right|, \text{ for } x = \widehat{D}^* f(z)(y).$$

Then, we calculate

$$\hat{\alpha}(f, \bar{z}, \bar{w}) = \sup_{\eta > 0} \inf\{\|x\|: x \in \widehat{D}^*f(z)(y), z \in \mathbb{B}(\bar{z}, \eta), f(z) \in \mathbb{B}(\bar{w}, \eta), \|y\| = 1\}$$

$$= \sup_{\eta > 0} \inf\left\{\frac{2\sqrt{z_1^2+z_2^2}}{1+z_1^2+z_2^2}\left|y_1 - \frac{y_2}{1+z_1^2+z_2^2}\right|: x \in \widehat{D}^*f(z)(y), z \in \mathbb{B}(\bar{z}, \eta), f(z) \in \mathbb{B}(\bar{w}, \eta), \|y\| = 1\right\}$$

$$\leq \sup_{\eta > 0} \inf\left\{0: x \in \widehat{D}^*f(z)(y), z \in \mathbb{B}(\bar{z}, \eta), f(z) \in \mathbb{B}(\bar{w}, \eta), \|y\| = 1, y_1 - \frac{y_2}{1+z_1^2+z_2^2} = 0\right\}$$

$$= 0. \qquad \square$$

**6.10. A rational mapping with radical.** We define $f: \mathbb{R}^2 \to \mathbb{R}^2$, for $x = (x_1, x_2) \in \mathbb{R}^2$, by

$$f(x) = \begin{cases} \left(\dfrac{x_1}{\sqrt{x_1^2+x_2^2}}, \dfrac{x_2}{\sqrt{x_1^2+x_2^2}}\right), & \text{for } (x_1, x_2) \neq \theta, \\ \theta, & \text{for } (x_1, x_2) = \theta. \end{cases}$$

$f$ is a continuous mapping on $\mathbb{R}^2\setminus\{\theta\}$, which has the following properties.

(a) $f$ is Fréchet differentiable and Mordukhovich differentiable on $\mathbb{R}^2\setminus\{\theta\}$. For any $z = (z_1, z_2) \in \mathbb{R}^2\setminus\{\theta\}$, we have

$$\nabla f(z) = \begin{pmatrix} \dfrac{z_2^2}{(z_1^2+z_2^2)\sqrt{x_1^2+x_2^2}} & \dfrac{-z_1 z_2}{(z_1^2+z_2^2)\sqrt{x_1^2+x_2^2}} \\ \dfrac{-z_1 z_2}{(z_1^2+z_2^2)\sqrt{x_1^2+x_2^2}} & \dfrac{z_1^2}{(z_1^2+z_2^2)\sqrt{x_1^2+x_2^2}} \end{pmatrix} \text{ and } \widehat{D}^*f(z) = \begin{pmatrix} \dfrac{z_2^2}{(z_1^2+z_2^2)\sqrt{x_1^2+x_2^2}} & \dfrac{-z_1 z_2}{(z_1^2+z_2^2)\sqrt{x_1^2+x_2^2}} \\ \dfrac{-z_1 z_2}{(z_1^2+z_2^2)\sqrt{x_1^2+x_2^2}} & \dfrac{z_1^2}{(z_1^2+z_2^2)\sqrt{x_1^2+x_2^2}} \end{pmatrix}.$$

(b) The covering constant for $f$ on $\mathbb{R}^2\setminus\{\theta\}$ satisfies that

$$\hat{\alpha}(f, \bar{z}, \bar{w}) = 0, \text{ for any } \bar{z} \in \mathbb{R}^2\setminus\{\theta\} \text{ with } \bar{w} = f(\bar{z}).$$

*Proof.* By Lemma 3.1 and Theorem 3.2, the proof of (a) is straight forward and it is omitted here.

Proof of (b). Let $z = (z_1, z_2) \in \mathbb{R}^2\setminus\{\theta\}$, and $x = (x_1, x_2)$, $y = (y_1, y_2) \in \mathbb{R}^2$. Suppose that $x = \widehat{D}^*f(z)(y)$. By Part (a), we have that

$$x_1 = y_1 \frac{z_2^2}{(z_1^2+z_2^2)\sqrt{x_1^2+x_2^2}} + y_2 \frac{-z_1 z_2}{(z_1^2+z_2^2)\sqrt{x_1^2+x_2^2}} \text{ and } x_2 = y_1 \frac{-z_1 z_2}{(z_1^2+z_2^2)\sqrt{x_1^2+x_2^2}} + y_2 \frac{z_1^2}{(z_1^2+z_2^2)\sqrt{x_1^2+x_2^2}}.$$

This implies that

$$x_1^2 + x_2^2 = \left(y_1 \frac{z_2^2}{(z_1^2+z_2^2)\sqrt{x_1^2+x_2^2}} + y_2 \frac{-z_1 z_2}{(z_1^2+z_2^2)\sqrt{x_1^2+x_2^2}}\right)^2 + \left(y_1 \frac{-z_1 z_2}{(z_1^2+z_2^2)\sqrt{x_1^2+x_2^2}} + y_2 \frac{z_1^2}{(z_1^2+z_2^2)\sqrt{x_1^2+x_2^2}}\right)^2$$

$$= y_1^2 \frac{z_2^4}{(z_1^2+z_2^2)^3} - 2y_1 y_2 \frac{z_1 z_2^3}{(z_1^2+z_2^2)^3} + y_2^2 \frac{z_1^2 z_2^2}{(z_1^2+z_2^2)^3} + y_1^2 \frac{z_1^2 z_2^2}{(z_1^2+z_2^2)^3} - 2y_1 y_2 \frac{z_1^3 z_2}{(z_1^2+z_2^2)^3} + y_2^2 \frac{z_1^4}{(z_1^2+z_2^2)^3}$$

$$= y_1^2 \frac{z_2^2(z_1^2+z_1^2)}{(z_1^2+z_2^2)^3} - 2y_1y_2 \frac{z_1z_2(z_1^2+z_1^2)}{(z_1^2+z_2^2)^3} + y_2^2 \frac{z_1^2(z_1^2+z_1^2)}{(z_1^2+z_2^2)^3}$$

$$= \frac{(y_1z_2-y_2z_1)^2}{(z_1^2+z_2^2)^2}.$$

By the above equations, for $z = (z_1, z_2) \in \mathbb{R}^2\setminus\{\theta\}$ and $x = (x_1, x_2)$, $y = (y_1, y_2) \in \mathbb{R}^2$, we obtain

$$\|x\| = \frac{|y_1z_2-y_2z_1|}{z_1^2+z_2^2}, \text{ if } x = \widehat{D}^*f(z)(y).$$

For $\bar{z} = (\bar{z}_1, \bar{z}_2) \neq \theta$ with $\bar{w} = f(\bar{z})$, by the above equation, we calculate

$$\hat{\alpha}(f, \bar{z}, \bar{w}) = \sup_{\eta>0} \inf\{\|x\|: x \in \widehat{D}^*f(z)(y), z \in \mathbb{B}(\bar{z},\eta)\setminus\{\theta\}, f(z) \in \mathbb{B}(\bar{w},\eta), \|y\| = 1\}$$

$$= \sup_{\eta>0} \inf\left\{\frac{|y_1z_2-y_2z_1|}{z_1^2+z_2^2}: z \in \mathbb{B}(\bar{z},\eta)\setminus\{\theta\}, f(z) \in \mathbb{B}(\bar{w},\eta), \|y\| = 1\right\}$$

$$\leq \sup_{\eta>0} \inf\left\{\frac{|y_1z_2-y_2z_1|}{z_1^2+z_2^2}: z \in \mathbb{B}(\bar{z},\eta)\setminus\{\theta\}, f(z) \in \mathbb{B}(\bar{w},\eta), \|y\| = 1, y_1z_2 - y_2z_1 = 0\right\}$$

$$= \sup_{\eta>0} \inf\{0: z \in \mathbb{B}(\bar{z},\eta)\setminus\{\theta\}, f(z) \in \mathbb{B}(\bar{w},\eta), \|y\| = 1, y_1z_2 - y_2z_1 = 0\}$$

$$= 0. \qquad \square$$

For the above single-valued mappings from $\mathbb{R}^n$ to $\mathbb{R}^m$ with $n \geq m$, we able to find the exact covering constants for the considered mappings. However, in general, it may be very difficulty to precisely calculate the covering constants for single-valued mappings in finite dimensional Hilbert spaces. For example, if we slightly modify the mapping in 6.10 as below, we will see that it is extremely complicated to find the exact value of the covering constant for it.

**6.11. Another rational mapping with radical.** We define $f: \mathbb{R}^2 \to \mathbb{R}^2$, for $x = (x_1, x_2) \in \mathbb{R}^2$, by

$$f(x) = \begin{cases} \left(\dfrac{x_1^2}{\sqrt{x_1^2+x_2^2}}, \dfrac{x_2^2}{\sqrt{x_1^2+x_2^2}}\right), & \text{for } (x_1, x_2) \neq \theta, \\ \theta, & \text{for } (x_1, x_2) = \theta. \end{cases}$$

$f$ is a continuous mapping on $\mathbb{R}^2$, which has the following properties.

(a) $f$ is Fréchet differentiable and Mordukhovich differentiable on $\mathbb{R}^2\setminus\{\theta\}$. For $z = (z_1, z_2) \in \mathbb{R}^2\setminus\{\theta\}$, we have

$$\nabla f(z) = \begin{pmatrix} \dfrac{z_1z_1^2+2z_1z_2^2}{(z_1^2+z_2^2)\sqrt{x_1^2+x_2^2}} & \dfrac{-z_2^2z_1}{(z_1^2+z_2^2)\sqrt{x_1^2+x_2^2}} \\ \dfrac{-z_1^2z_2}{(z_1^2+z_2^2)\sqrt{x_1^2+x_2^2}} & \dfrac{z_2z_2^2+2z_2z_1^2}{(z_1^2+z_2^2)\sqrt{x_1^2+x_2^2}} \end{pmatrix} \text{ and } \widehat{D}^*f(z) = \begin{pmatrix} \dfrac{z_1z_1^2+2z_1z_2^2}{(z_1^2+z_2^2)\sqrt{x_1^2+x_2^2}} & \dfrac{-z_1^2z_2}{(z_1^2+z_2^2)\sqrt{x_1^2+x_2^2}} \\ \dfrac{-z_2^2z_1}{(z_1^2+z_2^2)\sqrt{x_1^2+x_2^2}} & \dfrac{z_2z_2^2+2z_2z_1^2}{(z_1^2+z_2^2)\sqrt{x_1^2+x_2^2}} \end{pmatrix}.$$

(b) For any $\bar{z} = (\bar{z}_1, \bar{z}_2) \neq \theta$ with $\bar{w} = f(\bar{z})$, the covering constant for $f$ at $\bar{z}$ has the following properties.

(i) $\hat{\alpha}(f, \bar{z}, \bar{w}) \leq \frac{1}{\sqrt{2}}$.

(ii) $\hat{\alpha}(f, \bar{z}, \bar{w}) = 0$, if $\bar{z}_1 \bar{z}_2 = 0$.

(iii) $\hat{\alpha}(f, \bar{z}, \bar{w}) \leq \frac{2|\bar{z}_1 \bar{z}_2|}{\sqrt{\bar{z}_1^4 + \bar{z}_2^4}}$.

Notice that (ii) can be immediately reduced by (iii). Since the ideas of the direct proof of (ii) may be interesting for some readers, so we provide the direct proof of (i).

*Proof.* By Lemma 3.1 and Theorem 3.2, the proof of (a) is straight forward and it is omitted here.

Proof of (b). Let $z = (z_1, z_2) \in \mathbb{R}^2 \setminus \{\theta\}$, and $x = (x_1, x_2)$, $y = (y_1, y_2) \in \mathbb{R}^2$. Suppose that $x = \hat{D}^* f(z)(y)$. By Part (a), we have that

$$x_1 = y_1 \frac{z_1 z_1^2 + 2 z_1 z_2^2}{(z_1^2 + z_2^2)\sqrt{x_1^2 + x_2^2}} + y_2 \frac{-z_2^2 z_1}{(z_1^2 + z_2^2)\sqrt{x_1^2 + x_2^2}},$$

$$x_2 = y_1 \frac{-z_1^2 z_2}{(z_1^2 + z_2^2)\sqrt{x_1^2 + x_2^2}} + y_2 \frac{z_2 z_2^2 + 2 z_2 z_1^2}{(z_1^2 + z_2^2)\sqrt{x_1^2 + x_2^2}}.$$

This implies that

$$x_1^2 + x_2^2 = \left( y_1 \frac{z_1 z_1^2 + 2 z_1 z_2^2}{(z_1^2 + z_2^2)\sqrt{x_1^2 + x_2^2}} + y_2 \frac{-z_2^2 z_1}{(z_1^2 + z_2^2)\sqrt{x_1^2 + x_2^2}} \right)^2 + \left( y_1 \frac{-z_1^2 z_2}{(z_1^2 + z_2^2)\sqrt{x_1^2 + x_2^2}} + y_2 \frac{z_2 z_2^2 + 2 z_2 z_1^2}{(z_1^2 + z_2^2)\sqrt{x_1^2 + x_2^2}} \right)^2$$

$$= y_1^2 \frac{z_1^2 z_1^4 + 4 z_1^4 z_2^2 + 4 z_1^2 z_2^4}{(z_1^2 + z_2^2)^3} - 2 y_1 y_2 \frac{z_1^4 z_2^2 + 2 z_1^2 z_2^4}{(z_1^2 + z_2^2)^3} + y_2^2 \frac{z_1^2 z_2^4}{(z_1^2 + z_2^2)^3} + y_1^2 \frac{z_1^4 z_2^2}{(z_1^2 + z_2^2)^3} - 2 y_1 y_2 \frac{z_1^2 z_2^4 + 2 z_1^4 z_2^2}{(z_1^2 + z_2^2)^3} + y_2^2 \frac{z_2^2 z_2^4 + 4 z_1^2 z_2^4 + 4 z_1^4 z_2^2}{(z_1^2 + z_2^2)^3}$$

$$= y_1^2 \frac{z_1^2 z_1^4 + 5 z_1^4 z_2^2 + 4 z_1^2 z_2^4}{(z_1^2 + z_2^2)^3} - 2 y_1 y_2 \frac{3 z_1^4 z_2^2 + 3 z_1^2 z_2^4}{(z_1^2 + z_2^2)^3} + y_2^2 \frac{z_2^2 z_2^4 + 5 z_1^2 z_2^4 + 4 z_1^4 z_2^2}{(z_1^2 + z_2^2)^3}$$

$$= \frac{1}{(z_1^2 + z_2^2)^3} \left( y_1^2 z_1^2 (z_1^4 + 5 z_1^2 z_2^2 + 4 z_2^4) - 6 y_1 y_2 z_1^2 z_2^2 (z_1^2 + z_2^2) + y_2^2 z_2^2 (z_2^4 + 5 z_1^2 z_2^2 + 4 z_1^4) \right)$$

$$= \frac{1}{(z_1^2 + z_2^2)^3} \left( y_1^2 z_1^2 (z_1^2 + 4 z_2^2)(z_1^2 + z_2^2) - 6 y_1 y_2 z_1^2 z_2^2 (z_1^2 + z_2^2) + y_2^2 z_2^2 (z_2^2 + 4 z_1^2)(z_1^2 + z_2^2) \right)$$

$$= \frac{1}{(z_1^2 + z_2^2)^2} \left( y_1^2 z_1^2 (z_1^2 + 4 z_2^2) - 6 y_1 y_2 z_1^2 z_2^2 + y_2^2 z_2^2 (z_2^2 + 4 z_1^2) \right)$$

$$= \frac{1}{(z_1^2 + z_2^2)^2} \left( y_1^2 z_1^4 + 4 y_1^2 z_1^2 z_2^2 - 6 y_1 y_2 z_1^2 z_2^2 + y_2^2 z_2^4 + 4 y_2^2 z_1^2 z_2^2 \right)$$

$$= \frac{1}{(z_1^2 + z_2^2)^2} \left( y_1^2 z_1^4 + 2 y_1 y_2 z_1^2 z_2^2 + y_2^2 z_2^4 + 4 y_1^2 z_1^2 z_2^2 - 8 y_1 y_2 z_1^2 z_2^2 + 4 y_2^2 z_1^2 z_2^2 \right)$$

$$= \frac{1}{(z_1^2 + z_2^2)^2} \left( (y_1 z_1^2 + y_2 z_2^2)^2 + 4 z_1^2 z_2^2 (y_1 - y_2)^2 \right). \tag{6.14}$$

By (6.14), for $z = (z_1, z_2) \in \mathbb{R}^2 \setminus \{\theta\}$ and $x = (x_1, x_2)$, $y = (y_1, y_2) \in \mathbb{R}^2$, we obtain

$$\|x\| = \frac{\sqrt{(y_1 z_1^2 + y_2 z_2^2)^2 + 4z_1^2 z_2^2 (y_1 - y_2)^2}}{z_1^2 + z_2^2}, \text{ if } x = \widehat{D}^* f(z)(y).$$

In particular, if we let $y_1 = y_2 = \pm \frac{1}{\sqrt{2}}$ in (6.14), we obtain

$$\frac{1}{(z_1^2 + z_2^2)^2}\left(\left(\frac{1}{\sqrt{2}} z_1^2 + \frac{1}{\sqrt{2}} z_2^2\right)^2 + 4z_1^2 z_2^2 \left(\frac{1}{\sqrt{2}} - \frac{1}{\sqrt{2}}\right)^2\right) = \frac{1}{2}. \tag{6.15}$$

For $\bar{z} = (\bar{z}_1, \bar{z}_2) \neq \theta$ with $\bar{w} = f(\bar{z})$, by (6.14) and (6.15), we have that

$$\hat{\alpha}(f, \bar{z}, \bar{w}) = \sup_{\eta > 0} \inf\{\|x\|: x \in \widehat{D}^* f(z)(y), z \in \mathbb{B}(\bar{z}, \eta) \setminus \{\theta\}, f(z) \in \mathbb{B}(\bar{w}, \eta), \|y\| = 1\}$$

$$= \sup_{\eta > 0} \inf\left\{\frac{\sqrt{(y_1 z_1^2 + y_2 z_2^2)^2 + 4z_1^2 z_2^2 (y_1 - y_2)^2}}{z_1^2 + z_2^2}: z \in \mathbb{B}(\bar{z}, \eta) \setminus \{\theta\}, f(z) \in \mathbb{B}(\bar{w}, \eta), \|y\| = 1\right\}$$

$$\leq \sup_{\eta > 0} \inf\left\{\frac{\sqrt{(y_1 z_1^2 + y_2 z_2^2)^2 + 4z_1^2 z_2^2 (y_1 - y_2)^2}}{z_1^2 + z_2^2}: z \in \mathbb{B}(\bar{z}, \eta) \setminus \{\theta\}, f(z) \in \mathbb{B}(\bar{w}, \eta), \|y\| = 1, y_1 = y_2 = \pm \frac{1}{\sqrt{2}}\right\}$$

$$= \sup_{\eta > 0} \inf\left\{\frac{1}{\sqrt{2}}: z \in \mathbb{B}(\bar{z}, \eta) \setminus \{\theta\}, f(z) \in \mathbb{B}(\bar{w}, \eta), \|y\| = 1, y_1 = y_2 = \pm \frac{1}{\sqrt{2}}\right\}$$

$$= \frac{1}{\sqrt{2}}.$$

This proves (i) in Part (b). Next, we suppose $\bar{z}_1 \bar{z}_2 = 0$ with $\bar{z}_i = 0$, for some $k = 1, 2$. By assumption that $\bar{z} = (\bar{z}_1, \bar{z}_2) \neq \theta$, we have $\bar{z}_{3-k} \neq 0$. By the continuity of $f$ around $\bar{z} = (\bar{z}_1, \bar{z}_2) \neq \theta$, for any $\eta > 0$, we can have $z = (z_1, z_2) \in \mathbb{B}(\bar{z}, \eta) \setminus \{\theta\}$ satisfying $f(z) \in \mathbb{B}(\bar{w}, \eta)$ and $z_k = 0, z_{3-k} \neq 0$.

$$\hat{\alpha}(f, \bar{z}, \bar{w}) = \sup_{\eta > 0} \inf\{\|x\|: x \in \widehat{D}^* f(z)(y), z \in \mathbb{B}(\bar{z}, \eta) \setminus \{\theta\}, f(z) \in \mathbb{B}(\bar{w}, \eta), \|y\| = 1\}$$

$$= \sup_{\eta > 0} \inf\left\{\frac{\sqrt{(y_1 z_1^2 + y_2 z_2^2)^2 + 4z_1^2 z_2^2 (y_1 - y_2)^2}}{z_1^2 + z_2^2}: z \in \mathbb{B}(\bar{z}, \eta) \setminus \{\theta\}, f(z) \in \mathbb{B}(\bar{w}, \eta), \|y\| = 1\right\}$$

$$\leq \sup_{\eta > 0} \inf\left\{\frac{\sqrt{(y_1 z_1^2 + y_2 z_2^2)^2 + 4z_1^2 z_2^2 (y_1 - y_2)^2}}{z_1^2 + z_2^2}: z \in \mathbb{B}(\bar{z}, \eta) \setminus \{\theta\}, f(z) \in \mathbb{B}(\bar{w}, \eta), z_k = 0, z_{3-k} \neq 0, \|y\| = 1\right\}$$

$$\leq \sup_{\eta > 0} \inf\left\{\frac{|y_k 0 + 0 z_{3-k}^2|}{z_1^2 + z_2^2}: z \in \mathbb{B}(\bar{z}, \eta) \setminus \{\theta\}, f(z) \in \mathbb{B}(\bar{w}, \eta), z_k = 0, z_{3-k} \neq 0. \|y\| = 1, y_k = 1, y_{3-k} = 0\right\}$$

$$\leq \sup_{\eta > 0} \inf\{0: z \in \mathbb{B}(\bar{z}, \eta) \setminus \{\theta\}, f(z) \in \mathbb{B}(\bar{w}, \eta), z_k = 0, z_{3-k} \neq 0. \|y\| = 1, y_k = 1, y_{3-k} = 0\}$$

$$= 0.$$

This proves (ii) in Part (b). Next, we prove (iii) in (b). By (4.12), we estimate

$$\hat{\alpha}(f,\bar{z},\bar{w}) = \sup_{\eta>0}\inf\{\|x\|: x \in \widehat{D}^*f(z)(y), z \in \mathbb{B}(\bar{z},\eta)\setminus\{\theta\}, f(z) \in \mathbb{B}(\bar{w},\eta), \|y\| = 1\}$$

$$= \liminf_{\eta\downarrow 0}\left\{\frac{\sqrt{(y_1z_1^2+y_2z_2^2)^2+4z_1^2z_2^2(y_1-y_2)^2}}{z_1^2+z_2^2}: z \in \mathbb{B}(\bar{z},\eta)\setminus\{\theta\}, f(z) \in \mathbb{B}(\bar{w},\eta), \|y\| = 1\right\}$$

$$\leq \liminf_{\eta\downarrow 0}\left\{\frac{\sqrt{(y_1z_1^2+y_2z_2^2)^2+4z_1^2z_2^2(y_1-y_2)^2}}{z_1^2+z_2^2}: z \in \mathbb{B}(\bar{z},\eta)\setminus\{\theta\}, f(z) \in \mathbb{B}(\bar{w},\eta), y_1 = -\frac{z_2^2}{\sqrt{z_2^4+z_2^4}}, y_2 = \frac{z_1^2}{\sqrt{z_2^4+z_2^4}}\right\}$$

$$\leq \liminf_{\eta\downarrow 0}\left\{\frac{\sqrt{4z_1^2z_2^2}}{\sqrt{z_2^4+z_2^4}}: z \in \mathbb{B}(\bar{z},\eta)\setminus\{\theta\}, f(z) \in \mathbb{B}(\bar{w},\eta), y_1 = -\frac{z_2^2}{\sqrt{z_2^4+z_2^4}}, y_2 = \frac{z_1^2}{\sqrt{z_2^4+z_2^4}}\right\}$$

$$\leq \liminf_{\eta\downarrow 0}\left\{\frac{\sqrt{4\bar{z}_1^2\bar{z}_2^2}}{\sqrt{\bar{z}_2^4+\bar{z}_2^4}}: \bar{z} \in \mathbb{B}(\bar{z},\eta)\setminus\{\theta\}, f(\bar{z}) \in \mathbb{B}(\bar{w},\eta), y_1 = -\frac{\bar{z}_2^2}{\sqrt{\bar{z}_2^4+\bar{z}_2^4}}, y_2 = \frac{\bar{z}_1^2}{\sqrt{\bar{z}_2^4+\bar{z}_2^4}}\right\}$$

$$= \frac{2|\bar{z}_1\bar{z}_2|}{\sqrt{\bar{z}_2^4+\bar{z}_2^4}}.$$

This proves (iii) in (b). □

## 7. An Upper Bound of Covering Constants for Mappings from $\mathbb{R}^n$ to $\mathbb{R}^m$ with $n \geq m$

From the two previous sections, we see that it is difficult to precisely calculate the covering constants for single-valued mappings from $\mathbb{R}^n$ to $\mathbb{R}^m$ with $n \geq m$. In this section, we find an upper bound of covering constants.

**Theorem 7.1**. *Let $n \geq m$. Let $f = (f_1, f_2, \ldots, f_m): \mathbb{R}^n \to \mathbb{R}^m$ be a single-valued mapping. Let $\bar{z} = (\bar{z}_1, \bar{z}_2, \ldots, \bar{z}_n) \in \mathbb{R}^n$ with $\bar{w} = f(\bar{z}) \in \mathbb{R}^m$. Suppose that $f$ is Fréchet differentiable around $\bar{z}$ such that the Fréchet derivative and Mordukhovich derivative of $f$ at $z$ are the following $n\times m$ and $m\times n$ matrices, respectively*

$$\nabla f(z) = \begin{pmatrix} \frac{\partial f_1}{\partial x_1}(z_1,z_2,\ldots,z_n) & \cdots & \frac{\partial f_m}{\partial x_1}(z_1,z_2,\ldots,z_n) \\ \vdots & \ddots & \vdots \\ \frac{\partial f_1}{\partial x_n}(z_1,z_2,\ldots,z_n) & \cdots & \frac{\partial f_m}{\partial x_n}(z_1,z_2,\ldots,z_n) \end{pmatrix},$$

and 
$$\widehat{D}^*f(z) = \nabla f(z)^T = \begin{pmatrix} \frac{\partial f_1}{\partial x_1}(z_1,z_2,\ldots,z_n) & \cdots & \frac{\partial f_1}{\partial x_n}(z_1,z_2,\ldots,z_n) \\ \vdots & \ddots & \vdots \\ \frac{\partial f_m}{\partial x_1}(z_1,z_2,\ldots,z_n) & \cdots & \frac{\partial f_m}{\partial x_n}(z_1,z_2,\ldots,z_n) \end{pmatrix}. \quad (7.1)$$

*Then,*

$$\hat{\alpha}(f,\bar{z},\bar{w}) \leq \sqrt{\sum_{i=1}^n\sum_{j=1}^m\left(\frac{\partial f_j}{\partial x_i}(\bar{z}_1,\bar{z}_2,\ldots,\bar{z}_n)\right)^2}.$$

*Proof.* By Theorem 1.38 in [18] and the condition in this proposition, we have that $f$ is also Mordukhovich differentiable at point $z$ and the Mordukhovich derivative of $f$ at $z$ is represented by (7.1). Let

$$\|\nabla f(z)\|_{n\times m} = \|\widehat{D}^*f(z)\|_{n\times m} = \sqrt{\sum_{i=1}^{n}\sum_{j=1}^{m}\left(\frac{\partial f_j}{\partial x_i}(\bar{z}_1,\bar{z}_2,\ldots,\bar{z}_n)\right)^2}$$

Let $z = (z_1, z_2, \ldots, z_n) \in \mathbb{R}^n$. Let $v = (v_1, v_2, \ldots, v_n) \in \mathbb{R}^n$ and $y = (y_1, y_2, \ldots, y_m) \in \mathbb{S}^m$. If $\widehat{D}^*f(z) = x$, then

$$(v_1, v_2, \ldots, v_n) = (y_1, y_2, \ldots, y_m) \begin{pmatrix} \frac{\partial f_1}{\partial x_1}(z_1, z_2, \ldots, z_n) & \cdots & \frac{\partial f_1}{\partial x_n}(z_1, z_2, \ldots, z_n) \\ \vdots & \ddots & \vdots \\ \frac{\partial f_m}{\partial x_1}(z_1, z_2, \ldots, z_n) & \cdots & \frac{\partial f_m}{\partial x_n}(z_1, z_2, \ldots, z_n) \end{pmatrix}.$$

This implies that

$$v_i = \sum_{j=1}^{m} y_j \frac{\partial f_j}{\partial x_i}(z_1, z_2, \ldots, z_n), \text{ for } i = 1, 2, \ldots, n.$$

Then, by $y = (y_1, y_2, \ldots, y_m) \in \mathbb{S}^m$, we have

$$\|v\|_n^2 = \sum_{i=1}^{n}\left(\sum_{j=1}^{m} y_j \frac{\partial f_j}{\partial x_i}(z_1, z_2, \ldots, z_n)\right)^2$$

$$= \sum_{i=1}^{n}\left(\sum_{j=1}^{m} y_j \frac{\partial f_j}{\partial x_i}(z_1, z_2, \ldots, z_n)\right)^2$$

$$\leq \sum_{i=1}^{n}\left(\sum_{j=1}^{m} y_j^2 \sum_{j=1}^{m}\left(\frac{\partial f_j}{\partial x_i}(z_1, z_2, \ldots, z_n)\right)^2\right)$$

$$= \sum_{1\leq i\leq n, 1\leq j\leq m}\left(\frac{\partial f_j}{\partial x_i}(z_1, z_2, \ldots, z_n)\right)^2.$$

Since $\bar{z} \in \mathbb{B}(\bar{z}, \eta), f(\bar{z}) = \bar{w} \in \mathbb{B}(\bar{w}, \eta)$, for any $\eta > 0$, this implies that

$$\hat{\alpha}(f, \bar{z}, \bar{w}) = \sup_{\eta>0}\inf\{\|x\|_n : v \in \widehat{D}^*f(z)(y), z \in \mathbb{B}(\bar{z}, \eta), f(z) \in \mathbb{B}(\bar{w}, \eta), \|y\|_m = 1\}$$

$$\leq \sup_{\eta>0}\inf\left\{\sqrt{\sum_{1\leq i\leq n, 1\leq j\leq m}\left(\frac{\partial f_j}{\partial x_i}(z_1, z_2, \ldots, z_n)\right)^2} : v \in \widehat{D}^*f(z)(y), z \in \mathbb{B}(\bar{z}, \eta), f(z) \in \mathbb{B}(\bar{w}, \eta), \|y\|_m = 1\right\}$$

$$\leq \sup_{\eta>0}\left\{\sqrt{\sum_{1\leq i\leq n, 1\leq j\leq m}\left(\frac{\partial f_j}{\partial x_i}(\bar{z}_1, \bar{z}_2, \ldots, \bar{z}_n)\right)^2} : z \in \mathbb{B}(\bar{z}, \eta), f(z) \in \mathbb{B}(\bar{w}, \eta)\right\}$$

$$= \sqrt{\sum_{i=1}^{n}\sum_{j=1}^{m}\left(\frac{\partial f_j}{\partial x_i}(\bar{z}_1, \bar{z}_2, \ldots, \bar{z}_n)\right)^2}. \qquad \square$$

## 8. Some Applications to Parameterized Coincidence Point Problems

Mordukhovich derivatives (coderivatives) of mappings (for both of set-valued and single-valued) have been widely applied to variational analysis. One of the important applications of the Mordukhovich derivatives is to define and to calculate the covering constants for both set-valued and single-valued mappings in Banach spaces. By the concept of covering constants, in [1], the authors Arutyunov Mordukhovich and Zhukovskiy proved an important theorem, which is called Arutyunov Mordukhovich Zhukovskiy Parameterized Coincidence Point Theorem (It is simply named as AMZ Theorem). The AMZ Theorem provides a very powerful tool in set-valued and variational analysis (see [1, 11−15]). The underlying spaces of the AMZ Theorem are Asplund spaces, which includes Euclidean spaces as special cases (See [1, 23] for more details). In this section, we use the results obtained in the previous section about the covering constants of mappings in $\mathbb{R}^2$ and by applying the AMZ Theorem to solve some parameterized equations in $\mathbb{R}^2$. We recall the AMZ Theorem.

**(AMZ Theorem)** *Let the Banach spaces X and Y in be Asplund and let P be a topological space. Let $F: X \rightrightarrows Y$ and $G(\cdot, \cdot): X \times P \rightrightarrows Y$ be set-valued mappings. Let $\bar{x} \in X$ and $\bar{y} \in Y$ with $\bar{y} \in F(\bar{x})$. Suppose that the following conditions are satisfied:*

(A1) *The multifunction $F: X \rightrightarrows Y$ is closed around $(\bar{x}, \bar{y})$.*

(A2) *There are neighborhoods $U \subset X$ of $\bar{x}$, $V \subset Y$ of $\bar{y}$, and $O$ of $\bar{p} \in P$ as well as a number $\beta \geq 0$ such that the multifunction $G(\cdot, p): X \rightrightarrows Y$ is Lipschitz-like on U relative to V for each $p \in O$ with the uniform modulus $\beta$, while the multifunction $p \to G(\bar{x}, p)$ is lower/inner semicontinuous at $\bar{p}$.*

(A3) *The Lipschitzian modulus $\beta$ of $G(\cdot, p)$ is chosen as $\beta < \hat{\alpha}(F, \bar{x}, \bar{y})$, where $\hat{\alpha}(F, \bar{x}, \bar{y})$ is the covering constant of F around $(\bar{x}, \bar{y})$ taken from (2.5).*

*Then for each $\alpha > 0$ with $\beta < \alpha < \hat{\alpha}(F, \bar{x}, \bar{y})$, there exist a neighborhood $W \subset P$ of $\bar{p}$ and a single-valued mapping $\sigma: W \to X$ such that whenever $p \in W$ we have*

$$F(\sigma(p)) \cap G(\sigma(p), p) \neq \emptyset \quad \text{and} \quad \|\sigma(p) - \bar{x}\|_X \leq \frac{\text{dist}(\bar{y}, G(\bar{x}, p))}{\alpha - \beta}.$$

In particular, if both *F* and *G* are single-valued mappings, then we get the following corollary of the AMZ Theorem.

**Corollary 2.2 in [12].** *Let the Banach spaces X and Y be Asplund and let P be a topological space. Let $F: X \to Y$ and $G(\cdot, \cdot): X \times P \to Y$ be single-valued mappings. Let $\bar{x} \in X$ and $\bar{y} \in Y$ with $\bar{y} = F(\bar{x})$. Suppose that the following conditions are satisfied:*

(A1) *The mapping $F: X \to Y$ is continuous around $(\bar{x}, \bar{y})$.*

(A2) *There are neighborhoods $U \subset X$ of $\bar{x}$, $V \subset Y$ of $\bar{y}$, and $O$ of $\bar{p} \in P$ as well as a number $\beta \geq 0$ such that the mapping $G(\cdot, p): X \to Y$ satisfies the Lipschitz condition on U relative to V for each $p \in O$ with the uniform modulus $\beta$, while the mapping $p \to G(\bar{x}, p)$ is lower semicontinuous at $\bar{p}$.*

(A3) *The Lipschitzian modulus $\beta$ of $G(\cdot, p)$ is chosen as $\beta < \hat{\alpha}(F, \bar{x}, \bar{y})$, where $\hat{\alpha}(F, \bar{x}, \bar{y})$ is the covering constant of F around $(\bar{x}, \bar{y})$ taken from (2.5).*

*Then for each $\alpha > 0$ with $\beta < \alpha < \hat{\alpha}(F, \bar{x}, \bar{y})$, there exist a neighborhood $W \subset P$ of $\bar{p}$ and a single-valued mapping $\sigma: W \to X$ such that whenever $p \in W$ we have*

$$F(\sigma(p)) = G(\sigma(p), p) \quad \text{and} \quad \|\sigma(p) - \bar{x}\|_X \leq \frac{\|G(\bar{x}, p) - \bar{y}\|_Y}{\alpha - \beta}.$$

**Theorem 8.1.** *Let $C$ be a topological space. Let $h(\cdot, \cdot) = (h_1(\cdot, \cdot), h_2(\cdot, \cdot)): \mathbb{R}^2 \times C \to \mathbb{R}^2$ be a single-valued mapping. Let $\omega = (\omega_1, \omega_2): C \to \mathbb{R}^2$ be a single-valued lower semicontinuous mapping. Let $f: \mathbb{R}^2 \to \mathbb{R}^2$ be a single-valued mapping defined by*

$$f(x_1, x_2) = (x_1^2 - x_2^2, 2x_1 x_2), \text{ for any } (x_1, x_2) \in \mathbb{R}^2.$$

*Let $\bar{x}$ and $\bar{y} \in \mathbb{R}^2\{\theta\}$ with $\bar{y} = f(\bar{x})$. Let $\bar{s} \in C$. Let $\mathbb{B}(\bar{x}, \|\bar{x}\|)$ be the closed ball in $\mathbb{R}^2$ with radius $\|\bar{x}\|$ and centered at $\bar{x}$. Suppose that the following conditions are satisfied:*

*There is a number $\beta \in [0, \|\bar{x}\|)$ such that the mapping $h(\cdot, s): \mathbb{R}^2 \to \mathbb{R}^2$ satisfies the Lipschitz condition on $\mathbb{B}\left(\bar{x}, \frac{\|\bar{x}\|}{2}\right)$ for each $s \in C$ with the uniform modulus $\beta$, while the mapping $s \to h(\bar{x}, s)$ is lower semicontinuous at $\bar{s}$.*

*Then for each $\alpha > 0$ with $\beta < \alpha < \|\bar{x}\|$, there exist a neighborhood $W \subset C$ of $\bar{s}$ and a single-valued mapping $\sigma = (\sigma_1, \sigma_2): W \to \mathbb{R}^2$ such that whenever $s \in W$ we have*

$$\begin{cases} \sigma_1(s)^2 - \sigma_2(s)^2 = h_1\left((\sigma_1(s), \sigma_2(s)), s\right) + \omega_1(s) \\ 2\sigma_1(s)\sigma_2(s) = h_2\left((\sigma_1(s), \sigma_2(s)), s\right) + \omega_2(s) \end{cases},$$

*and*

$$\|(\sigma_1(s), \sigma_2(s)) - (\bar{x}_1, \bar{x}_2)\| \leq \frac{\|h((\bar{x}_1, \bar{x}_2), s) + (\omega_1(s), \omega_2(s)) - (\bar{x}_1^2 - \bar{x}_2^2, 2\bar{x}_1\bar{x}_2)\|}{\alpha - \beta}. \tag{8.1}$$

*In particular, $\sigma$ satisfies that, whenever $s \in W$,*

$$\left(h_1\left((\sigma_1(s), \sigma_2(s)), s\right) + \omega_1(s)\right)^2 + \left(h_2\left((\sigma_1(s), \sigma_2(s)), s\right) + \omega_2(s)\right)^2 = (\sigma_1(s)^2 + \sigma_2(s)^2)^2. \tag{8.2}$$

*Proof.* The given polynormal mapping $f$ in this theorem is studied in mapping 6.7, from which we have

$$\hat{\alpha}(f, x, f(x)) = 2\|x\|, \text{ for any } x \in \mathbb{R}^2. \tag{8.3}$$

Notice that, for any $x \in \mathbb{B}\left(\bar{x}, \frac{\|\bar{x}\|}{2}\right)$, we have that $\|x\| \geq \frac{\|\bar{x}\|}{2}$. By (8.3), we have that

$$\hat{\alpha}(f, x, f(x)) = 2\|x\| \geq \|\bar{x}\|, \text{ for any } x \in \mathbb{B}\left(\bar{x}, \frac{\|\bar{x}\|}{2}\right). \tag{8.4}$$

In AMZ Theorem, we particularly take

$$F(x) = f(x) \text{ and } G(x, s) = h(x, s) + \omega(s), \text{ for all } x \in \mathbb{R}^2, s \in C.$$

As single-valued mappings, we can verify that the mappings, which are defined in (8.4), $F: \mathbb{R}^2 \to \mathbb{R}^2$ and $G(\cdot, \cdot): \mathbb{R}^2 \times C \to \mathbb{R}^2$ satisfy all conditions in the AMZ Theorem, from which the conclusions (8.1) are immediately proved. (8.2) is reduced by (8.1) immediately. □

**Theorem 8.2.** *Let $C$ be a topological space. Let $h(\cdot, \cdot) = (h_1(\cdot, \cdot), h_2(\cdot, \cdot)): \mathbb{R}^2 \times C \to \mathbb{R}^2$ be a single-valued mapping. Let $\omega = (\omega_1, \omega_2): C \to \mathbb{R}^2$ be a single-valued lower semicontinuous mapping. Let $f: \mathbb{R}^2 \to \mathbb{R}^2$ be the single-valued mapping defined in (5.1),*

$$f(x_1, x_2) = \left( \frac{x_1^2 - x_2^2}{\sqrt{x_1^2 + x_2^2}}, \frac{2x_1 x_2}{\sqrt{x_1^2 + x_2^2}} \right), \text{ for } (x_1, x_2) \in \mathbb{R}^2 \setminus \{\theta\} \text{ with } f(\theta) = \theta.$$

*Let $\bar{x}$ and $\bar{y} \in \mathbb{R}^2\{\theta\}$ with $\bar{y} = f(\bar{x})$. Let $\bar{s} \in C$. Suppose that the following conditions are satisfied:*

*There is a number $\beta \in [0, 1)$ such that the mapping $h(\cdot, s): \mathbb{R}^2 \to \mathbb{R}^2$ satisfies the Lipschitz condition on $\mathbb{R}^2$ for each $s \in C$ with the uniform modulus $\beta$ and the mapping $s \to h(\bar{x}, s)$ is lower semicontinuous at $\bar{s}$.*

*Then for each $\alpha > 0$ with $\beta < \alpha < 1$, there exist a neighborhood $W \subset C$ of $\bar{s}$ and a single-valued mapping $\sigma = (\sigma_1, \sigma_2): W \to \mathbb{R}^2$ such that whenever $s \in W$ we have*

$$\begin{cases} \frac{\sigma_1(s)^2 - \sigma_2(s)^2}{\sqrt{\sigma_1(s)^2 + \sigma_2(s)^2}} = h_1\left((\sigma_1(s), \sigma_2(s)), s\right) + \omega_1(s) \\ \frac{2\sigma_1(s)\sigma_2(s)}{\sqrt{\sigma_1(s)^2 + \sigma_2(s)^2}} = h_2\left((\sigma_1(s), \sigma_2(s)), s\right) + \omega_2(s) \end{cases},$$

*and*
$$\|(\sigma_1(s), \sigma_2(s)) - (\bar{x}_1, \bar{x}_2)\| \leq \frac{\left\| h((\bar{x}_1, \bar{x}_2), s) + (\omega_1(s), \omega_2(s)) - \left( \frac{\bar{x}_1^2 - \bar{x}_2^2}{\sqrt{\bar{x}_1^2 + \bar{x}_2^2}}, \frac{2\bar{x}_1 \bar{x}_2}{\sqrt{\bar{x}_1^2 + \bar{x}_2^2}} \right) \right\|}{\alpha - \beta}. \tag{8.5}$$

*In particular, $\sigma$ satisfies that, whenever $s \in W$,*

$$\left( h_1\left((\sigma_1(s), \sigma_2(s)), s\right) + \omega_1(s) \right)^2 + \left( h_2\left((\sigma_1(s), \sigma_2(s)), s\right) + \omega_2(s) \right)^2 = \sigma_1(s)^2 + \sigma_2(s)^2. \tag{8.6}$$

*Proof.* The given polynormal mapping $f$ is defined by (5.1) and is studied in Section 5. From Theorem 5.5, we have

$$\hat{\alpha}(f, x, f(x)) = 1, \text{ for any } x \in \mathbb{R}^2. \tag{8.7}$$

By (8.7) and the conditions in this theorem, we have that

$$\hat{\alpha}(f, x, f(x)) = 1 > \beta \geq 0, \text{ for any } x \in \mathbb{R}^2.$$

In AMZ Theorem, we particularly take

$$F(x) = f(x) \text{ and } G(x, s) = h(x, s) + \omega(s), \text{ for all } x \in \mathbb{R}^2, s \in C.$$

As single-valued mappings, we can verify that the mappings, which are defined in (8.4), $F: \mathbb{R}^2 \to \mathbb{R}^2$ and $G(\cdot, \cdot): \mathbb{R}^2 \times C \to \mathbb{R}^2$ satisfy all conditions in the AMZ Theorem, from which the conclusions (8.5) are immediately proved. Since $f$ is norm preserving, it implies that

$$\left( \frac{\sigma_1(s)^2 - \sigma_2(s)^2}{\sqrt{\sigma_1(s)^2 + \sigma_2(s)^2}} \right)^2 + \left( \frac{2\sigma_1(s)\sigma_2(s)}{\sqrt{\sigma_1(s)^2 + \sigma_2(s)^2}} \right)^2 = \sigma_1(s)^2 + \sigma_2(s)^2. \tag{8.8}$$

Then, (8.6) is reduced by (8.5) and (8.8) immediately. □

## 9. Conclusion and Remarks

In this paper, we prove the guidelines for calculating the Fréchet derivatives of single-valued mappings in

Euclidean spaces, which is represented by the linear approximation. With great help of the relationship between Fréchet derivatives and Mordukhovich derivatives (It is given by Theorem 1.38 in [18]), we derive the algorithm for calculating the Mordukhovich derivatives of single-valued mappings in $\mathbb{R}^2$.

For the purpose to find some formulas of Fréchet derivatives and Mordukhovich derivatives of single-valued mappings in Euclidean spaces, we provide examples of polynomial mappings, rational mappings, exponential mappings and logarithm mappings in $\mathbb{R}^2$. which includes. For each given mapping, we find it's Fréchet derivatives, Mordukhovich derivatives, and covering constants.

From the examples studied in this paper, we see the difficulty and complexity for finding the covering constants of single-valued mappings in Euclidean spaces. It is easy to understand that will be much more difficult and more complicated to find the covering constants of single-valued mappings in Banach spaces. If we consider Hilbert spaces to be special cases of Banach spaces, we have the following problems for consideration by interested readers.

(i) Find the principles for the existence of the Fréchet derivatives of single-valued mappings in Hilbert spaces;
(ii) Find the guidelines for calculating the Mordukhovich derivatives of single-valued mappings in Hilbert spaces;
(iii) Find an applicable algorithm to calculate the covering constants for single-valued mappings in Hilbert spaces;
(iv) Consider the above problems for set-valued mappings in Hilbert spaces;
(v) Consider the above problems for single valued and set-valued mappings in Asplund spaces.